# LIMIT OF NORMALIZED QUADRANGULATIONS: THE BROWNIAN MAP

By Jean-François Marckert and Abdelkader Mokkadem

*Université de Versailles Saint-Quentin*

Consider $q_n$ a random pointed quadrangulation chosen equally likely among the pointed quadrangulations with $n$ faces. In this paper we show that, when $n$ goes to $+\infty$, $q_n$ suitably normalized converges weakly in a certain sense to a random limit object, which is continuous and compact, and that we name the *Brownian map*. The same result is shown for a model of rooted quadrangulations and for some models of rooted quadrangulations with random edge lengths. A metric space of rooted (resp. pointed) abstract maps that contains the model of discrete rooted (resp. pointed) quadrangulations and the model of the Brownian map is defined. The weak convergences hold in these metric spaces.

**1. Introduction.** A planar map is a proper embedding without edge crossing of a connected graph in the sphere. Two planar maps are identical if one of them can be mapped to the other by a homeomorphism that preserves the orientation of the sphere. A planar map is a quadrangulation if all faces have degree four. A quadrangulation is bipartite, does not contain any loop, but may contain some multiple edges. Any quadrangulation with $n$ faces has $2n$ edges and $n+2$ vertices. Notice that there is a difference between planar maps and planar graphs since a planar graph can have several nonidentical representations on the sphere as planar map. On Figure 1, one finds two representations of the same planar graph on the sphere. The right one is a quadrangulation (by convention, an edge that lies entirely in a face is counted twice in the degree of the face).

A planar map is said to be pointed (resp. rooted) if one node, called the origin or the root-vertex (resp. one oriented edge, called the root or root-edge) is distinguished. Two pointed (resp. rooted) quadrangulations are identical if the homeomorphism preserves also the distinguished node (resp. oriented edge). We denote by $\mathcal{Q}_n^\bullet$ (resp. $\overrightarrow{\mathcal{Q}}_n$) the set of pointed (resp. rooted) quadrangulations with $n$ faces.

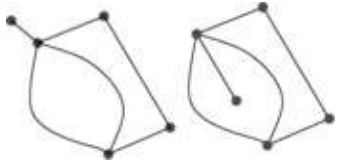

Fig. 1.　*Two different maps.*

Since the pioneer work of Tutte [34], the combinatorial study of planar maps has received considerable attention. Many statistical properties have been obtained [9, 10, 11, 20, 32, 34, 35] for a number of classes of finite planar maps. Among the classes of planar maps, the best known is the class of rooted finite planar trees; Aldous [1, 2] built a mathematical object called a "continuum random tree," which is the limit of random rooted finite planar trees under an appropriate scaling. A question arises: does it exist a similar (continuous) limit object for some other classes of planar maps? This question is important not only in combinatorics and in probability but also in theoretical physics. As a matter of fact, it has been realized in these last years that random planar structures have a leading role in quantum field theory, string theory and quantum gravity [6, 7, 12, 17, 37]. Following the algebraic topology point of view, the physicists consider triangulations and quadrangulations (or other classes of maps) as discretized versions of 2-dimensional manifolds; they are mainly interested in a continuous limit for suitably normalized discretization. A limit behavior without any scaling has been investigated by Angel and Schramm [8] see also Krikun [24]. They show that the uniform law on the set of finite planar triangulations with $n$ faces converges to a law on the set of infinite planar triangulations (endowed with a non-Archimedean metric). They obtain a limit behavior of the triangulations in the ball of fixed radius $k$ around the origin. Chassaing and Durhuus [14] show a similar result for the convergence of unscaled random quadrangulations with a different approach.

With the topology used in [8, 14], the limit random metric space can hardly be continuous and bounded. Physicists [7, 37] give consistent arguments to show that if such a limit object exists, then the scaling should be $n^{1/4}$. Results in this direction have been obtained by Chassaing and Schaeffer [15]; in particular, they show that the radius of a random rooted quadrangulation taken uniformly in $\overrightarrow{\mathcal{Q}}_n$ and scaled by $n^{1/4}$ converges in distribution, up to a multiplicative constant, to the range of the Brownian snake.

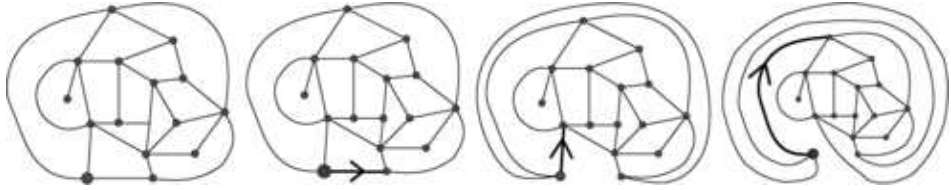

FIG. 2. *A pointed quadrangulation from $\mathcal{Q}_{14}^\bullet$ and the canonical representations of the three rooted quadrangulations from $\vec{\mathcal{Q}}_{14}$ in the fiber by $K$.*

Our purpose in the present paper is to show that suitably scaled random quadrangulations, uniformly chosen in $\mathcal{Q}_n^\bullet$ or chosen in $\vec{\mathcal{Q}}_n$ endowed with the distribution $\mathbb{P}_D^n$ defined below, converge to a limit object, "a continuum random map"; we name this object the Brownian map. As expected by the physicists, the adequate scaling is $n^{1/4}$.

*Models.* We are mainly interested in the limit of two random models of quadrangulations:

- $(\mathcal{Q}_n^\bullet, \mathbb{P}_U^n)$, where $\mathbb{P}_U^n$ is the uniform distribution on $\mathcal{Q}_n^\bullet$.
- $(\vec{\mathcal{Q}}_n, \mathbb{P}_D^n)$, where $\mathbb{P}_D^n$ is defined for each $q \in \vec{\mathcal{Q}}_n$ with root degree $\deg(q)$ by

$$\mathbb{P}_D^n(q) = \frac{c_n}{\deg(q)}, \qquad \text{where } c_n = \left( \sum_{q' \in \vec{\mathcal{Q}}_n} \frac{1}{\deg(q')} \right)^{-1}.$$

The probability $\mathbb{P}_D^n$ gives to each rooted quadrangulation a weight proportional to the inverse of its root degree (it is not the uniform distribution on $\vec{\mathcal{Q}}_n$ which is the law studied in [15]).

Denote by $K$ the canonical surjection $K$ from $\vec{\mathcal{Q}}_n$ onto $\mathcal{Q}_n^\bullet$ (see Figure 2). For any $q'$ in $\vec{\mathcal{Q}}_n$ with root-edge $\overrightarrow{vw}$, the pointed quadrangulation $K(q')$ is the planar pointed map whose origin is $v$ and which is identical to $q'$ as unrooted map. It will be shown in Lemma 4.19 that the distance in variation between the image of $\mathbb{P}_D^n$ by $K$ and the uniform distribution $\mathbb{P}_U^n$ on $\mathcal{Q}_n^\bullet$ goes to 0.

*Contents.* In Section 2 we gather some elements concerning rooted quadrangulations. The Schaeffer's bijection $Q$ between $\vec{\mathcal{Q}}_n$ and $\mathcal{W}_n^+$, the set of well-labeled trees with $n$ edges, is presented in Section 2.2. In Section 2.3 we describe the application $Q$ in a slightly different way. We exhibit two trees, the doddering tree $\mathcal{D}_n$ and the gluer tree $\mathcal{G}_n$, naturally associated with rooted quadrangulations (Section 2.4). This leads us to a new description of $Q$: a rooted quadrangulation is shown to be "$\mathcal{D}_n$ folded around $\mathcal{G}_n$," in other words, a rooted quadrangulation is shown to be $\mathcal{D}_n$ together with an identification of its nodes, with the help of $\mathcal{G}_n$.

This is the starting point of the notion of the rooted abstract map (Section 3). The leading idea is to construct a notion of maps sufficiently robust to be compatible with rooted quadrangulations described with the normalized version of $(\mathcal{D}_n, \mathcal{G}_n)$ and their limits, which are shown to exist. This leads us to present a notion of abstract trees in Section 3.2. An abstract tree is described in terms of a depth first walk and a measure. The convergence of $(\mathbf{D}_n, \mathbf{G}_n)$, normalized doddering tree and gluer tree under $\mathbb{P}_D^n$, is stated in Section 3.3. The notion of rooted abstract map is then presented in Section 3.4; then follows some elements on the topology and geometry of abstract maps. The convergence of normalized rooted quadrangulations under $\mathbb{P}_D^n$, presented as rooted abstract maps, is given in Section 3.5.1. The limit, that we name the Brownian map, is described with the help of the Brownian snake with lifetime process the normalized Brownian excursion. A model of rooted quadrangulation with random edge lengths is also shown to converge to the Brownian map (Section 3.5.2).

Using the surjection $K$, a pointed quadrangulation may be seen as an equivalence class of rooted quadrangulations. This is the point of view we use to build the notion of the pointed abstract map in Section 3.6. At the end of Section 3.6, the convergence of normalized pointed quadrangulations under $\mathbb{P}_U^n$ in the space of pointed abstract maps is given. The limit is still the Brownian map.

The remaining part of the paper is mainly devoted to the proofs. Each labeled tree $T$ with $n$ edges is encoded with the help of two discrete processes $R_n$ and $V_n$ deeply related with $(\mathcal{D}_n, \mathcal{G}_n)$; $V_n$ is the depth first walk of $T$ and $R_n$ is its labels process. Using the application $K$, each pointed quadrangulation is naturally encoded by a class of well-labeled trees (Section 4.2.1). A rerooting operation on well-labeled trees plays an important role: well-labeled trees associated with rooted quadrangulations in the same class modulo $K$ are equal up to a rerooting (Proposition 4.2). In order to prove the convergence of well-labeled trees under $\mathbb{P}_D^n$, or of classes of well-labeled trees modulo a rerooting, we adopt the following process: we introduce a family of labeled trees (well-labeled trees are "positive" labeled trees). We state the convergence of uniform normalized labeled trees (Proposition 4.8). We construct some classes of labeled trees corresponding to pointed quadrangulations, the classes of labeled trees being in bijection with the classes of well-labeled trees modulo rerooting (Theorem 4.5). The convergence of rescaled classes of labeled trees is stated in Theorem 4.10. It remains to deduce from this convergence, the convergence of the class of well-labeled trees (Proposition 4.17), and the convergence of rescaled well-labeled trees under $\mathbb{P}_D^n$; this is Theorem 3.3, proved in Section 4.7.

In Section 5 are shown the convergence of the radius and of the profile of rooted and pointed quadrangulation. In Section 6 a conclusion of the paper is given. The Appendix contains some postponed proofs.

**2. Combinatorics of rooted quadrangulations.** We begin with some considerations on discrete trees and on their encodings.

2.1. *Labeled trees and encoding of labeled trees.* A tree is a planar map with one face. A tree with a distinguished oriented edge $\overrightarrow{u_0 v}$ is called a rooted tree; $\overrightarrow{u_0 v}$ is the root-edge (or simply root) and $u_0$ is the root-vertex. The adjacent nodes of $u_0$ are called the children of $u_0$. If $u$ is a node different from $u_0$, and if $(u_0, u_1, \ldots, u_l, u)$ is the unique geodesic between $u$ and $u_0$, then the node $u_l$ is the father of $u$; the other adjacent nodes of $u$ are the children of $u$. The root-edge induces the notion of subtrees rooted at a node.

Around each node, there are two circular orders: the clockwise order and the reverse order. When a root $\overrightarrow{u_0 v}$ is given, the two circular orders around each node $u$ induce two corresponding total orders between the adjacent nodes (and edges) of $u$:

- if $u = u_0$, $(u_0, v)$ is the smallest incident edge of $u_0$,
- if $u \neq u_0$ and if $u_l$ is the father of $u$, the edge $(u, u_l)$ is the smallest incident edge of $u$.

We denote by $\Omega_n$ the set of rooted trees with $n$ edges. Its cardinality is $C_n = \binom{2n}{n}/(n+1)$, the $n$th Catalan number. A *labeled tree* with $n$ edges is a tree of $\Omega_n$ in which the $n+1$ nodes are labeled by integers that satisfy the two following conditions:

– the label of the root-vertex is 1,
– the difference between the labels of two adjacent nodes is 1, 0 or $-1$.

If all labels are positive, such a tree is called *well-labeled*. We denote by $\mathcal{W}_n$ (resp. $\mathcal{W}_n^+$) the set of labeled (resp. well-labeled) trees with $n$ edges. The well-labeled tree on Figure 3 belongs to $\mathcal{W}_{14}^+$. The unlabeled rooted tree that "supports" the labels is called the underlying tree.

The clockwise depth first traversal (CDFT) of the rooted tree $t \in \Omega_n$ is a function:

$$F : [\![0, 2n]\!] \longrightarrow \text{Nodes}(t) := \{\text{nodes of } t\},$$

which we regard as a walk around $t$. First, $F(0) = u_0$. For $i$ from 0 to $2n-1$, given $F(i) = z$, choose, if possible, and according to the clockwise order around $z$, the smallest child $w$ of $z$ which has not already been visited, and set $F(i+1) = w$. If not possible, let $F(i+1)$ be the parent of $z$ (see Figures 4 and 5 for illustrations).

The root-edge of $t$ is the oriented edge $\overrightarrow{F(0)F(1)} = \overrightarrow{u_0 u_1}$. The CDFT induces a total order on Nodes($t$), that we call the clockwise order (CO) (see the first picture of Figure 4). A *corner* is a sector between two consecutive edges around a vertex. For any node $u \neq u_0$, the successive times $i_1, \ldots, i_k$ such that $F(i_j) = u$ are in one-to-one correspondence with the corners around $u$. Thus $i_1, \ldots, i_k$ encode and order these corners ([1, 3, 5] and [7] encode the four corners around the node $F(1)$ in Figure 5). For the root-vertex, 0 and $2n$ encode the same corner. The clockwise depth first order defined here is also called in the literature the lexicographical order or the prefix order.

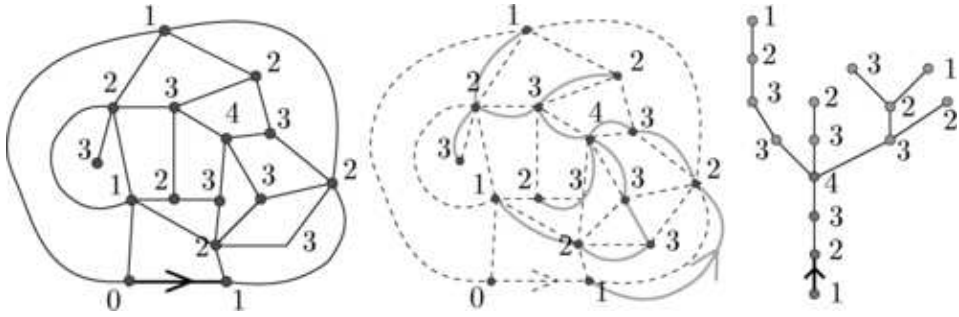

FIG. 3. *A rooted quadrangulation from $\overrightarrow{\mathcal{Q}}_{14}$ and the associated well-labeled tree belonging to $\mathcal{W}_{14}^+$.*

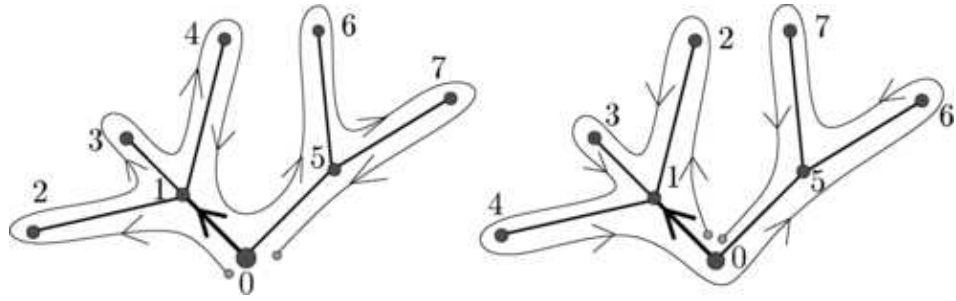

FIG. 4. *CDFT and next, the RDFT; the numbering of the nodes is done according to their first visit in the clockwise order, and in the second picture, in the reverse order.*

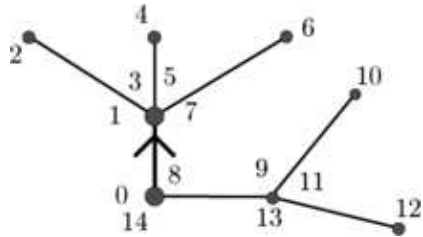

FIG. 5. *Clockwise depth first traversal and the notion of corner.*

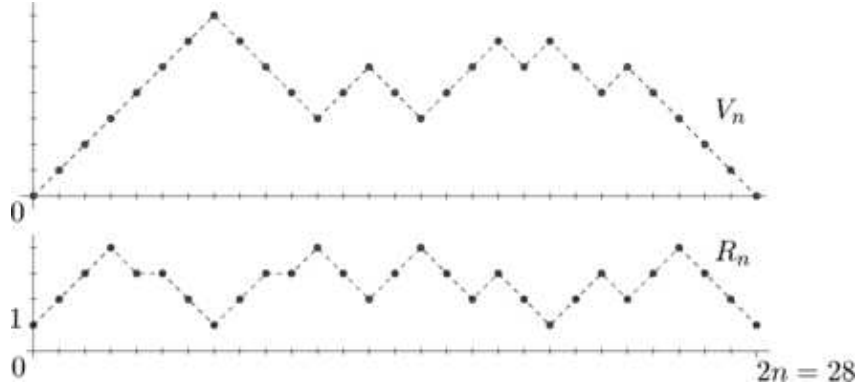

Fig. 6. *The processes $V_{14}$ and $R_{14}$ encoding the well-labeled tree of Figure 3.*

The clockwise depth first walk (CDFW) of $t \in \Omega_n$ is the process $V_n$:

$$(2.1) \qquad V_n(i) = d(F(0), F(i)), \qquad 0 \leq i \leq 2n,$$

where $d(u,v)$ is the number of edges in the unique shortest path between the nodes $u$ and $v$ (i.e., the graph distance between $u$ and $v$). The distance $V_n(i)$ is often called the height or the depth of the node $F(i)$. The process $V_n$ is also known as the Harris walk of $t$, or the tour of $t$ (see Figure 6).

Notice that $i$ and $j$ encode corners of the same node in $t$ iff

$$(2.2) \qquad \min\{V_n(u), u \in [\![i \vee j, i \wedge j]\!]\} = V_n(i) = V_n(j).$$

Let $T$ be an element of $\mathcal{W}_n$ with CDFT $F$ and CDFW $V_n$. The label process of $T$ is $(R_n(j))_{j \in [\![0,2n]\!]}$ defined by

$$(2.3) \qquad R_n(j) = \mathrm{label}(F(j)).$$

The bi-dimensional process $\{(R_n(k), V_n(k)), k \in [\![0, 2n]\!]\}$ uniquely determines $T$; we call it the encoding of $T$ (see Figure 6).

We define the reverse depth first traversal (RDFT) $\widetilde{F}$ of $t$ in the same way as the CDFT, except that the order used around each node is the reverse order; the total order induced on Nodes$(t)$ will be called the reverse order (RO) (see the second picture of Figure 4). The reverse depth first walk (RDFW) $\widetilde{V}_n$ is

$$(2.4) \qquad \widetilde{V}_n(i) = d(\widetilde{F}(0), \widetilde{F}(i)), \qquad 0 \leq i \leq 2n.$$

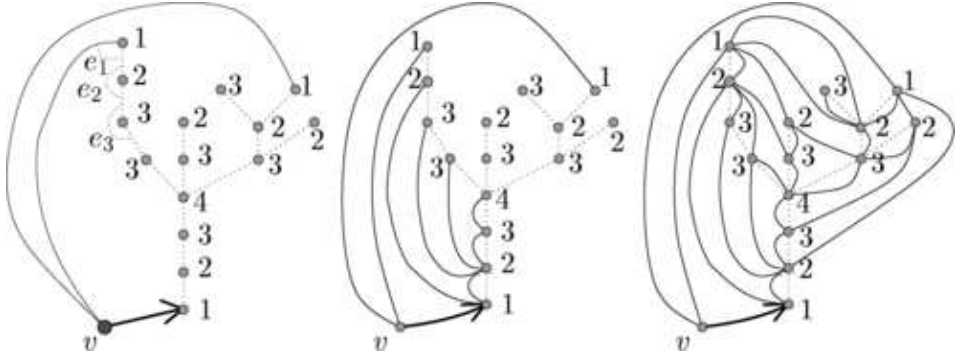

FIG. 7. *Construction of the rooted map associated with the well-labeled tree of Figure* 3.

2.2. *Schaeffer's bijection between* $\overrightarrow{\mathcal{Q}}_n$ *and* $\mathcal{W}_n^+$. The content of this Section 2.2 can be found in [33] or in [15].

THEOREM 2.1 ([16, 33]).　*There exists a bijection* $Q$ *from* $\mathcal{W}_n^+$ *onto* $\overrightarrow{\mathcal{Q}}_n$.

The bijection $Q$ was discovered by Cori and Vauquelin [16]. We present the construction of $Q$ given by Schaeffer [33].

*Description of* $Q$. Any $T$ in $\mathcal{W}_n^+$ has a unique face, the infinite face. A vertex of $T$ with degree $k$ defines $k$ corners and so the total number of corners is $2n$. We label each corner of $T$ with the label of its vertex. To build $Q(T)$, two steps are needed:

*First step* (see an illustration on Figure 7):

(a) Dot the edges of $T$.

(b) A vertex $v$, called the origin, with label 0 is placed in the infinite face.

(c) One edge is added between $v$ and each of the $l$ corners with label 1 (as on Figure 7).

(d) The root of $Q(T)$ is chosen as the added edge between $v$ and the first corner of the root-vertex of $T$.

After this first step one has a rooted map $T_0$ with $l$ faces. The next step takes place independently in each of these $l$ faces and is thus described for a generic face $F$ of $T_0$.

*Second step*: Let $k$ be the degree of $F$. Among the corners of $F$, only one is a corner of $v$ (and has label 0). Let the corners be numbered from 1 to $k$ in the clockwise order along the border, starting right after $v$. Let $e_i$ be the label of corner $i$ (one has $e_1 = e_{k-1} = 1$ and $e_k = 0$) (the numbering of the corners is started on a face in Figure 7). For the infinite face, the clockwise order is obtained by letting the infinite face on the right. The function successor $s$ is defined for all corners $2, \ldots, k-2$ by

(2.5) $$s(i) = \inf\{j > i | e_j = e_i - 1\}.$$

For each corner $i \in [\![2, k-2]\!]$, a chord $(i, s(i))$ is added inside the face $F$, in such a way that the various chords do not intersect. Remove the doted edges (the edges of $T$). The resulting map is a rooted quadrangulation $Q(T) \in \overrightarrow{\mathcal{Q}}_n$, with set of nodes $\{v\} \cup \mathrm{Nodes}(T)$.

REMARK 2.2. In [15, 33], the construction is a little bit different but equivalent: a chord $(i, s(i))$ is added only if $(i, s(i))$ is not already an edge of $T$. Then, not all the edges are removed but only the edges of $T$ that begin and end with the same label.

*Description of $Q^{-1}$.* Take a rooted quadrangulation $\omega$. Label the nodes of $\omega$ with their distances to the root-vertex. Consider a face $F$ and denote by $e_1, e_2, e_3, e_4$ the labels of the nodes of $F$ clockwise ordered. Two cases appear (up to a rotation of the indices):

- If $e_1 = e_3 = e_2 + 1 = e_4 + 1$, then add a blue edge between the nodes with label $e_1$ and $e_3$.
- If $e_4 = e_2 = e_1 + 1, e_3 = e_1 + 2$, then color the edge between $e_2$ and $e_3$ in blue.

The graph whose set of vertices is the set of vertices of the quadrangulation (minus its root-vertex) and whose edges are the blue edges turns out to be a tree [15]. The root of this blue tree is the first selected edge around the endpoint of the root of $\omega$. It is shown in [15] that the blue tree is $Q^{-1}(\omega)$. An example of this *tree-extraction* is made on the second picture of Figure 3.

2.3. *Construction of $Q(T)$ using the CDFT.* We give here a new presentation of the construction of $Q(T)$ more adapted to the present paper. We start with a definition of the predecessor function.

2.3.1. Let $N$ be a positive integer and let $\mathcal{R}$ be a process defined on $[\![0, N]\!]$ satisfying the conditions

$$\mathcal{R}(0) = 1,$$
(2.6) $$\mathcal{R}(j) \geq 1 \qquad \text{for } 1 \leq j \leq N,$$
$$\mathcal{R}(j+1) - \mathcal{R}(j) \in \{+1, 0, -1, -2, \ldots\} \qquad \text{for } j \in [\![0, N-1]\!].$$

The predecessor function $\mathcal{P}$ (we should write $\mathcal{P}_\mathcal{R}$) associated with $\mathcal{R}$ is defined for $i \in [\![0, N]\!]$ and takes its values in $[\![-1, N-1]\!]$: set $\mathcal{R}(-1) = 0$ and

(2.7) $$\mathcal{P}(i) = \max\{k \in [\![-1, i-1]\!], \mathcal{R}(k) = \mathcal{R}(i) - 1\}.$$

We say that $\mathcal{P}(i)$ is the predecessor of $i$. Thanks to (2.6), the predecessor function is well defined and has the following straightforward property.

LEMMA 2.3. *Let $i$ and $j$ be two integers such that $\mathcal{P}(i) < j < i$; then $\mathcal{P}(i) \leq \mathcal{P}(j) < j$. Thus, two cases arise: either $\mathcal{P}(j+1) = j$ or $\mathcal{P}(j+1) \leq \mathcal{P}(j)$.*

2.3.2. Now, we present a slight modification of Schaeffer's algorithm for the construction of $Q(T)$; the main point is that, thanks to this modification, the algorithm follows the CDFT of $T$.

Let $T$ be an element of $\mathcal{W}_n^+$ with CDFT $F$ and encoded by $(R_n^+, V_n^+)$. Since the labels of $T$ are positive, $R_n^+$ is a positive process satisfying (2.6) for $N = 2n-1$ (the point $2n$ is excluded in the construction). Let $P$ be the predecessor function associated with $R_n^+$.

(a) Dot the edges of $T$. Add a vertex $v$ in the unique face. Set $F(-1) = v$ and consider $-1$ as the single corner of $v$.

(b) Visit $T$ according to the CDFT from time $0$ to time $2n-1$. At time $i \in [\![0, 2n-1]\!]$, draw a chord $\widehat{(i, P(i))}$ starting from the corner $i$ and ending in the corner $P(i)$ such that:

(1) $\widehat{(i, P(i))}$ surrounds all the trajectory of the CFDT between $P(i)$ and $i$ (as drawn on Figure 8).

(2) $\widehat{(i, P(i))}$ surrounds all the chords $\widehat{(j, P(j))}$ such that $P(i) < j < i$. This is possible in virtue of Lemma 2.3.

The rooted planar map whose root-edge is the oriented chord $(v, F(0))$, whose edges are the chords $\widehat{(i, P(i))}, i \in [\![0, 2n-1]\!]$, and whose vertices are Nodes$(T) \cup \{v\}$, is exactly $Q(T)$. Indeed, let us examine why our construction builds the same quadrangulation as the one of Schaeffer. Let $i_1 < \cdots < i_k$ be the times such that $R_n^+(i_l) = 1$; for any $s$ such that $i_l < s < i_{l+1}$, the chord $\widehat{(i_{l+1}, v)}$ surrounds the chord $\widehat{(s, P(s))}$. The chord $\widehat{(s, P(s))}$ does not intersect the chord $\widehat{(i_l, v)}$, so that we can begin our construction by drawing all the chords $\widehat{(i_l, v)}, l = 1, \ldots, k$. We can continue the construction independently in each interval $]\!]i_l, i_{l+1}[\![$. The only difference with Schaeffer's procedure is that we work from $i_l + 1$ to $i_{l+1} - 1$, while Schaeffer works from $i_{l+1} - 1$ to $i_l + 1$ (see Figure 8 for an illustration).

2.4. *Construction of $Q(T)$ with a doddering tree and a gluer tree.* Our idea now is to perform the previous construction on the graphs of the encoding $(R_n^+, V_n^+)$. This is done in two steps. The first step is the construction of a rooted tree, *the doddering tree* $\mathcal{D}(R_n^+)$, containing once each edge of $Q(T)$.

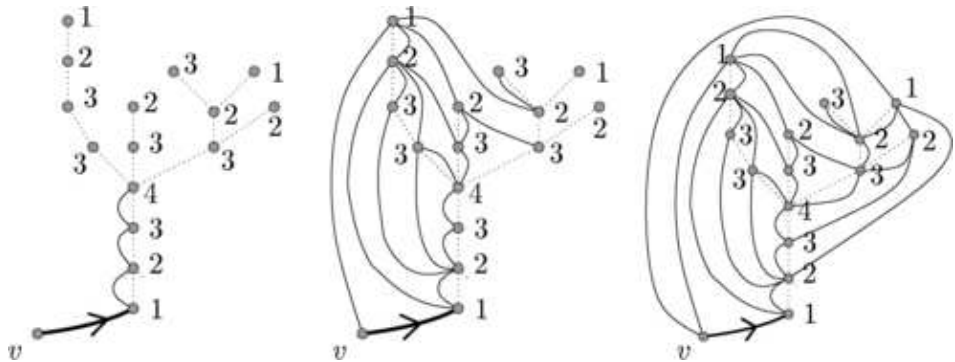

FIG. 8. *Construction of $Q(T)$ using the CDFT; on the first figure, the four first edges are drawn and the construction is completed in the next pictures.*

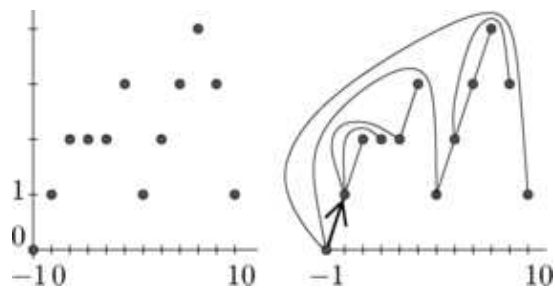

FIG. 9. *A process $\mathcal{R}$ and the associated doddering tree $\mathcal{D}(\mathcal{R})$.*

Since the nodes of $Q(T)$ are encoded several times in $\mathcal{D}(R_n^+)$, the second step consists in the gluing of the nodes of $\mathcal{D}(R_n^+)$ using $V_n^+$.

This representation of rooted quadrangulations using two trees is fundamental to understand the construction of our notion of Brownian map.

2.4.1. *The doddering tree.* Let $N$ be a positive integer, $\mathcal{R}$ a process defined on $[\![0, N]\!]$ satisfying the conditions (2.6), and $\mathcal{P}$ the predecessor function associated with $\mathcal{R}$.

For each $i$ in $[\![-1, N]\!]$, draw the point with coordinates $(i, \mathcal{R}(i))$ in the plane, with the convention $\mathcal{R}(-1) = 0$ (see Figure 9). For each $i$ from 0 to $N$, draw a chord from the point $(i, \mathcal{R}(i))$ to the point $(\mathcal{P}(i), \mathcal{R}(\mathcal{P}(i)))$ that goes above the chords drawn from the points $(j, \mathcal{R}(j))$ such that $\mathcal{P}(i) < j < i$ (this is allowed by Lemma 2.3). Let us call $\mathcal{D}(\mathcal{R})$ the planar map whose edges are the $N + 1$ drawn chords, whose vertices are the $N + 2$ points $(i, \mathcal{R}(i))$, for $i \in [\![-1, N]\!]$, and whose root is the oriented chord from $(-1, 0)$ to $(0, 1)$.

PROPOSITION 2.4. *$\mathcal{D}(\mathcal{R})$ is a rooted tree (we call it the doddering tree).*

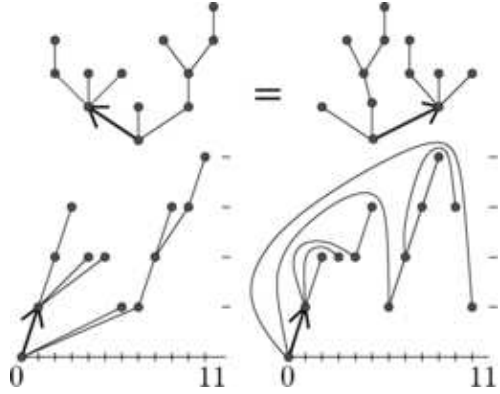

Fig. 10. *On the first line, a rooted tree $\tau$ is drawn in two ways (via a circular permutation of the edges adjacent to the root-vertex). The points of the two pictures in the second line represent respectively the CHP and the RHP of $\tau$. The reconstruction of $\tau$ is done on the CHP thanks to a "lying tree," and in the RHP by the "doddering procedure."*

PROOF. Let us denote by $m(i)$ the vertex $(i, \mathcal{R}(i))$. For each $i \geq 0$, $m(i)$ has one and only one adjacent vertex $m(j)$ such that $\mathcal{R}(j) = \mathcal{R}(i) - 1$. Hence, for each $i \geq 0$, there exists a path between $m(i)$ and $m(-1)$ and the planar map is connected. Since there are $N + 2$ vertices and $N + 1$ edges, it is a tree (see Figure 9).  □

*Doddering tree and height process.* Let $\tau$ be a rooted tree with $N + 2$ nodes $v_0, \ldots, v_{N+1}$ sorted according to the clockwise order. For each $k$ in $[\![0, N+1]\!]$, set

$$h'(k) = d(v_0, v_k),$$

the depth of the node $v_k$ in the tree (the root-vertex is $v_0$). The process $(h'(k))_k$ is called the clockwise height process (CHP) of $\tau$. The CHP characterizes the rooted tree $\tau$ (see, e.g., [28]). On the first column of Figure 10, we show how to rebuild $\tau$ given its CHP. Formally, we consider the points $(k, h'(k))$ as the nodes of $\tau$; then the father of $(k, h'(k))$ (for $k$ from 1 to $N+1$) is the node $(p(k), h'(p(k)))$, where

$$p(i) = \max\{j \in [\![0, i-1]\!], h'(j) = h'(i) - 1\} \qquad \text{for } i \in [\![1, N+1]\!].$$

In the CHP, the nodes are represented according to their clockwise orders in the tree. To reconstruct $\tau$ with the good order of the edges, we have to add these edges *below* the CHP (see column 1 of Figure 10).

The reverse height process (RHP) of a tree $\tau$ is the sequence $h''(k) = d(w_0, w_k)$, where the nodes of $\tau$ are sorted according to the reverse order. On the second column of Figure 10, the RHP associated with the tree is drawn. The construction of a tree with the doddering procedure on the

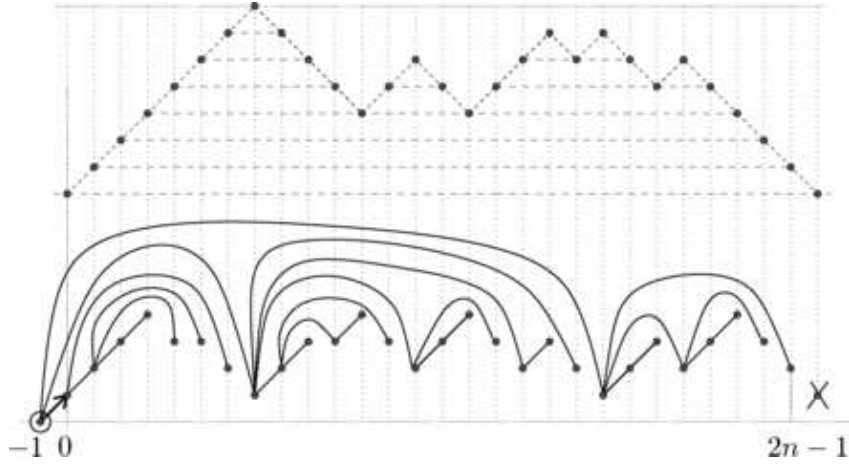

FIG. 11. *Encoding and doddering tree corresponding to the tree drawn in Figure 6. The tree $\mathcal{D}(R_{14}^+)$ is drawn under the graph of $V_{14}^+$ in such a way that the nodes of $\mathcal{D}(R_{14}^+)$ to be glued are easily identified.*

RHP produces the tree drawn in the second column of Figure 10. This is, somehow, the result of a double inversion: the first one is the reverse traversal. The second one is the construction of the doddering tree: indeed, if $k$ and $k'$ have the same predecessor $j$ and if $k < k'$, then $k'$ is a left brother for $k$ in $\mathcal{D}(\mathcal{R})$. As a consequence, we have the following:

PROPOSITION 2.5. *The process $(\mathcal{R}(-1+k))_{k\in[\![0,N+1]\!]}$ is the RHP of the doddering tree $\mathcal{D}(\mathcal{R})$.*

REMARK 2.6. A RHP is a nonnegative process $(H_k)_{k=0,\ldots,n}$ that satisfies $H_{k+1} - H_k \in \{+1, 0, -1, -2, \ldots\}$. The doddering tree $\mathcal{D}(R_n^+)$ used to build $Q(T)$ has a nonusual RHP since $R_n^+$ is a Motzkin type walk: its increments belong to $\{+1, 0, -1\}$ and $R_n^+(2n) = 1$.

2.4.2. *Construction of $Q(T)$ with the doddering tree.* Let $T$ be an element of $\mathcal{W}_n^+$ and $(R_n^+, V_n^+)$ its encoding. Recall that each $k$ in $[\![1, 2n-1]\!]$ encodes on $V_n^+$ exactly one corner of $T$, whereas 0 and $2n$ encode the same corner of the root-vertex of $T$. We call gluer tree, and denote by $\mathcal{G}(V_n^+)$, the tree with CDFW $V_n^+$, that is, the underlying tree of $T$.

Here is now, in two steps, the new procedure to construct $Q(T)$:

(I) *Construction of a doddering tree with $2n$ edges.* Consider $R_n^+$ as a process on $[\![0, 2n-1]\!]$. Draw $\mathcal{D}(R_n^+)$. (See Figure 11.)

(II) *Gluing of the nodes of $\mathcal{D}(R_n^+)$.* In view of the description of $Q$ given in Section 2.3, the edges of $\mathcal{D}(R_n^+)$ are exactly the edges of $Q(T)$. The root-edge of $Q(T)$ is the oriented edge $\overrightarrow{(-1,0),(0,1)}$. The point $(-1,0)$ of the plane is the root-vertex of $Q(T)$. The other vertices of $Q(T)$ are represented by one or several nodes of $\mathcal{D}(R_n^+)$. In order to build $Q(T)$, we have to glue some nodes of $\mathcal{D}(R_n^+)$. The set of nodes of $T$ is exactly the set of vertices of $Q(T)$ different from $(-1,0)$. Two nodes of $\mathcal{D}(R_n^+)$ with abscissas $i$ and $j$ must be glued iff $i$ and $j$ are corners of the same node of $\mathcal{G}(V_n^+)$. If $i$ and $j$ encode the same node $u$ of $\mathcal{G}(V_n^+)$, the abscissas in $\rrbracket i,j\llbracket$ encode the nodes of $T$ that are descendants of $u$. According to the description of $Q$ (in Section 2.3), the nodes $(i, R_n^+(i))$ and $(j, R_n^+(j))$ must be glued in such a way to envelop by below all the nodes represented in $V_n^+$ by abscissas in $\rrbracket i,j\llbracket$.

There are at least two ways to see how the gluing of the nodes works:

• The first one is illustrated in Figure 12 and 13: by a homeomorphism of the plane, we send the doddering tree built on the points $(k, R_n^+(k))$ on the points $(k, V_n^+(k))$ [this can be done by drawing the chords directly from the points $(k, V_n^+(k))$ to the points $(\mathcal{P}(k), V_n^+(\mathcal{P}(k)))$]. In that way, the points to be glued are clearly characterized: they are the points corresponding to

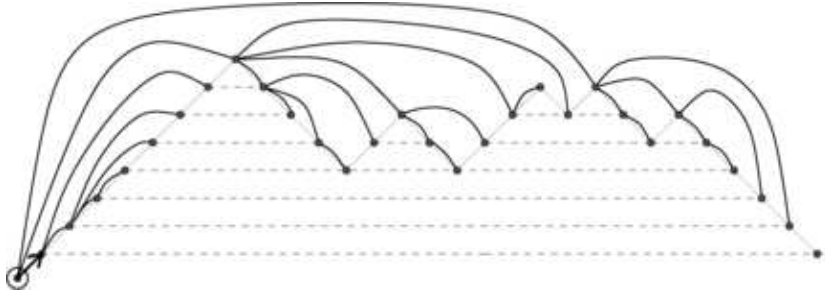

FIG. 12.  *Displacement of $\mathcal{D}(R_{14}^+)$ for the identification of the nodes.*

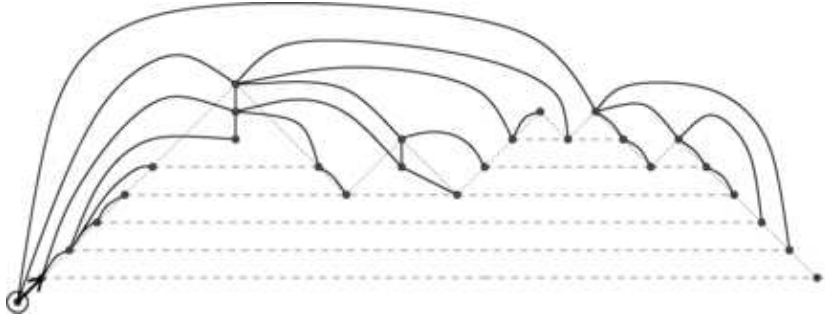

FIG. 13.  *Illustration of three gluings corresponding to Figure 12. The quadrangulation is obtained by the identification of some nodes of $\mathcal{D}(R_n^+)$; seen as a graph, the quadrangulation it then a quotient graph. As illustrated in this figure, the identification (gluing) operation is planar and induces a suitable embedding of the quotient graph in the plane.*

the same node of the gluer tree $\mathcal{G}(V_n^+)$. (The gluer tree is the dotted tree in Figures 7 and 8.)

• The second one consists in doing the gluings (as on Figure 14) below the doddering tree. Suppose that $j$ and $k$ (with $j < k$) must be glued and that all the gluings between these two nodes have been done. Then, we take the node $k$ and pull it below the nodes present in $\rrbracket j, k \llbracket$ until its position equals the one of the node $j$. We do this job for all couples of points to be glued.

## 3. Notion of rooted abstract maps and main results.

3.1. *Introduction.* We saw in the previous section that each quadrangulation is a quotient space: the doddering tree $\mathcal{D}(R_n^+)$ appears to be the quadrangulation unfolded, and the gluer tree $\mathcal{G}(V_n^+)$ characterizes the nodes to be glued. We want now to pass at the limit in this construction.

We denote by $(\mathbf{R}_n^+, \mathbf{V}_n^+)$ the encoding under the distribution induced by $\mathbb{P}_D^n$, and by $(\mathcal{D}(\mathbf{R}_n^+), \mathcal{G}(\mathbf{V}_n^+))$ the corresponding pair of trees. We consider now $\mathbf{R}_n^+$ and $\mathbf{V}_n^+$ as continuous processes on $[0, 2n]$ by interpolating linearly between integer abscissa. We introduce

$$\mathbf{v}_n^+(s) = \frac{\mathbf{V}_n^+(2ns)}{\sqrt{n}} \quad \text{and} \quad \mathbf{r}_n^+(s) = \frac{\mathbf{R}_n^+(2ns) - 1}{n^{1/4}} \quad \text{for } s \in [0, 1],$$

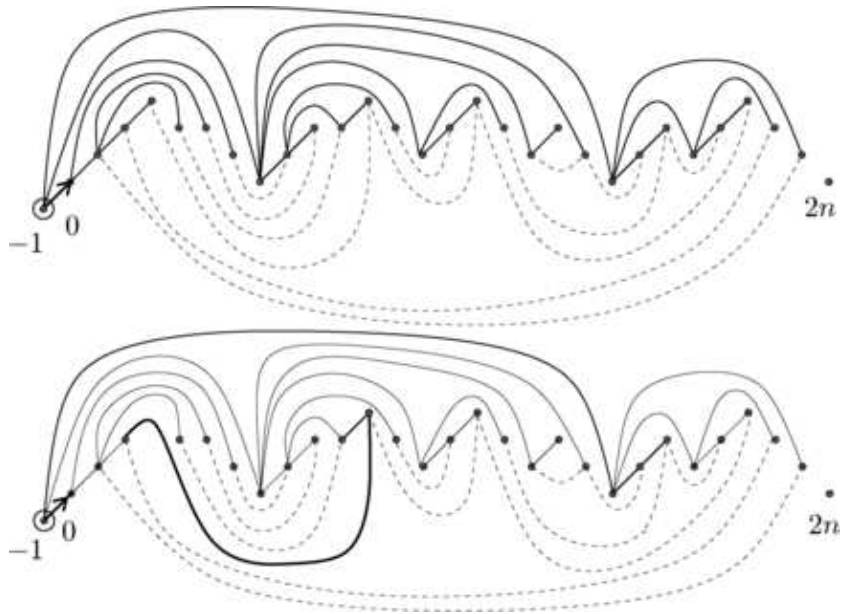

FIG. 14. *The identifications to be made are drawn under the doddering tree on the first picture. The second picture shows a first gluing.*

the "normalized" version of the encoding under $\mathbb{P}_D^n$. The process $(\mathbf{r}_n^+, \mathbf{v}_n^+)$ takes its values in $\mathbb{T}$:

DEFINITION 3.1. For any function $g$ defined on an interval $I$ of $\mathbb{R}$, and for any $x, y \in I$, set $\check{g}(x,y) = \min\{g(u), u \in [x \wedge y, x \vee y]\}$. We denote by $\mathbb{T}$ the subspace of $(C[0,1])^2$ of functions $(f, \zeta)$ satisfying

$$\begin{cases} \zeta(0) = \zeta(1) = 0, \quad \zeta \geq 0, \\ f(0) = f(1) = 0, \\ \text{for any } 0 \leq s \leq s' \leq 1, \text{ if } \zeta(s) = \zeta(s') = \check{\zeta}(s, s'), \quad \text{then } f(s) = f(s'). \end{cases}$$

The space $\mathbb{T}$ is the states space of the tour of the Brownian snake with lifetime process the normalized Brownian excursion (see [29]). We refer also to [18, 25, 26] for considerations on the Brownian snake. The interpretation of the third condition is the following. The function $\zeta$ encodes a tree $\mathcal{T}$ (see Section 3.2) and $f$ is a function compatible with $\zeta$: if $s$ and $s'$ are such that $\check{\zeta}(s, s') = \zeta(s) = \zeta(s')$, then $s$ and $s'$ encode the same point of $\mathcal{T}$ and the third condition ensures that $f$ is a function of the points of $\mathcal{T}$; this property is called the snake property.

We endow $\mathbb{T}$ with the metric

$$d_\mathbb{T}((f_1, \zeta_1), (f_2, \zeta_2)) = \|f_1 - f_2\|_\infty + \|\zeta_1 - \zeta_2\|_\infty.$$

Consider the $\mathbb{T}$-valued random variable $(\mathbf{r}, \mathbf{v})$ whose distribution follows:

- $(\mathbf{v}(t))_{t \in [0,1]} \stackrel{law}{=} \sqrt{2}(\mathbf{e}(t))_{t \in [0,1]}$, where $\mathbf{e}$ is the normalized Brownian excursion.
- Knowing $\mathbf{v}$, the process $\mathbf{r}$ is a Gaussian process with mean 0 and covariance function

$$\mathbb{E}(\mathbf{r}(s)\mathbf{r}(t)) = \sqrt{2/3}\check{\mathbf{v}}(s, t).$$

We denote by $\mathbb{P}_\mathcal{S}$ the law on $\mathbb{T}$ of $(\mathbf{r}, \mathbf{v})$. The process $(\mathbf{r}, \mathbf{v})$ is, up to a scale factor, the tour of the Brownian snake with lifetime process $\sqrt{2}(\mathbf{e}(t))_{t \in [0,1]}$.

We consider the group $\mathcal{O} = ([0, 1), \oplus)$, where $\oplus$ is the addition modulo 1. For $\theta \in \mathcal{O}$, the *rerooting operator* $J_\theta$ is

$$J^{(\theta)} : \mathbb{T} \longrightarrow C[0,1] \times C[0,1],$$
(3.1)
$$(f, \zeta) \longmapsto (f, \zeta)^{(\theta)} = (f^{[\theta]}, \zeta^{(\theta)}),$$

where $(f^{[\theta]}, \zeta^{(\theta)})$ is defined by

(3.2) $\begin{cases} f^{[\theta]}(x) = f(\theta \oplus x) - f(\theta), & \text{for any } x \in [0, 1], \\ \zeta^{(\theta)}(x) = \zeta(\theta \oplus x) + \zeta(\theta) - 2\check{\zeta}(\theta \oplus x, \theta), & \text{for any } x \in [0, 1]. \end{cases}$

The following property is proven in the Appendix:

PROPOSITION 3.2. (i) *For any $\theta \in \mathcal{O}$, $J^{(\theta)}$ takes its values in $\mathbb{T}$.*
(ii) *The operators $J^{(\theta)}$ define a group action of $\mathcal{O}$ on $\mathbb{T}$; that is, for any $(\theta, \theta') \in \mathcal{O}^2$,*

$$J^{(\theta)} \circ J^{(\theta')} = J^{(\theta \oplus \theta')}.$$

Let $\theta^\star = \theta^\star_{\mathbf{r}} = \min \operatorname{argmin} \mathbf{r}$ be the first time where $\mathbf{r}$ reaches its minimum. The proof of the following theorem is given in Section 4.7.

THEOREM 3.3. *The following weak convergence holds in $\mathbb{T}$:*

$$(\mathbf{r}_n^+, \mathbf{v}_n^+) \xrightarrow[n]{weakly} (\mathbf{r}^{[\theta^\star]}, \mathbf{v}^{(\theta^\star)}).$$

(Recall that the process $\mathbf{r}_n^+$ is nonnegative.)

NOTATION. We write $\mathbb{P}^+$ the law of $(\mathbf{r}^{[\theta^\star]}, \mathbf{v}^{(\theta^\star)})$. In the sequel $(\mathbf{r}^+, \mathbf{v}^+)$ will denote a random variable $\mathbb{P}^+$-distributed.

Our aim is to construct a limit map with the help of two continuous trees encoded by $(\mathbf{r}^+, \mathbf{v}^+)$. For this purpose, we need first to carefully (re)define the notion of tree. Usually, we represent a discrete rooted tree with $n$ nodes in $\mathbb{R}^2$ as a continuous planar map. This representation is not adequate for the limit trees. We define here a notion of *tree* that covers the model of scaled gluer trees, the model of scaled doddering trees and their continuous limits, but also the model of finite trees with random edges length. In the present work the notion of nodes, the notion of cyclic orders around nodes and the notion of tree traversals are particularly important. We need a definition of trees that takes into account these notions. We choose to encode the presence and the "quantity" of nodes in a region of the tree with the help of a measure. Hence, each tree will be encoded by a measure and a DFW.

A large part of our description of trees is inspired from [1, 2, 18, 19, 31]. For other considerations on trees, see also [3, 4, 5] and [13].

When a suitable definition of ordered rooted trees will be given, we will introduce our notion of abstract maps. It is quite important to have in mind that the goal is to define an abstract map with a pair of trees $(\mathcal{D}, \mathcal{G})$, in order to generalize the description of quadrangulations.

3.2. *Notion of abstract trees.* For $a > 0$, consider $C^+[0, a]$ the set of continuous functions $g$ from $[0, a]$ to $\mathbb{R}^+$ that satisfy $g(0) = g(a) = 0$. For any $g \in C^+[0, a]$, we introduce the equivalence relation in $[0, a]$,

$$x \underset{g}{\sim} y \iff g(x) = g(y) = \check{g}(x, y).$$

We denote by $E_g$ the quotient space $[0,a]/\underset{g}{\sim}$ and by $F_g$ the canonical surjection from $[0,a]$ onto $E_g$. For short, we write sometimes $\dot{x}$ instead of $F_g(x)$ and we say that $x$ is a representative of $\dot{x}$. Let $\mathcal{M}(a)$ be the set of finite measures on $[0,a]$. For $\mu \in \mathcal{M}(a)$, set $E_\mu = F_g(\operatorname{supp}(\mu))$ the image by $F_g$ of the support of $\mu$.

DEFINITION 3.4. A pair $(g,\mu) \in C^+[0,a] \times \mathcal{M}(a)$ is said to be a tree-encoding if $E_g^{(T)} \stackrel{def}{=} \{u \in E_g, \#F_g^{-1}(u) \neq 2\} \cup \{\dot{0}\}$ satisfies:

$$E_g^{(T)} \subset E_\mu. \tag{3.3}$$

Let $(g,\mu)$ be a tree encoding. For any $\dot{x}$ and $\dot{y}$ in $E_g$, set

$$d_{E_g}(\dot{x},\dot{y}) = g(x) + g(y) - 2\check{g}(x,y).$$

It is not difficult to check that $d_{E_g}$ is a metric on $E_g$ and that, for any $x \in [0,a]$, $g(x) = d_{E_g}(\dot{0},\dot{x})$.

DEFINITION 3.5. Let $(g,\mu)$ be a tree encoding. The rooted tree $\mathcal{T}$ clockwise encoded by $(g,\mu)$ [we write $\mathcal{T} = \operatorname{CTree}(g,\mu)$] is the metric space $\mathcal{T} = (E_g, d_{E_g})$. The function $F_g$ is called the CDFT of $\mathcal{T}$, the elements of $E_g$ are called points of $\mathcal{T}$, the elements of $E_\mu$ are called nodes of $\mathcal{T}$, the class $F_g(0) = \dot{0}$ is called the root-vertex of $\mathcal{T}$, and the function $g$, the CDFW of $\mathcal{T}$.

We often use the notation $E_\mathcal{T}$, $F_\mathcal{T}$, $d_\mathcal{T}$ instead of $E_g$, $F_g$, $d_g$.

REMARK 3.6. Condition (3.3) imposes to consider as nodes the root and the points giving some topological informations: the points with degree 1 (the leaves) and the points with degree larger or equal to 3. The measure $\mu$ gives information on the repartition of the nodes in the tree.

Since $F_\mathcal{T}$ is continuous, and since $[0,a]$ is compact and path-connected, we have the following:

LEMMA 3.7. $E_\mathcal{T}$ *is a path-connected compact metric space, and it is parameterized by* $[0,a]$.

Since $\operatorname{supp}(\mu)$ is compact, the set of nodes $E_\mu$ is compact and thus closed. The measure $\mu \circ F_\mathcal{T}^{-1}$ is a measure on $E_\mathcal{T}$ and its support is $E_\mu$.

The parameter $a$ expresses in some sense the size or the total weight of the nodes of the tree.

*Arborescent structure of $\mathcal{T}$ and order on $\mathcal{T}$.* We define some notions related to this construction of trees:

DEFINITION 3.8.  (i) The set of corners of $\mathcal{T}$ is $[0, a)$. The set of corners around a point $\dot{x}$ is $F_{\mathcal{T}}^{-1}(\dot{x}) \cap [0, a)$. The corner 0 is the root-corner.

(ii) For $u \in E_{\mathcal{T}}$, $\deg(u) := \# F_{\mathcal{T}}^{-1}(u) \cap [0, a)$ is called the (total) degree of $u$.

The following proposition is classical:

PROPOSITION 3.9.  *Let $x$ and $y$ be representatives of $\dot{x}$ and $\dot{y}$. If $z$ is such that $z \in [x, y]$ and $g(z) = \check{g}(x, y)$, then the class $\dot{z}$ does not depend on the representatives of $\dot{x}$ and $\dot{y}$. The point $\dot{z}$ is called the deepest common ancestor of $\dot{x}$ and $\dot{y}$.*

DEFINITION 3.10.  Let $\dot{x} \in E_{\mathcal{T}}$ and we denote the smallest and largest representatives of $\dot{x}$ in $\mathcal{T}$ by

$$\underline{\rho}(\dot{x}) = \inf\{y \in [0, a], F_{\mathcal{T}}(y) = \dot{x}\}$$

and

$$\overline{\rho}(\dot{x}) = \sup\{y \in [0, a], F_{\mathcal{T}}(y) = \dot{x}\}.$$

The interval $[\underline{\rho}(\dot{x}), \overline{\rho}(\dot{x})]$ is called the subtree rooted at $\dot{x}$.

It is straightforward that:

LEMMA 3.11.  (i) *If $[\underline{\rho}(\dot{x}), \overline{\rho}(\dot{x})] \cap [\underline{\rho}(\dot{y}), \overline{\rho}(\dot{y})] \neq \varnothing$, then $[\underline{\rho}(\dot{x}), \overline{\rho}(\dot{x})] \subset [\underline{\rho}(\dot{y}), \overline{\rho}(\dot{y})]$ or $[\underline{\rho}(\dot{y}), \overline{\rho}(\dot{y})] \subset [\underline{\rho}(\dot{x}), \overline{\rho}(\dot{x})]$. In the first case we say that $\dot{x}$ is a descendant of $\dot{y}$ or that $\dot{y}$ is an ancestor of $\dot{x}$.*

(ii) *If ($\dot{x} = \dot{y}$, $\dot{z} \neq \dot{x}$ and $x < z < y$), then $\dot{z} \subset (x, y)$. In particular, $\dot{z}$ is a descendant of $\dot{x}$.*

DEFINITION 3.12.  The clockwise order $\preccurlyeq_{CO}$ is defined by

$$\dot{x} \preccurlyeq_{CO} \dot{y} \Leftrightarrow \underline{\rho}(\dot{x}) \leq \underline{\rho}(\dot{y}).$$

It is a total order on $E_{\mathcal{T}}$. If $\dot{x}$ is an ancestor of $\dot{y}$, then $\dot{x} \preccurlyeq_{CO} \dot{y}$. Around each point $\dot{x}$, the clockwise cyclic order is defined as follows: first, the points represented in $[0, \underline{\rho}(\dot{x}))$ with the clockwise order, then those represented in $[\underline{\rho}(\dot{x}), \overline{\rho}(\dot{x})]$ and then those represented in $(\overline{\rho}(\dot{x}), 1)$.

DEFINITION 3.13.  If $x_1 < x_2 < x_3$ are representatives of $\dot{x}$, then $[x_1, x_2]$ and $[x_2, x_3]$ are called subtrees of $[\underline{\rho}(\dot{x}), \overline{\rho}(\dot{x})]$ and we say that $[x_1, x_2]$ is before $[x_2, x_3]$ (with respect to the clockwise order).

*Geodesics, branches and cycles.*

DEFINITION 3.14. Let $u \in E_\mathcal{T}$. We call branch $\mathcal{S}_u$ between $u$ and the root-vertex the set of ancestors of $u$.

The branch $\mathcal{S}_u$ is the continuous curve parameterized as follows. Set $x$ such that $\dot{x} = u$ and for $s \in [0, g(x)]$, set $m(s) = \sup\{y \in [0, x], g(y) = s\}$. The function $s \longmapsto F_\mathcal{T}(m(s))$ is a continuous bijection between $[0, g(x)]$ and $\mathcal{S}_u$. Clearly, $\mathcal{S}_u$ is a geodesic in the metric space $E_\mathcal{T}$. In the same way, one can see that between any two points $\dot{x}$ and $\dot{y}$, there is a geodesic parameterized by $[\check{g}(x,y), g(x)] \cup (\check{g}(x,y), g(y)]$.

LEMMA 3.15. $\mathcal{T}$ *has no cycle (i.e., no subset homeomorphic to a circle) and so between any two points $u$ and $v$ in $E_\mathcal{T}$, there is a unique geodesic.*

PROOF. First, a cycle can not be included in a branch. So, a cycle $C$ must contain two points $u$ and $v$ such that $u$ is not an ancestor of $v$ and $v$ is not an ancestor of $u$ (they are different from their deepest common ancestor $w$). One can show by the connectedness argument that $w \in C$, and then that $C \setminus w$ is disconnected. So there is no cycle. □

*Reverse order and reverse depth first walk (RDFW) of $\mathcal{T}$.* Set $\theta = \inf\{t | t \in (0, a], F_\mathcal{T}(t) = F_\mathcal{T}(0)\}$. Let

$$\Psi_\mathcal{T} : [0, a] \longrightarrow [0, a],$$
(3.4)
$$x \longmapsto \begin{cases} \theta - x, & \text{if } x \in [0, \theta], \\ a + \theta - x, & \text{if } x \in ]\theta, a]. \end{cases}$$

DEFINITION 3.16. The function $G_\mathcal{T} : [0, a] \longrightarrow E_\mathcal{T}$, defined by
$$G_\mathcal{T}(x) = F_\mathcal{T}(\Psi_\mathcal{T}(x)) \qquad \text{for } x \in [0, a],$$
is called RDFT of $\mathcal{T}$. The function $l_\mathcal{T}$ defined by
$$l_\mathcal{T}(x) = d(G_\mathcal{T}(x), G_\mathcal{T}(0)) = g(\Psi_\mathcal{T}(x))$$
is called the RDFW of $\mathcal{T}$. The reverse measure $\mu'$ is defined by $\mu'(A) = \mu(\Psi_\mathcal{T}(A))$ for all Borelian of $[0, a]$.

Notice that $\text{supp}(\mu) = \Psi_\mathcal{T}(\text{supp}(\mu'))$ and thus,
$$G_\mathcal{T}(\text{supp}(\mu')) = F_\mathcal{T}(\Psi_\mathcal{T}(\text{supp}(\mu'))) = E_\mu.$$
We say that $u$ is smaller than $v$ for the reverse order (we write $u \preccurlyeq_{RO} v$) if $\min\{G_\mathcal{T}^{-1}(u)\} \leq \min\{G_\mathcal{T}^{-1}(v)\}$. Hence, $\preccurlyeq_{RO}$ is a total order on $\mathcal{T}$; this also induces a reverse cyclic order around each point.

In some sense, $\Psi_\mathcal{T}$ reorders the corners of $\mathcal{T}$ in the reverse order.

*Reversely encoded tree.*

DEFINITION 3.17. The rooted tree $\mathcal{T}$ reversely-encoded by the tree-encoding $(g, \mu)$ [we use the notation $\mathcal{T} = \text{RTree}(g, \mu)$] is defined and described in the same way as $\text{CTree}(g, \mu)$: replace everywhere each occurrence of the word "clockwise" by "reverse," and vice versa.

For the interested reader, relations between CTrees and RTrees are given in Appendix A.1.

*Set of trees with size $a$.* We denote by $C\Gamma(a)$ the set of clockwise-encoded trees with size $a \in \mathbb{R}^+$. Let $d_{C\Gamma(a)} : C\Gamma(a)^2 \to \mathbb{R}^+$ be the application defined for $(\mathcal{T}_1, \mathcal{T}_2) = (\text{CTree}(g_1, \mu_1), \text{CTree}(g_2, \mu_2))$ element of $C\Gamma(a)^2$ by

$$d_{C\Gamma(a)}(\mathcal{T}_1, \mathcal{T}_2) = \|g_1 - g_2\|_\infty^{(a)} + d_{\mathcal{M}(a)}(\mu_1, \mu_2),$$

with $\|g_1 - g_2\|_\infty^{(a)} = \sup\{|g_1(x) - g_2(x)|, x \in [0, a]\}$ and

$$d_{\mathcal{M}(a)}(\mu_1, \mu_2) = \sup_{x \in \mathbb{R}} |C_{\mu_1}(x) - C_{\mu_2}(x)|,$$

where $C_\mu$ is the repartition function of $\mu$ [i.e., $C_\mu(x) = \mu((-\infty, x])$]. The application $d_{C\Gamma(a)}$ is a metric on $C\Gamma(a)$.

In the same manner, we define $R\Gamma(a)$ the set of reversely-encoded tree with size $a$ and metric $d_{R\Gamma(a)}$.

3.3. *The normalized discrete doddering and gluer trees.*

- The rescaled gluer tree $\mathbf{G}_n \in C\Gamma(1)$ is derived from $\mathcal{G}(\mathbf{V}_n^+)$: we set $\mathbf{G}_n = \text{CTree}(\mathbf{v}_n^+, \mu_{\mathbf{G}_n})$, where $\mu_{\mathbf{G}_n} = \frac{1}{2n}\sum_{k=0}^{2n-1} \delta_{k/(2n)}$ (where $\delta_x$ is the Dirac measure on the point $x$). The set of its nodes is then $E_{\mu_{\mathbf{G}_n}} = \{\dot{x}, x = k/(2n), k \in [\![0, 2n]\!]\}$.
- The RHP of the doddering tree is $\mathbf{R}_n^+(-1 + \cdot)$; denote by $\mathbf{C}_n^+$ its RDFW (linearized between integer points). Let $\mathbf{c}_n^+$ be the rescaled version:

$$\mathbf{c}_n^+(t) = n^{-1/4} \mathbf{C}_n^+(2nt) \qquad \text{for } t \in [0, 2].$$

The rescaled doddering tree is $\mathbf{D}_n = \text{RTree}(\mathbf{c}_n^+, \mu_{\mathbf{D}_n}) \in R\Gamma(2)$, where

$$\mu_{\mathbf{D}_n} = \frac{2}{2n+1}\left(\sum_{l=1}^{2n} \delta_{m(l)/(2n)} + \frac{1}{d} \sum_{t, 0 \le t < 1, \mathbf{c}_n^+(t) = 0} \delta_t\right),$$

where $d$ is the degree of the root of $\mathbf{D}_n$, $m(0) = 0$, and for $l \in [\![1, 2n]\!]$,

(3.5) $$m(l) = \inf\left\{j, j > m(l-1), \mathbf{c}_n^+\left(\frac{j}{2n}\right) > \mathbf{c}_n^+\left(\frac{j-1}{2n}\right)\right\}.$$

In other words, $m(l)/(2n)$ encodes the first corner of the $(l+1)$th node (according to the RO) of $\mathbf{D}_n$. The set of the $2n+1$ nodes of $\mathbf{D}_n$ is $E_{\mu_{\mathbf{D}_n}} = \{\dot{x}, x = k/(2n), k \in [\![0, 4n]\!]\}$. We attribute to each nonroot node the same weight, and we place it on its first corner because the nodes are glued in their first corner during the gluing procedure. The root of $\mathbf{D}_n$ is not glued. We choose to share its weight among its corners for reasons that will appear clearer when studying the pointed quadrangulations. We choose to encode $\mathbf{D}_n$ on $[0, 2]$ in order to stress that it contains twice as many nodes as $\mathbf{G}_n$.

The pair $(\mathbf{D}_n, \mathbf{G}_n)$ is a random variable taking its values in $R\Gamma(2) \times C\Gamma(1)$.

REMARK 3.18. The asymptotics of $\mu_{\mathbf{D}_n}$ is the same if it is only subject to give weight $1/(2n+1)$ to each node, whatever is the repartition between the corners (see the proof of Proposition 3.20).

3.3.1. *Convergence of* $(\mathbf{D}_n, \mathbf{G}_n)$. The following lemma, proved in the Appendix, illustrates the proximity of the RHP and the RDFW processes associated with discrete trees:

LEMMA 3.19. *For $n \geq 0$, let $A_n$ denote a set of rooted trees with $n$ edges, endowed with a probability $\mu_n$. Let $\widetilde{H}_n$ and $\widetilde{V}_n$ be the associated RHP and RDFW. Assume there exists a sequence of real numbers $(c_n)_{n \geq 0}$ such that $c_n \longrightarrow +\infty$, $c_n = o(n)$, and*

$$(c_n^{-1} \widetilde{H}_n(nt))_{t \in [0,1]} \xrightarrow[n]{weakly} (h(t))_{t \in [0,1]},$$

*where the process $h$ is a.s. nonnegative and continuous on $[0, 1]$; then*

$$\sup_{t \in [0,1]} c_n^{-1} |\widetilde{V}_n(2nt) - \widetilde{H}_n(nt)| \longrightarrow 0 \text{ in probability.}$$

*In particular, this yields*

$$\left( \frac{\widetilde{V}_n(2nt_1)}{c_n}, \frac{\widetilde{H}_n(nt_2)}{c_n} \right)_{(t_1, t_2) \in [0,1]^2} \xrightarrow[n]{weakly} (h(t_1), h(t_2))_{(t_1, t_2) \in [0,1]^2}.$$

Of course, the same result holds for the CHP and the CDFW.

This strong relation between the DFW and the height process legitimates the following point of view. The process $\mathbf{r}_n^+$ is the RHP of the normalized doddering tree. By Lemma 3.19, the weak convergence of $(\mathbf{r}_n^+, \mathbf{v}_n^+)$ to $(\mathbf{r}^+, \mathbf{v}^+)$ implies

$$(3.6) \qquad (\mathbf{c}_n^+, \mathbf{v}_n^+) \xrightarrow[n]{weakly} (\mathbf{c}^+, \mathbf{v}^+),$$

where $\mathbf{c}^+(\cdot) = \mathbf{r}^+(\cdot/2)$, and where the convergence holds in $C[0,2] \times C[0,1]$ endowed with the topology of uniform convergence.

Let $(\mathbf{D}_\infty, \mathbf{G}_\infty)$ be the $R\Gamma(2) \times C\Gamma(1)$-valued random variable defined by

(3.7)
$$\mathbf{D}_\infty = \mathrm{RTree}(\mathbf{c}^+, \mathrm{Leb}^{[0,2]}) \in R\Gamma(2),$$
$$\mathbf{G}_\infty = \mathrm{CTree}(\mathbf{v}^+, \mathrm{Leb}^{[0,1]}) \in C\Gamma(1).$$

PROPOSITION 3.20. *The following weak convergence holds:*

$$(\mathbf{D}_n, \mathbf{G}_n) \xrightarrow[n]{weakly} (\mathbf{D}_\infty, \mathbf{G}_\infty) \quad in \ R\Gamma(2) \times C\Gamma(1).$$

PROOF. First, $(\mathbf{c}_n^+, \mathbf{v}_n^+) \xrightarrow[n]{weakly} (\mathbf{c}^+, \mathbf{v}^+)$ and $d_{\mathcal{M}(a)}(\mu_{\mathbf{G}_n}, \mathrm{Leb}^{[0,1]}) \to 0$. For the convergence of $\mu_{\mathbf{D}_n}$, we prove a more general result: if $\mu_{\mathbf{D}_n}$ is only subject to put a weight $1/(2n+1)$ on each node, then

$$\|C_{\mu_{\mathbf{D}_n}} - C_{\mathrm{Leb}^{[0,2]}}\|_\infty \xrightarrow[n]{proba} 0.$$

To this end, we first consider the particular case where $\mu_{\mathbf{D}_n}$ puts all the weight of each node, including the root, on its first corner. For $k \in [\![0, 4n]\!]$, let $\mathcal{N}_k$ be the number of nodes visited before time $k/(2n)$ in $\mathbf{D}_n$ with the RDFT. This is also the number of increasing steps on the nonnormalized walk $\mathbf{C}_n^+$ before time $k$. Hence, $\mathcal{N}_k = (k + \mathbf{C}_n^+(k))/2$ and $m(l)$, defined in (3.5), satisfies $m(l) = \inf\{k \,|\, \mathcal{N}_k = l\}$. One has immediately

$$\sup\{|m(l) - 2l|, l \in [\![0, 2n]\!]\} \leq \|\mathbf{C}_n^+\|_\infty.$$

For any $l \in [\![0, 2n]\!]$, we have $C_{\mu_{\mathbf{D}_n}}(m(l)/(2n)) = 2l/(2n+1)$; using that at most $|m(l) - 2l|$ nodes are visited for the first time in the interval $[m(l) \wedge 2l, m(l) \vee 2l]$, for any $l \in [\![0, 2n]\!]$,

(3.8) $$\left|C_{\mu_{\mathbf{D}_n}}\left(\frac{2l}{2n}\right) - \frac{2l}{2n+1}\right| \leq |m(l) - 2l|\frac{2}{2n} \leq \frac{\|\mathbf{C}_n^+\|_\infty}{n} = \frac{\|\mathbf{c}_n^+\|_\infty}{n^{3/4}}.$$

Since $\mathbf{c}_n^+ \xrightarrow[n]{weakly} \mathbf{c}^+$ the right-hand side of (3.8) goes to 0 in probability, and using that $C_{\mu_{\mathbf{D}_n}}$ is nondecreasing, $\|C_{\mu_{\mathbf{D}_n}} - C_{\mathrm{Leb}^{[0,2]}}\|_\infty \xrightarrow[n]{proba} 0$. Now, by symmetry, the same result holds if one places the mass on the last corner of each node. Hence, if each node has mass $(2n+1)^{-1}$, whatever is the mass repartition on the corners, $\mu_{\mathbf{D}_n}$ converges to $\mathrm{Leb}^{[0,2]}$. □

3.4. *Rooted abstract maps.* We construct an abstract map thanks to a doddering tree $\mathcal{D}$, a gluer tree $\mathcal{G}$ and a gluing function $b$.

DEFINITION 3.21. Let $a \in \mathbb{R}^+$,
$$(\mathcal{D}, \mathcal{G}) = (\mathrm{RTree}(\xi_\mathcal{D}, \mu_\mathcal{D}), \mathrm{CTree}(\zeta_\mathcal{G}, \mu_\mathcal{G})) \in R\Gamma(2a) \times C\Gamma(a)$$
and $b\colon E_{\mu_\mathcal{D}} \setminus \{\textit{root-vertex}\} \to [0, a)$ be an application that sends the nodes of $\mathcal{D}$ (but its root), onto the corners of $\mathcal{G}$. The 3-tuple $(\mathcal{D}, \mathcal{G}, b)$ is said to be $a$-admissible if the three following conditions are satisfied:

(i) $b$ is an injection.
(ii) $b$ is increasing: if $u \preccurlyeq_{RO} v$ in $\mathcal{D}$, then $b(u) \leq b(v)$ in $[0, a)$.
(iii) If $u$ and $v$ are two nodes in $\mathcal{D}$ such that $b(u) \underset{\zeta_\mathcal{G}}{\sim} b(v)$, then $\xi_\mathcal{D}(u) = \xi_\mathcal{D}(v)$.

Let $(\mathcal{D}, \mathcal{G}, b)$ be $a$-admissible. We define an equivalence relation on $E_\mathcal{D}$: for $x, y \in E_\mathcal{D}$, we say that

(3.9) $\quad x \underset{b}{\sim} y \Leftrightarrow (x = y) \quad \text{or} \quad \left( \{x, y\} \subset \mathbb{E}_{\mu_\mathcal{D}} \setminus \{\text{root-vertex}\} \text{ and } b(x) \underset{\zeta_\mathcal{G}}{\sim} b(y) \right).$

For $x \in E_\mathcal{D}$, we set $\hat{x} = \{y \in E_\mathcal{D}, y \underset{b}{\sim} x\}$. A class $\hat{x}$ is either a point of $\mathcal{D}$, or the set of the nodes of $\mathcal{D}$ glued with $x$ ($x$ included), or the root of $\mathcal{D}$.

Let $M$ be the set
$$M = \{\hat{x}, x \in E_\mathcal{D}\}.$$

For any $\hat{u}, \hat{w} \in M$ and any $k > 0$, set
$$d^{(k)}(\hat{u}, \hat{w}) = \inf \sum_{i=0}^{k} d_\mathcal{D}(u_{2i}, u_{2i+1}),$$
where the infimum is taken on the set $\Gamma = \{(u_0, \ldots, u_{2k+1}) \in E_\mathcal{D}^{2k+2} \text{ s.t. } \hat{u}_0 = \hat{u}, \hat{u}_{2k+1} = \hat{w}, \hat{u}_{2i+1} = \hat{u}_{2i+2}\}$ and where $d_\mathcal{D}$ is the metric in $\mathcal{D}$. The application $d_M \colon M^2 \to \mathbb{R}^+$, defined for any $\hat{u}, \hat{w} \in M$ by
$$d_M(\hat{u}, \hat{v}) = \inf_{k \geq 0} d^{(k)}(\hat{u}, \hat{v}),$$
is a metric on $M$.

DEFINITION 3.22. The metric space $(M, d_M)$ is called the rooted map encoded by $(\mathcal{D}, \mathcal{G}, b)$. We denote this space $\mathrm{Map}(\mathcal{D}, \mathcal{G}, b)$. The elements of $M$ are called points. The real $a > 0$ is called the size of $M$.

The reader will find in Section 3.5 and in the sections that follow the representation of normalized rooted quadrangulations as abstract maps, the definition of the Brownian map and the convergence of rescaled quadrangulations to the Brownian map. For the moment, we give some properties of abstract rooted maps.

We have a canonical surjection

$$S: E_\mathcal{D} \longrightarrow M,$$
$$x \longmapsto \hat{x}.$$

REMARK 3.23 (About the metric $d_M$). Each element $\gamma \in \Gamma$ defines a path in the map:

– between $\hat{u}_{2i}$ and $\hat{u}_{2i+1}$ it is the image by $S$ of the geodesic between $u_{2i}$ and $u_{2i+1}$ in $\mathcal{D}$,
– $\hat{u}_{2i+1} = \hat{u}_{2i+2}$ (since $u_{2i+1}$ and $u_{2i+2}$ are identified).

LEMMA 3.24. (1) *The canonical surjection $S$ is $1$-Lipschitz.* (2) *$M$ is a compact path-connected metric space.*

PROOF. (1) by definition of $d_M$, $d_M(\hat{x}, \hat{y}) \leq d_\mathcal{D}(x, y)$. (2) $S$ is continuous and $E_\mathcal{D}$ is compact path-connected. □

3.4.1. *Corners, degree, DFT, order in abstract maps.* We define now some notions, and prove some properties related to the abstract map

(3.10) $\quad M = \mathrm{Map}(\mathcal{D}, \mathcal{G}, b) = \mathrm{Map}(\mathrm{RTree}(\xi_\mathcal{D}, \mu_\mathcal{D}), \mathrm{CTree}(\zeta_\mathcal{G}, \mu_\mathcal{G}), b).$

(a) Let $F_\mathcal{D}: [0, 2a] \longrightarrow E_\mathcal{D}$ the RDFT of $\mathcal{D}$. We call *depth first traversal (DFT)* of $M$ the application $S \circ F_\mathcal{D}: [0, 2a] \longrightarrow M$. It is a continuous parameterization of $M$ by $[0, 2a]$.

(b) The set of nodes of $M$ is $S(E_{\mu_\mathcal{D}})$ the image by $S$ of the nodes of $\mathcal{D}$. The root-vertex of $M$ is $S(F_\mathcal{D}(0))$ (the image by $S$ of the root-vertex of $\mathcal{D}$).

(c) The set of corners of $M$ is the set $[0, 2a)$. The set of corners around the point $\hat{x}$ is the set $F_\mathcal{D}^{-1}(S^{-1}(\hat{x})) \cap [0, 2a)$. The root-corner of $M$ is 0.

(d) The degree of a point $\hat{x}$ of $M$ is

$$\deg \hat{x} = \# F_\mathcal{D}^{-1}(S^{-1}(\hat{x})) \cap [0, 2a).$$

(e) The cyclic order around $\hat{x}$ is the following order ($\preccurlyeq$) on the set of corners: set $s, t \in F_\mathcal{D}^{-1}(S^{-1}(\hat{x})) \cap [0, 2a)$.
– If $s$ and $t$ are corners of the same node in $\mathcal{D}$, then $s \preccurlyeq t$ in $M$ if $s \preccurlyeq_{CO} t$ in $\mathcal{D}$.

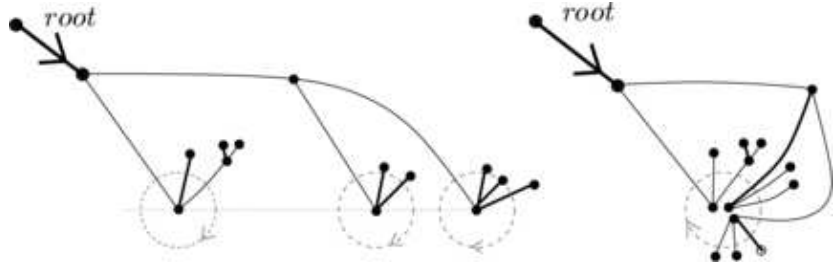

Fig. 15. *Cyclic orders around the nodes in the map induced by the order in the doddering tree and the gluings. This is consistent with the gluing procedure of the discrete doddering tree.*

– If $s$ and $t$ are corners of different nodes $u$ and $v$ in $\mathcal{D}$, then $s \preccurlyeq t$ in $M$ if $b(u) \leq b(v)$ in $[0, a]$.

(See an illustration on Figure 15.)

(f) Set $\underline{\rho}(\hat{x})$ and $\bar{\rho}(\hat{x})$ the smallest and largest corners of $\hat{x}$:

$$\underline{\rho}(\hat{x}) = \inf\{s \in [0, 2a], S(F_\mathcal{D}(s)) = \hat{x}\}, \bar{\rho}(\hat{x}) = \sup\{s \in [0, 2a], S(F_\mathcal{D}(s)) = \hat{x}\}.$$

We define a total order on $M$ by setting

$$\hat{x} \leq \hat{y} \Leftrightarrow \underline{\rho}(\hat{x}) \leq \underline{\rho}(\hat{y}).$$

3.4.2. *Topology and geometry of abstract maps.* The geometry of the abstract map $M$ in (3.10) will be described with the help of "simple geodesic," "cycles" and "submaps." We give up the notions of edges and faces that seem to be not suitable for continuous maps.

*Simple geodesics.*

LEMMA 3.25. *Let $x, y \in E_\mathcal{D}$ and $s, t \in [0, 2a]$ such that $x = F_\mathcal{D}(s)$ and $y = F_\mathcal{D}(t)$.*

(1) *Assume that $\xi_\mathcal{D}(t) - \xi_\mathcal{D}(s) \geq 0$, then $d_M(\hat{x}, \hat{y}) \geq \xi_\mathcal{D}(t) - \xi_\mathcal{D}(s)$.*

(2) *If $y$ is a descendant of $x$ in $\mathcal{D}$, then $d_M(\hat{x}, \hat{y}) = \xi_\mathcal{D}(t) - \xi_\mathcal{D}(s)$ and the continuous curve*

$$[\xi_\mathcal{D}(s), \xi_\mathcal{D}(t)] \longrightarrow M,$$
$$u \longmapsto S(F_\mathcal{D}(\alpha(u))),$$

*where $\alpha(u) = \sup\{s' \in [s, t], \xi_\mathcal{D}(s') = u\}$, is a geodesic in $M$ between $\hat{x}$ and $\hat{y}$ (we denote it by $\langle x, y \rangle$).*

We call such geodesic in $M$ simple geodesic. For $x \in E_\mathcal{D}$ the simple geodesic $\langle F_\mathcal{D}(0), x \rangle$ is called a *branch*.

PROOF OF LEMMA 3.25. By the definition of $d_\mathcal{D}$, we have $d_\mathcal{D}(x,y) \geq \xi_\mathcal{D}(t) - \xi_\mathcal{D}(s)$ and $d_\mathcal{D}(x,y) = \xi_\mathcal{D}(t) - \xi_\mathcal{D}(s)$ if $y$ is a descendant of $x$. Now, since for $u_1, u_2 \in E_\mathcal{D}$ and $s_1$ and $s_2$ such that $u_1 = F_\mathcal{D}(s_1)$, $u_2 = F_\mathcal{D}(s_2)$, $u_1 \underset{b}{\sim} u_2 \Rightarrow \xi_\mathcal{D}(s_1) = \xi_\mathcal{D}(s_2)$, we have, for any $x, y \in E_\mathcal{D}$ and any $k \geq 0$, $d^{(k)}(\hat{x}, \hat{y}) \geq \xi_\mathcal{D}(t) - \xi_\mathcal{D}(s)$ and, thus,

$$d_M(\hat{x}, \hat{y}) \geq \xi_\mathcal{D}(t) - \xi_\mathcal{D}(s).$$

When $y$ is a descendant of $x$ in $\mathcal{D}$, $d^{(0)}(x,y) = \xi_\mathcal{D}(t) - \xi_\mathcal{D}(s)$ and then

$$d_M(\hat{x}, \hat{y}) = \xi_\mathcal{D}(t) - \xi_\mathcal{D}(s) = d_\mathcal{D}(x,y). \qquad \square$$

Notice that $(\hat{x}', \hat{y}' \in \langle x, y \rangle$ and $x' \neq y') \Rightarrow \xi_\mathcal{D}(x') \neq \xi_\mathcal{D}(y')$.

*Cycles.* $E_\mathcal{D}$ has no cycle. The map $M$ is obtained from $E_\mathcal{D}$ by the gluing of some nodes. It turns out that each pair of nodes that are glued gives one cycle in $M$.

LEMMA 3.26. *Let $x, y \in E_{\mu_\mathcal{D}}$ such that $x \neq y$, $x \preccurlyeq_{RO} y$ and $\hat{x} = \hat{y}$; let $z$ be the deepest common ancestor of $x$ and $y$ in $\mathcal{D}$. The set $\langle z, x \rangle \cup \langle z, y \rangle$ denoted by $\langle x, z, y \rangle$ is a cycle in $M$. We will call $\langle x, z, y \rangle$ a simple cycle and the node $\hat{x}$ its origin.*

Notice that $z \in E_{\mu_\mathcal{D}}$ since either $z = \text{root}$ or $\deg z \geq 3$.

PROOF OF LEMMA 3.26. We have to show that if $u$ and $v$ (points of $\mathcal{D}$) are such that $\hat{u}$ and $\hat{v}$ belong to $\langle z, x \rangle \cup \langle z, y \rangle$ and $u \neq v$, then $\hat{u} \neq \hat{v}$. We can assume that $u$ and $v$ are different from $x$, $y$ and $z$ since the depth of $z$ (in $\mathcal{D}$) is strictly smaller than the depth of any node (of $\mathcal{D}$) in $\langle z, x \rangle \cup \langle z, y \rangle \setminus \{z\}$ and since the depth of $x$ (and $y$) is strictly larger than the depth of any node in $\langle z, x \rangle \cup \langle z, y \rangle \setminus \{x, y\}$. Now, there are two cases:

– If $\hat{u}$ and $\hat{v}$ are in the same simple geodesic, $f(u) \neq f(v)$ and $\hat{u} \neq \hat{v}$.
– If $\hat{u}$ and $\hat{v}$ are not in the same simple geodesic, assume that $\hat{u} \in \langle z, x \rangle$ and $\hat{v} \in \langle z, y \rangle$. Since $x \preccurlyeq_{RO} y$ and $v \neq z$, we have $b(x) < b(v) < b(y)$ and $b(u) < b(x)$. It follows that $v$ can be glued only with nodes $v'$ such that $b(x) < b(v') < b(y)$ and $u$ only with nodes $u'$ such that $b(u') < b(x)$ or $b(u') > b(y)$; this implies that $\hat{u} \neq \hat{v}$. $\square$

By construction:

- If $\#S^{-1}(\hat{x}) \neq 1$, then the set $\{\{y, z\} \in E_{\mu_\mathcal{D}}, z \neq y, \hat{z} = \hat{y} = \hat{x}\}$ is in bijection with the set of cycles with origin $\hat{x}$.
- If $\#S^{-1}(\hat{x}) = 1$, then $\hat{x}$ is not the origin of a simple cycle.
- The root of the map is not the origin of a simple cycle.

*Submaps.*

DEFINITION 3.27. Let $\langle x, z, y \rangle$ be a simple cycle with origin $\hat{x}$. The set
$$\text{Int}\langle x, z, y \rangle = \{\hat{x}' \in M, \underline{\rho}(x) < \underline{\rho}(x') < \underline{\rho}(y)\} \setminus \{\hat{x}\}$$
is called the *interior* of $\langle x, z, y \rangle$.

If we denote
$$U_{x,y} = \Big\{x' \in E_{\mathcal{D}}, x \underset{RO}{<} x' \underset{RO}{<} y\Big\},$$
then $\text{Int}\langle x, z, y \rangle = S(U_{x,y}) \setminus \hat{x}$ (here $x \underset{RO}{<} x'$ means that $x \underset{RO}{\preccurlyeq} x'$ and $x \neq x'$). The important point is that a node $x' \in U_{x,y}$ can be glued only with a node in $U_{x,y}$ because $b(x) \underset{\zeta_{\mathcal{G}}}{\sim} b(y)$ and $x \underset{RO}{<} x' \underset{RO}{<} y$. Notice that $z$ is a common ancestor in $\mathcal{D}$ of any point of $U_{x,y}$.

DEFINITION 3.28. We call simple submap with origin $\hat{x}$ the set
$$M_{x,y} = \langle x, z, y \rangle \cup \text{Int}\langle x, z, y \rangle.$$

PROPOSITION 3.29. (i) $M_{x,y}$ *is path-connected and compact.*
(ii) *Set* $\hat{x}' \in \text{Int}\langle x, z, y \rangle$ *and* $\hat{x}'' \notin M_{x,y}$; *then any continuous path in* $M$ *between* $\hat{x}'$ *and* $\hat{x}''$ *intersects* $\langle x, z, y \rangle$.

PROOF. (i) We have $M_{x,y} = S(F_{\mathcal{D}}([\underline{\rho}(x), \underline{\rho}(y)]))$.
(ii) Let $\gamma$ be a continuous path in $M$ between $\hat{x}'$ and $\hat{x}''$. Let $K$ be the compact $S(F_{\mathcal{D}}([0, \underline{\rho}(x)] \cup [\underline{\rho}(y), 2a]))$. We have
$$K \cap M_{x,y} \subset \langle x, z, y \rangle.$$
Set $\gamma_1 = \gamma \cap K$ and $\gamma_2 = \gamma \cap M_{x,y}$. The two sets $\gamma_1$ and $\gamma_2$ are nonempty, closed and $\gamma = \gamma_1 \cup \gamma_2$. Since $\gamma$ is connected, $\gamma_1 \cap \gamma_2 \neq \varnothing$. This implies that $\gamma \cap \langle x, z, y \rangle \neq \varnothing$. □

Set $\hat{x}$ such that $\#S^{-1}(\hat{x}) \neq 1$; the set of simple cycles (and simple submaps) with origin $\hat{x}$ is ordered around $\hat{x}$ according to the total CO of $\mathcal{G}$.

There is also an inclusion order described by the following straightforward proposition:

PROPOSITION 3.30. (i) *Let* $\hat{x}_1 = \hat{x}_2 = \hat{x}_3$ *such that* $b(x_1) < b(x_2) < b(x_3)$. *Let* $\langle x_1, z, x_2 \rangle$ *and* $\langle x_2, z', x_3 \rangle$ *be the two simple cycles built by* $(x_1, x_2)$ *and* $(x_2, x_3)$ *and let* $z_0$ *be the deepest node between* $z$ *and* $z'$ *(it is* $z$ *or* $z'$*), then*
$$\text{Int}\langle x_1, z, x_2 \rangle \cap \text{Int}\langle x_2, z', x_3 \rangle = \varnothing$$

*and*

$$M_{x_1,x_2} \cap M_{x_2,x_3} = \langle z_0, x_2 \rangle \quad \text{and} \quad M_{x_1,x_3} = M_{x_1,x_2} \cup M_{x_2,x_3}.$$

(ii) *Let $\hat{x}_1 = \hat{x}_4$ and $\hat{x}_2 = \hat{x}_3$ such that $b(x_1) < b(x_2) < b(x_3) < b(x_4)$, then*

$$M_{x_2,x_3} \subset M_{x_1,x_4}.$$

See an illustration in Figure 16.

The maximal simple submap (with respect to the inclusion order) with origin $\hat{x}$ is the simple submap $M_{x_1,x_2}$, where $x_1, x_2 \in E_{\mu_\mathcal{D}}$ and

$$\hat{x}_1 = \hat{x}, \qquad \underline{\rho}(x_1) = \underline{\rho}(\hat{x}),$$
$$\hat{x}_2 = \hat{x}, \qquad \bar{\rho}(x_2) = \bar{\rho}(\hat{x}).$$

A simple submap $M_{x_1,x_2}$ is minimal (or is said to be a simple face) if there do not exist $x_1', x_2'$ such that $M_{x_1',x_2'} \subset M_{x_1,x_2}$ and $M_{x_1',x_2'} \neq M_{x_1,x_2}$.

3.5. *The set of maps of size $a$.* Let $\mathfrak{M}(a)$ be the set of maps of size $a$, and $d_{\mathfrak{M}(a)} : \mathfrak{M}(a)^2 \to \mathbb{R}^+$ be the application defined by

$$d_{\mathfrak{M}(a)}(M_1, M_2) = d_{R\Gamma(2a)}(\mathcal{D}_1, \mathcal{D}_2)$$
$$+ d_{C\Gamma(a)}(\mathcal{G}_1, \mathcal{G}_2) + \|C_{\mu_{\mathcal{D}_1} \circ b_1^{-1}} - C_{\mu_{\mathcal{D}_2} \circ b_2^{-1}}\|_\infty,$$

where, for $i \in \{1, 2\}$,

$$M_i = \text{Map}(\mathcal{D}_i, \mathcal{G}_i, b_i) = \text{Map}(\text{RTree}(\xi_{\mathcal{D}_i}, \mu_{\mathcal{D}_i}), \text{CTree}(\zeta_{\mathcal{G}_i}, \mu_{\mathcal{G}_i}), b_i)$$

and where the function

$$x \mapsto C_{\mu_{\mathcal{D}_1} \circ b_1^{-1}}(x) = \mu_{\mathcal{D}_1}(b_1^{-1}(-\infty, x]) = \mu_{\mathcal{D}_1}(\{y \in E_{\mu_{\mathcal{D}_1}}, b_1(y) \in (-\infty, x]\})$$

measures the amount of nodes of the doddering trees glued in the corners interval $(-\infty, x]$ of the gluer tree $\mathcal{G}_1$. (Here and in the sequel, for simplicity of the notation, we have denoted $\mu_{\mathcal{D}_i} \circ F_{\mathcal{D}_i}^{-1}$ by $\mu_{\mathcal{D}_i}$. The context suffices to decide which function is involved.) The application $d_{\mathfrak{M}(a)}$ is a metric on $\mathfrak{M}(a)$.

Two maps are close according to $d_{\mathfrak{M}(a)}$ if they are constructed with close trees, and if, moreover, the functions $b_1$ and $b_2$ induce close distributions of the nodes of $\mathcal{D}_1$ and $\mathcal{D}_2$ in the corners of $\mathcal{G}_1$ and $\mathcal{G}_2$.

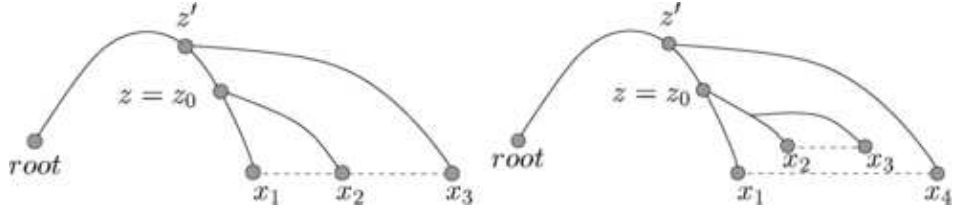

FIG. 16. *Illustration of Proposition* 3.30: *relative positions of the cycles and inclusion order.*

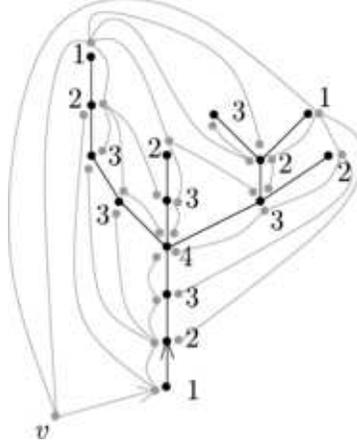

Fig. 17. *The bijection between the set of nonroot nodes of $\mathbf{D}_n$ and the set of corners of the nodes of $\mathbf{G}_n$.*

3.5.1. *Convergence of discrete rooted quadrangulations to the Brownian map.*

*Discrete quadrangulations seen as abstract maps.* Let $(\mathbf{D}_n, \mathbf{G}_n)$ in $R\Gamma(2) \times C\Gamma(1)$ be the scaled doddering tree and the scaled gluer tree associated with $q$ under $\mathbb{P}_D^n$ as defined in Section 3.3. A corner $x \in [0,1)$ is called a *node-corner* if $F_{\mathbf{G}_n}(x) \in E_{\mu_{\mathbf{G}_n}}$. Let $b_n$ be the application that sends the $(k+1)$th node of $\mathbf{D}_n$ according to the RO on the $k$th node-corner of $\mathbf{G}_n$ according to the CO (see Figure 17). The root of $\mathbf{D}_n$ has no image by $b_n$. The application $b_n$ satisfies the three conditions of Definition 3.21, and then $(\mathbf{D}_n, \mathbf{G}_n, b_n)$ is 1-admissible. The rescaled quadrangulation under $\mathbb{P}_D^n$ is the $\mathfrak{M}(1)$-valued random variable:

$$\mathbf{M}_n = \mathrm{Map}(\mathbf{D}_n, \mathbf{G}_n, b_n).$$

We denote by $\mathfrak{M}_n$ the support of the random variable $\mathbf{M}_n$ in $\mathfrak{M}(1)$. We identify $\mathfrak{M}_n$ with the set of rescaled quadrangulations with $n$ faces. This identification is allowed because the function $b_n$ and the measures $\mu_{\mathbf{D}_n}$ and $\mu_{\mathbf{G}_n}$ are completely determined by $(r_n^+, v_n^+)$.

*Brownian map.* We define the Brownian map with the help of the limit trees $(\mathbf{D}_\infty, \mathbf{G}_\infty)$ and the application $b_\infty$ defined by

(3.11) $$b_\infty(\dot{u}) = \rho(u)/2,$$

where

$$\rho : [0,2] \longrightarrow [0,2],$$
$$u \longmapsto \inf\Big\{x, x \underset{\mathbf{c}^+}{\sim} u\Big\}.$$

LEMMA 3.31. *The 3-tuple $(\mathbf{D}_\infty, \mathbf{G}_\infty, b_\infty)$ is 1-admissible.*

PROOF. We have to check the three conditions of Definition 3.21.

(i) is a direct consequence of (3.11).
(ii) follows from the fact that $\rho$ is increasing.
(iii) if $b_\infty(u) \underset{\mathbf{v}^+}{\sim} b_\infty(w)$, then $\rho(u)/2 \underset{\mathbf{v}^+}{\sim} \rho(w)/2$. Since $(\mathbf{r}^+, \mathbf{v}^+)$ is in $\mathbb{T}$, we have $\mathbf{r}^+(\rho(u)/2) = \mathbf{r}^+(\rho(w)/2)$; this implies that $\mathbf{c}^+(\rho(u)) = \mathbf{c}^+(\rho(w))$ and $\mathbf{c}^+(u) = \mathbf{c}^+(w)$. □

DEFINITION 3.32. We call the Brownian map the $\mathfrak{M}(1)$ valued random variable $\mathbf{M}_\infty = \mathrm{Map}(\mathbf{D}_\infty, \mathbf{G}_\infty, b_\infty)$.

*Convergence.* Finally, here is the convergence result:

THEOREM 3.33. *The following weak convergence holds in $(\mathfrak{M}(1), d_{\mathfrak{M}(1)})$:*

$$\mathbf{M}_n \xrightarrow[n]{weakly} \mathbf{M}_\infty.$$

PROOF. The nodes of $\mathbf{D}_n$ but the root are sent on different corners $\{k/(2n), k \in [\![0, 2n-1]\!]\}$ of $\mathbf{G}_n$, and then

$$(3.12) \quad C_{\mu_{\mathbf{D}_n} \circ b_n^{-1}}(x) = \begin{cases} \dfrac{2}{2n+1}\lfloor 2nx+1 \rfloor, & \text{on } [0, (2n-1)/(2n)], \\ 0, & \text{for } x < 0, \\ 4n/(2n+1), & \text{for } x > (2n-1)/(2n). \end{cases}$$

In other respects, we have $\mathbb{P}^+$ a.s.,

$$(3.13) \qquad C_{\mu_{\mathbf{D}_\infty} \circ b_\infty^{-1}} = C_{\mathrm{Leb}^{[0,2]} \circ b_\infty^{-1}} = C_{2\,\mathrm{Leb}^{[0,1]}}.$$

Indeed, the Lebesgue's measure of $L_{\mathbf{D}_\infty} = \{y \in [0,2], \# \deg F_{\mathbf{D}_\infty}(y) = 1\}$, that is, the set of corners of leaves in $\mathbf{D}_\infty$, is $\mathbb{P}^+$ a.s. equal to 2. For $y = 2x \in L_{\mathbf{D}_\infty}$, one has $C_{\mathrm{Leb}^{[0,2]} \circ b_\infty^{-1}}(x) = y = 2x$. This means that, on a dense subset of $[0,1]$, $C_{\mathrm{Leb}^{[0,2]} \circ b_\infty^{-1}}(x) = 2x$ and then (3.13) holds.

The sequence $(C_{\mu_{\mathbf{D}_n} \circ b_n^{-1}})_n$ is nonrandom and converges to $C_{2\,\mathrm{Leb}^{[0,1]}}$ uniformly. Now, the conclusion follows from (3.13) and Proposition 3.20. □

3.5.2. *Rooted quadrangulations with random edge lengths.* We denote by $\overrightarrow{E_n}$ the set of rooted quadrangulations with $n$ faces, in which the edges have real positive lengths. The set $\overrightarrow{E_n}$ can be represented by

$$\overrightarrow{E_n} \simeq \overrightarrow{\mathcal{Q}}_n \times ((0,+\infty))^{2n}.$$

We assume that the edge lengths are positive i.i.d. random variables, independent from the underlying quadrangulation. We denote by $\mu$ the distribution of the lengths and by $L$ a $\mu$-distributed random variable. Hence, $\overrightarrow{E_n}$ is naturally endowed with the law $\mathbb{P}_D^n \otimes \mu^{2n}$. We assume in the following that $\mathbb{E}(L) = 1$ and that there exists $\delta > 0$ such that $\mathbb{E}(L^{5+\delta}) < +\infty$.

Set $q = (U(q), (l(i))_{i \in [\![0,2n]\!]}) \in \overrightarrow{E_n}$. The quadrangulation $U(q)$ is the underlying rooted quadrangulation that is homeomorph (as a rooted map) to $q$. Our aim is to embed scaled version of elements of $\overrightarrow{E_n}$ in $\mathfrak{M}(1)$ and to show the weak convergence to the Brownian map in this space.

Let $(\mathbf{D}_n, \mathbf{G}_n)$ be the doddering tree and the gluer tree associated with $U(q)$ and let $\mathbf{M}_n = \mathrm{Map}(\mathbf{D}_n, \mathbf{G}_n, b_n)$ be as defined in Section 3.5.1; let $\widetilde{\mathbf{D}}_n$ be the doddering tree with random edge lengths that satisfies:

- $\widetilde{\mathbf{D}}_n$ and $\mathbf{D}_n$ have the same arborescent structure.
- The edge lengths of $\widetilde{\mathbf{D}}_n$ are independent and have the law of $L/n^{1/4}$.
- $\mu_{\widetilde{\mathbf{D}}_n} = \mu_{\mathbf{D}_n}$.

The 3-tuple $(\widetilde{\mathbf{D}}_n, \mathbf{G}_n, b_n)$ is 1-admissible. One sets $\widetilde{\mathbf{M}}_n = \mathrm{Map}(\widetilde{\mathbf{D}}_n, \mathbf{G}_n, b_n)$ the (normalized) abstract map corresponding to $q$.

THEOREM 3.34. *We have*

$$\widetilde{\mathbf{M}}_n \xrightarrow[n]{weakly} \mathbf{M}_\infty \qquad in \ (\mathfrak{M}(1), d_{\mathfrak{M}(1)}).$$

The proof is postponed to the Appendix.

REMARK 3.35. To obtain a continuous limit for the tour of $\widetilde{\mathbf{D}}_n$, the increments $L$ must satisfy:

$$(3.14) \qquad \text{for any fixed } \alpha > 0, \qquad \mathbb{P}\bigg(\sup_{i \in [\![1,4n]\!]} L_i \geq \alpha n^{1/4}\bigg) \longrightarrow 0.$$

If the distribution tail of $L$ is of the form $\mathbb{P}(L > x) \sim cx^{-\beta}$, then $\beta = 4$ appears to be a threshold: if $\beta > 4$, then (3.14) is true, if $\beta \leq 4$, then (3.14) is false (see analogous considerations for the limit of discrete snake in [22]). This means that one may extend our result to the case where the lengths own a moment of order $4 + \varepsilon$.

3.6. *Pointed abstract maps.* Our leading idea to define the notion of pointed abstract map is simple, even if the formal description will be quite involved: a pointed discrete map is an equivalence class of rooted discrete maps. In the case of quadrangulations, the equivalence classes are defined with the help of the canonical surjection $K$ from $\vec{\mathcal{Q}}_n$ onto $\mathcal{Q}_n^\bullet$. We choose here to define pointed abstract map in terms of classes of rooted abstract maps. Since these later are defined with the help of trees, the rerooting on map will be translated on rerooting on trees. A little technical difficulty arises here: the rerooting are not exactly in the same place on the doddering tree and on the gluer tree; this may be seen easily in the discrete case since $b_n$ owns a simple form when one deals with the RHP $r_n^+$ of $D_n$, and is not so simple, when one deals with its RDFW, what we do, to be consistent.

Let $\mathcal{D} = \mathrm{RTree}(\xi_\mathcal{D}, \mu_\mathcal{D})$, $\mathcal{G} = \mathrm{CTree}(\zeta_\mathcal{G}, \mu_\mathcal{G})$ and $M$ the rooted abstract map $M = \mathrm{Map}(\mathcal{D}, \mathcal{G}, b)$. For any corner $\theta \in [0, 2a)$ of the root vertex $\widehat{0}_M$, we will define a rooted map $M^{(\theta)}$ corresponding to $M$ rerooted in $\theta$. For this, we will reroot $\mathcal{D}$ in $\theta$ and define an associated rerooting for $\mathcal{G}$ (recall that the root $\widehat{0}_M$ of $M$ is the root of $\mathcal{D}$). We endow the set $[0, 2a)$ (of the corners of $\mathcal{D}$) with the total order $\leq_\theta$ that matches the usual order on $[0, \theta)$ and $[\theta, 2a)$ and such that, for any $(x, y) \in [\theta, 2a) \times [0, \theta)$, we have $x \leq_\theta y$. The order $\leq_\theta$ induces a total order $\ll_\theta$ on the image of $E_{\mu_\mathcal{D}}$ by $b$. Recall that $b(E_{\mu_\mathcal{D}})$ is included in the set of corners of $\mathcal{G}$, that is, $[0, a)$. We will reroot $\mathcal{G}$ in the "first" corner of $\mathcal{G}$ according to $\ll_\theta$. We set $A_\theta = \{\dot{x} \in E_{\mu_\mathcal{D}}, \underline{\rho}(\dot{x}) > \theta\}$ and

$$\tau_\theta = \begin{cases} \inf\{b(\dot{x}), \dot{x} \in A_\theta\}, & \text{if } A_\theta \neq \varnothing, \\ 0, & \text{if } A_\theta = \varnothing. \end{cases}$$

Consider $\mathcal{D}^{(\theta)} = \mathrm{RTree}(\xi_\mathcal{D}^{(\theta)}, \mu_\mathcal{D}^{(\theta)})$ and $\mathcal{G}^{(\tau_\theta)} = \mathrm{CTree}(\zeta_\mathcal{G}^{(\tau_\theta)}, \mu_\mathcal{G}^{(\tau_\theta)})$, where

$$\mu_\mathcal{D}^{(\theta)}(A) = \mu_\mathcal{D}(A + \theta \bmod 2a)$$

and

$$\mu_\mathcal{G}^{(\tau_\theta)}(B) = \mu_\mathcal{G}(B + \tau_\theta \bmod a).$$

We define the function $b_\theta$ by

$$b_\theta(\dot{x}) = b(\dot{x}_\theta) - \tau_\theta \bmod a,$$

where $\dot{x}_\theta$ is the node in $\mathcal{D}$ with representative $x + \theta \bmod 2a$.

LEMMA 3.36. *$(\mathcal{D}^{(\theta)}, \mathcal{G}^{(\tau_\theta)}, b_\theta)$ is $a$-admissible.*

PROOF. The checking of the conditions (i) and (iii) are simple. To show (ii), we endow the set $[0, a)$ (of the corners of $\mathcal{G}$) with the total order $\leq^{\tau_\theta}$ that matches the usual order on $[0, \tau_\theta)$ and $[\tau_\theta, a)$ and such that, for any $(x,y) \in [\tau_\theta, a) \times [0, \tau_\theta)$, we have $x \leq^{\tau_\theta} y$. Since $\theta$ is a root corner of $\mathcal{D}$, $x \geq \theta \Leftrightarrow \rho(\dot{x}) \geq \theta$. This equivalence implies that $b$ is increasing from $(E_{\mu_\mathcal{D}}, \leq_\theta)$ in $([0,a), \leq^{\tau_\theta})$. We set $\underline{\rho}^\theta(\dot{x})$ the smallest representative of $\dot{x}$ in $\mathcal{D}^{(\theta)}$.

For any $\dot{u}$ and $\dot{v}$ in $\overline{E}_{\mu_{\mathcal{D}^{(\theta)}}} \setminus \{\text{root vertex}\}$,

$$\underline{\rho}^\theta(\dot{u}) \leq \underline{\rho}^\theta(\dot{v}) \Leftrightarrow \rho(\dot{u}_\theta) \leq_\theta \rho(\dot{v}_\theta)$$
$$\Leftrightarrow b(\dot{u}_\theta) \leq^{\tau_\theta} b(\dot{v}_\theta)$$
$$\Leftrightarrow b_\theta(\dot{u}) \leq b_\theta(\dot{v}). \quad \square$$

We set $M^{(\theta)} = \text{Map}(\mathcal{D}^{(\theta)}, \mathcal{G}^{(\tau_\theta)}, b_\theta)$. We introduce an equivalence relation in $\mathfrak{M}$: two maps $M_1$ and $M_2$ are equivalent, we note $M_1 \underset{\bullet}{\simeq} M_2$, if there exists a corner $\theta \in \widehat{0}_{M_1}$, such that $M_2 = M_1^{(\theta)}$. We call the quotient space $\mathfrak{M}^\bullet = \mathfrak{M}/\underset{\bullet}{\simeq}$, the set of pointed abstract maps. We introduce $d_{\mathfrak{M}^\bullet}$:

$$d_{\mathfrak{M}^\bullet}(M_1^\bullet, M_2^\bullet) = \inf\{d_\mathfrak{M}(M_1^{\theta_1}, M_1^{\theta_2}), \ (\theta_1, \theta_2) \in \widehat{0}_{M_1} \times \widehat{0}_{M_2}\}.$$

For any $M^\bullet$ and $\alpha > 0$, we set

$$\mathcal{B}(M^\bullet, \alpha) = \{N^\bullet \in \mathfrak{M}^\bullet, \ d_{\mathfrak{M}^\bullet}(M^\bullet, N^\bullet) < \alpha\}.$$

We endow $\mathfrak{M}^\bullet$ with the topology generated by the family

$$\{\mathcal{B}(M^\bullet, \alpha), M^\bullet \in \mathfrak{M}^\bullet, \alpha > 0\}.$$

In Section 3.5.1, $\mathfrak{M}_n$ was defined as the set of rooted abstract maps corresponding to normalized quadrangulations from $\overrightarrow{\mathcal{Q}}_n$. We denote by $\mathfrak{M}_n^\bullet$ the quotient set $\mathfrak{M}_n/\underset{\bullet}{\simeq}$ and by $\mathbf{M}_n^\bullet$ a random variable uniformly distributed on $\mathfrak{M}_n^\bullet$. The elements of $\mathfrak{M}_n^\bullet$ are identified with normalized quadrangulations from $\mathcal{Q}_n^\bullet$ (see Section 4.1.3). Now, let $\mathfrak{M}^\star$ be the set of maps in $\mathfrak{M}$ with root degree 1, and by $\mathfrak{M}_0 = \mathfrak{M} \setminus \mathfrak{M}^\star$. We can consider that

$$\mathfrak{M}^\bullet = \mathfrak{M}^\star \cup \left(\mathfrak{M}_0/\underset{\bullet}{\simeq}\right).$$

It follows from Lemma 4.15 that the Brownian map has a.s. a root with degree 1, and so it is almost surely in $\mathfrak{M}^\star$; we then consider that the random variable $\mathbf{M}_\infty$ take its values in $\mathfrak{M}^\bullet$. We have the following:

THEOREM 3.37. *The following weak convergence holds in $(\mathfrak{M}^\bullet, d_{\mathfrak{M}^\bullet})$:*

$$\mathbf{M}_n^\bullet \xrightarrow[n]{weakly} \mathbf{M}_\infty.$$

Some additional considerations are needed to present a proof of this result. The proof of this theorem will be given at the end of Section 4.

REMARK 3.38. The weak convergence of the encodings of rooted quadrangulations $(\mathbf{r}_n^+, \mathbf{v}_n^+)$ will be shown via the convergence of the encoding of pointed quadrangulations. The convergence of rooted quadrangulations under $\mathbb{P}_D^n$ is then obtained as a consequence of the convergence of pointed quadrangulations under $\mathbb{P}_U^n$.

**4. Embeddings and convergence.** The aim of this section is first to prove the convergence of the encodings of pointed quadrangulations and then to prove Theorems 3.37 and 3.3.

4.1. *The space of pointed quadrangulations.* We present a combinatorial fact concerning pointed quadrangulations that may allow to better understand our definition of pointed abstract maps.

4.1.1. *Effect of the starting point in the construction of $Q(T)$.* Consider a well-labeled tree $T$ encoding by $(R_n^+, V_n^+)$, and let $i_1 < \cdots < i_k$ be the times such that $R_n^+(i_l) = 1$. The starting point of the construction of $Q(T)$ in Section 2.3.2 is the corner $i_1 = 0$; in other words, we began the CDFT by the root-edge of $T$. Let us examine what happens if we start the construction from the corner $i_l$ and then visit the times $i_l + 1, \ldots, 2n-1, 0, \ldots, i_l - 1$. The construction is the same in each integer's interval $[\![i_s, i_{s+1}[\![$. The only change is that the chords $\widehat{(i_l, v)}$ are not drawn in the same order: they are drawn according to their ranks in the circular permutation $(i_l, \ldots, 2n-1, 0, \ldots, i_l - 1)$ of $(0, \ldots, 2n-1)$. Thus, we obtain a rooted map $\overrightarrow{q_l}$. The only difference with $Q(T) = \overrightarrow{q_1}$ is that the adjacent edges of $v$ are circularly permuted; this means that $\overrightarrow{q_l}$ is identical to $Q(T)$ as unrooted map. The root-edge of $\overrightarrow{q_l}$ is $\widehat{(v, i_l)}$ instead of $\widehat{(v, i_1)}$ for $\overrightarrow{q_1}$. In other words, the algorithm starting from any time $i_l$, $l = 1, \ldots, k$, gives the same pointed quadrangulation (see Figure 18).

It turns out that the set of rooted quadrangulations $\{\overrightarrow{q_l}, l \in \{1, \ldots, k\}\}$ is exactly $K^{-1}(Q(T))$. Moreover, the map $\overrightarrow{q_l}$ is equal to $Q(T_l)$, where $T_l$ is obtained by rerooting $T$ in $(i_l, i_l + 1)$ and in keeping its labels.

We formalize now the notion of rerooting of trees and labeled trees.

4.1.2. *Rerooting in $\mathcal{W}_n$.* Let $t$ be an element of $\Omega_n$ and let $F$ be its CDFT. For each $\theta \in [\![0, 2n]\!]$, we define an application $t \mapsto t^{(\theta)}$ from $\Omega_n$ to $\Omega_n$, that we call rerooting; the rooted tree $t^{(\theta)}$ is identical as an unrooted tree to $t$, and the root-edge of $t^{(\theta)}$ is $\overrightarrow{F(\theta)F(\theta+1)}$.

Let $G_n$ be the group $([\![0, 2n[\![, \oplus)$, where $\oplus$ is the addition modulo $2n$. Consider $F^{(\theta)}$ and $V_n^{(\theta)}$ the CDFT and the CDFW of $t^{(\theta)}$. The function $F^{(\theta)}$ visits successively the nodes $F(\theta), F(\theta+1), \ldots, F(2n-1), F(0), \ldots, F(\theta-1)$, and then, it is straightforward that

$$F^{(\theta)}(x) = F(\theta \oplus x) \qquad \text{for any } x \in [\![0, 2n[\![$$

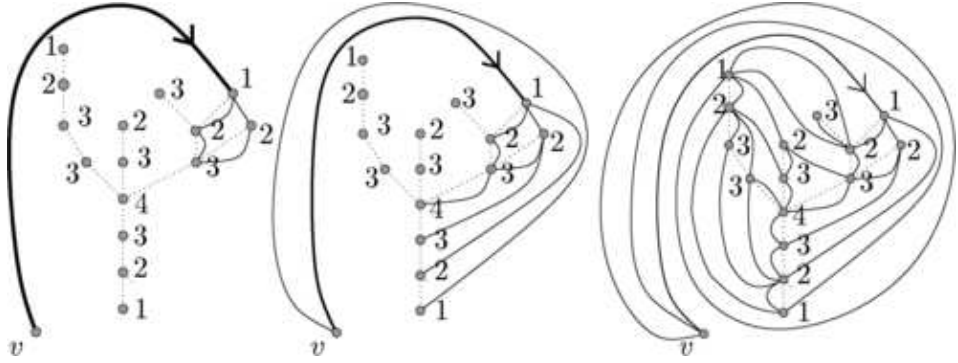

FIG. 18. *Construction of $Q(T)$ with the procedure starting with the third corner with label* 1.

and that the rerooting operation is an action of $G_n$ on $\Omega_n$ which lets invariant the unrooted tree structure (i.e., $t = t^{(\theta)}$ as unrooted trees). More precisely, the quotient set $\Omega_n/G_n$ is in bijection with the set of unrooted trees with $n$ edges. Moreover, $V_n$ and $V_n^{(\theta)}$ are related by

(4.1) $V_n^{(\theta)}(i) = V_n(\theta \oplus i) + V_n(\theta) - 2\check{V}_n(\theta \oplus i, \theta)$ for any $i \in [\![0, 2n]\!]$.

Indeed, since $V_n^{(\theta)}(i) = d(F^{(\theta)}(i), F^{(\theta)}(0)) = d(F(\theta \oplus i), F(\theta))$, the distance in the tree between the nodes $F(\theta \oplus i)$ and $F(\theta)$ is $V_n(\theta \oplus i) + V_n(\theta) - 2\check{V}_n(\theta \oplus i, \theta)$. We now extend the action of $G_n$ on $\mathcal{W}_n$ by defining the action on the labels. We set

(4.2) $\qquad R_n^{[\theta]}(i) = R_n(\theta \oplus i) - R_n(\theta) + 1$ for any $i \in [\![0, 2n[\![$.

This defines clearly an action of $G_n$ on $\mathcal{W}_n$.

LEMMA 4.1. *If $(R_n, V_n)$ is the encoding of a labeled tree $T \in \mathcal{W}_n$, then for any $\theta \in [\![0, 2n]\!]$, $(R_n^{[\theta]}, V_n^{(\theta)})$ is the encoding of a labeled tree belonging to $\mathcal{W}_n$, which we denote by $T^{(\theta)}$.*

PROOF. If a tree $t \in \Omega_n$ has CDFW $V_n$, then $V_n^{(\theta)}$ is the CDFW of $t^{(\theta)}$. It remains to prove that if $R_n$ encodes a labeling of $t$, then $R_n^{[\theta]}$ encodes a labeling of $t^{(\theta)}$. At first, note that $R_n^{[\theta]}(0) = 1$ and that $R_n^{[\theta]}(i+1) - R_n^{[\theta]}(i)$ is equal to $+1, -1$ or $0$. We have to show that if $i$ and $j$ are such that $F^{(\theta)}(i) = F^{(\theta)}(j)$, then $R_n^{[\theta]}(i) = R_n^{[\theta]}(j)$. If $F^{(\theta)}(i) = F^{(\theta)}(j)$, then $F(\theta \oplus i) = F(\theta \oplus j)$ and thus, $R_n(\theta \oplus i) = R_n(\theta \oplus j)$. □

4.1.3. *Normalized pointed quadrangulations seen as pointed abstract maps.* We first reinterpret the content of Section 4.1.1 in terms of rerooting.

According to Section 4.1.1, we have the following result:

PROPOSITION 4.2. *Let $\overrightarrow{q_1}$, $\overrightarrow{q_2}$ in $\overrightarrow{\mathcal{Q}}_n$ and $T_1 = Q^{-1}(\overrightarrow{q_1})$, $T_2 = Q^{-1}(\overrightarrow{q_2})$ the associated well-labeled tree in $\mathcal{W}_n^+$. We have*

$$K(\overrightarrow{q_1}) = K(\overrightarrow{q_2}) \quad \Leftrightarrow \quad \text{there exists } \tau \in [\![0, 2n]\!] \text{ such that } T_1^{(\tau)} = T_2.$$

REMARK 4.3. (1) Let $T \in \mathcal{W}_n^+$ be a well-labeled tree encoded by $(R_n^+, V_n^+)$. The labeled tree $T^{(\tau)}$ is well-labeled iff $R_n^+(\tau) = \min R_n^+ = 1$.

(2) Consider two well-labeled trees $T$ and $T'$ in $\mathcal{W}_n^+$ equal up to a rerooting, let $(R_n^+, V_n^+)$ and $(R_n^{+\prime}, V_n^{+\prime})$ be their encodings, and $(\mathcal{D}(R_n^+), \mathcal{G}(V_n^+))$ and $(\mathcal{D}(R_n^{+\prime}), \mathcal{G}(V_n^{+\prime}))$ the associated pairs of trees [note $C_n^+$ and $C_n^{+\prime}$ the RDFW of $\mathcal{D}(R_n^+)$ and $\mathcal{D}(R_n^{+\prime})$].

There exists a constant $\tau$ such that $(R_n^{+\prime}, V_n^{+\prime}) = (R_n^{[\tau]}, V_n^{(\tau)})$. Since $R_n^+(\tau) = 1$, we have $R_n^+(\cdot) = R_n^{+\prime}(\cdot + \tau \mod 2n)$. Since, up to the point added to encode the root, $R_n^+$ and $R_n^{+\prime}$ are the RHP of $\mathcal{D}(R_n^+)$ and $\mathcal{D}(R_n^{+\prime})$, it is immediate that $\mathcal{D}(R_n^+)$ and $\mathcal{D}(R_n^{+\prime})$ are equal as unrooted trees, have the same root-vertex, and that only the second extremities of their root-edges are different. As a consequence, the two RDFW $C_n^+$ and $C_n^{+\prime}$ satisfy $C_n^+(\cdot) = C_n^{+\prime}(\cdot + \theta \mod 4n)$, where $\theta \in [\![0, 4n]\!]$, is a corner of the root of $\mathcal{D}(R_n^{+\prime})$.

PROPOSITION 4.4. *Let $M_1 = \mathrm{Map}(D_n, G_n, b_n)$ and $M_2 = \mathrm{Map}(D_n', G_n', b_n')$ be two elements of $\mathfrak{M}_n$ corresponding to two rooted quadrangulations $\overrightarrow{q_1}$ and $\overrightarrow{q_2}$ in $\overrightarrow{\mathcal{Q}}_n$. We have*

$$K(\overrightarrow{q_1}) = K(\overrightarrow{q_2}) \Leftrightarrow M_1 \underset{\bullet}{\simeq} M_2.$$

This proposition allows to identify (the normalized quadrangulations from) $\mathcal{Q}_n^\bullet$ with $\mathfrak{M}_n^\bullet$. Indeed, for an element $M = \mathrm{Map}(D, G, b) \in \mathfrak{M}_n$, the injection $b$, the measures $\mu_\mathcal{D}$ and $\mu_\mathcal{G}$ are fixed knowing the doddering tree $\xi_\mathcal{D}$ and the gluer tree $\zeta_\mathcal{G}$.

PROOF OF PROPOSITION 4.4. We first prove the implication. This is mainly a consequence of Proposition 4.2 and of Remark 4.3. Let $M_1 = \mathrm{Map}(D_n, G_n, b_n)$ and $M_2 = \mathrm{Map}(D_n', G_n', b_n')$, corresponding to the same pointed quadrangulation. The well-labeled trees $T_1$ and $T_2$ are equal up to a rerooting, and this is also the case for $D_n$ and $D_n'$. Let $\theta$ be the unique real number such that $D_n^{(\theta)} = D_n'$. According to the above remark, $\theta$ is a corner of the root of $D_n$. In some sense $\theta$ is used as a shift to pass from $D_n$ to $D_n'$. The root of $D_n$ has no image by $b_n$. Hence, $b_n(\dot{\theta})$ does not exist, and then cannot be used as a shift to define $b_n'$ or $G_n'$. The good shift is $\tau_\theta = \inf\{b(\dot{x}), \rho(\dot{x}) > \theta\}$, as defined in Section 3.6. In other words, $\tau_\theta$ is the image of the first node of $D_n$ visited after the corner $\theta$ in $D_n$ (which corresponds to the corner with label 1 in the well-labeled tree). But, in the discrete case, the first node is visited at time $u := \theta + 1/2n$ (we work on normalized version), and then in the discrete case $\tau_\theta = b(\dot u)$. Hence, $G'_n = G_n^{(\tau_\theta)}$, and $b'_n$ is equal to $b_\theta$ as defined in Section 3.6.

To conclude that $M_1 \underset{\bullet}{\simeq} M_2$, it remains to consider the measures associated with the corners of the trees. The measures were defined to be invariant by the changing of roots involved here: the measure on $D_n$ is invariant by the changing of root on the corners of the root, the measure on $G_n$ is invariant by any change of root on a node-corner.

Conversely, starting from $M_1 \underset{\bullet}{\simeq} M_2$, we deduce, using similar arguments, that $(D_n, G_n)$ and $(D'_n, G'_n)$ correspond to two well-labeled trees equal up to a change of root. □

4.2. *Description of pointed quadrangulations with labeled trees.* Hence, rooted quadrangulations are defined with well-labeled trees and pointed quadrangulations with classes of well-labeled trees. Well-labeled trees appear as labeled trees conditioned to be positive, and this conditioning is quite difficult to handle. In the present section we introduce some combinatorial facts that will allows to work with labeled trees.

The rerooting operation defined by (4.1) and (4.2) is an action on $\mathcal{W}_n$; for $T$ and $T'$ in $\mathcal{W}_n$, we write $T = T' \bmod G_n$, if $T^{(\theta)} = T'$ for some $\theta \in [\![0, 2n]\!]$.

Let $T \in \mathcal{W}_n$, $F$ be its CDFT, $R_n$ be its label process, and $\theta_1 < \cdots < \theta_k$ be the times in $[\![0, 2n-1]\!]$ such that $R_n(\theta_l) = \inf_{s \in [\![0, 2n-1]\!]} R_n(s)$ [notice that $R_n(\theta_1)$ may be different to 1 since $T$ is not assumed to be well-labeled]. Consider the well-labeled tree

$$L(T) = T^{(\theta_1)}.$$

By formula (4.2), $R^{(\theta_1)}$ is positive and $T^{(\theta_1)}$ is an element of $\mathcal{W}_n^+$. Moreover, $T \in \mathcal{W}_n^+$ clearly implies that $L(T) = T$; we thus have defined a surjective application $L : \mathcal{W}_n \to \mathcal{W}_n^+$.

4.2.1. *$\mathcal{Q}_n^\bullet$ is in bijection with $\mathcal{W}_n/G_n$.* The application $\widetilde{Q}$ defined by

$$\widetilde{Q}(T) = K(Q(L(T))) \tag{4.3}$$

is a surjection from $\mathcal{W}_n$ onto $\mathcal{Q}_n^\bullet$.

The aim of this part is to show the following theorem:

THEOREM 4.5. *The application* $\Phi : \mathcal{W}_n/G_n \longrightarrow \mathcal{Q}_n^\bullet$, *defined by* $\Phi(x) = \widetilde{Q}(T)$, *where $T$ is any representative of $x$, is well defined and bijective.*

PROOF. Notice that, for any $T \in \mathcal{W}_n$,

$$L(T) = T \bmod G_n. \tag{4.4}$$

Now set $\theta \in [\![0, 2n-1]\!]$ and $\theta_i$ such that $\theta_i \leq \theta < \theta_{i+1}$; since $R^{[\theta]}(s) = R(\theta \oplus s) - R(\theta) + 1$, the successive times $s_1 < \cdots < s_k$ in $[\![0, 2n-1]\!]$ such that $R^{[\theta]}(s_l) = \inf_{s \in [\![0, 2n-1]\!]} R^{[\theta]}(s)$ are $\theta_{i+1} - \theta, \ldots, \theta_k - \theta, 2n - 1 + \theta_1 - \theta, \ldots, 2n - 1 + \theta_i - \theta$. Thus,

$$L(T^{(\theta)}) = T^{(\theta_{i+1})}. \tag{4.5}$$

It follows from (4.5) and Section 4.1.1 that, for any $T \in \mathcal{W}_n$ and any $\theta \in [\![0, 2n-1]\!]$,

$$\widetilde{Q}(T^{(\theta)}) = \widetilde{Q}(T). \tag{4.6}$$

Moreover, let $T$ and $T'$ in $\mathcal{W}_n^+$; if $K \circ Q(T) = K \circ Q(T')$, then the rooted maps $Q(T)$ and $Q(T')$ are equal as pointed maps, the only difference is that the root-edge is possibly not the same but the root-vertex (say, $v$) is the same; let $F$ be the CDFT of $T$, the root of $Q(T)$ is $\overrightarrow{vF(0)}$ and the root of $Q(T')$ is $\overrightarrow{vF(\theta)}$ for some $\theta \in [\![0, 2n-1]\!]$. Thus, $Q(T') = Q(T^{(\theta)})$ and, since $Q$ is a bijection, $T' = T^{(\theta)}$. We have thus proved that, for $T$ and $T'$ in $\mathcal{W}_n^+$, $\quad K \circ Q(T) = K \circ Q(T') \Rightarrow T = T' \bmod G_n$.

Taking into account (4.5), it follows that,

for $T$ and $T'$ in $\mathcal{W}_n$, $\quad \widetilde{Q}(T) = \widetilde{Q}(T') \Rightarrow T = T' \bmod G_n$.

This together with (4.6) means,

for $T$ and $T'$ in $\mathcal{W}_n$, $\quad \widetilde{Q}(T) = \widetilde{Q}(T') \Leftrightarrow T = T' \bmod G_n$.

Thus, the application $\Phi : \mathcal{W}_n/G_n \longrightarrow \mathcal{Q}_n^\bullet$, defined by $\Phi(x) = \widetilde{Q}(T)$, where $T$ is any representative of $x$, is well defined and is injective. Since $\widetilde{Q}$ is surjective, $\Phi$ is bijective. $\square$

REMARK 4.6. A bijection between $\mathcal{Q}_n^\bullet$ and a set of unrooted marked trees is presented in the Appendix.

4.2.2. *Elements on the distribution of pointed quadrangulations.* We endow $\mathcal{W}_n$ with the uniform law, denoted by $\mathbb{P}_{\mathcal{S}}^n$. The image law of $\mathbb{P}_{\mathcal{S}}^n$ by the surjection $\widetilde{Q}$ is the distribution on $\mathcal{Q}_n^\bullet$ denoted by $\mathbb{P}_{\widetilde{\mathcal{S}}}^n$, and defined for any $q \in \mathcal{Q}_n^\bullet$ by

$$\mathbb{P}_{\widetilde{\mathcal{S}}}^n(q) = \frac{\#\widetilde{Q}^{-1}(q)}{\#\mathcal{W}_n} = \frac{\#\widetilde{Q}^{-1}(q)}{C_n 3^n}. \tag{4.7}$$

The proof of the following proposition is given in the Appendix.

PROPOSITION 4.7. *The distance in variation between $\mathbb{P}_U^n$ and $\mathbb{P}_{\widetilde{\mathcal{S}}}^n$ goes to 0 when $n$ goes to $+\infty$.*

*Consequence.* The two models of pointed quadrangulations $(\mathcal{Q}_n^\bullet, \mathbb{P}_U^n)$ and $(\mathcal{Q}_n^\bullet, \mathbb{P}_{\tilde{\mathcal{S}}}^n)$ are asymptotically equivalent. In the sequel we will only consider the distribution $\mathbb{P}_{\tilde{\mathcal{S}}}^n$.

4.3. *Embedding in functional spaces and weak convergences.* With each element of $\mathcal{W}_n$ is associated one and only one normalized encoding $(r_n, v_n)$. By construction, $(r_n, v_n)$ is an element of $\mathbb{T}$ (it is called in [29] the normalized tour of the discrete snake).

The set of all normalized encodings $(r_n, v_n)$ of elements of $\mathcal{W}_n$ is denoted by $\mathbb{T}_n$. Since $\mathcal{W}_n \simeq \mathbb{T}_n$, we say that $\mathcal{W}_n$ is embedded in $\mathbb{T}$ and that its embedding is $\mathbb{T}_n$. The set $\mathbb{T}_n$ can be characterized as follows: $(f, \zeta) \in \mathbb{T}_n$ if and only if $(f, \zeta) \in \mathbb{T}$, and for $i \in [\![0, 2n-1]\!]$,

$$(4.8) \quad \begin{cases} f \text{ and } \zeta \text{ are linear in the intervals } \left[\dfrac{i}{2n}, \dfrac{i+1}{2n}\right], \\ \zeta((i+1)/(2n)) - \zeta(i/(2n)) = \pm n^{-1/2}, \\ f((i+1)/(2n)) - f(i/(2n)) \in \{0, \pm n^{-1/4}\}. \end{cases}$$

*Notation.* Since $\mathcal{W}_n$ is in bijection with $\mathbb{T}_n$, the image law of the uniform distribution on $\mathcal{W}_n$ is the uniform law on $\mathbb{T}_n$. By a slight abuse of notation, we denote by $\mathbb{P}_{\mathcal{S}}^n$ the uniform law on $\mathbb{T}_n$ (it is also a law on $\mathbb{T}$). We denote by $(\mathbf{r}_n, \mathbf{v}_n)$ a $\mathbb{P}_{\mathcal{S}}^n$-distributed random variable.

PROPOSITION 4.8. *The following weak convergence holds in $\mathbb{T}$*

$$(\mathbf{r}_n, \mathbf{v}_n) \xrightarrow[n]{weakly} (\mathbf{r}, \mathbf{v}).$$

PROOF. The elements of $\mathcal{W}_n$ may be seen as finite branching random walks: their underlying trees being chosen uniformly at random in $\Omega_n$ and the displacements are i.i.d., uniform on $\{-1, 0, +1\}$, independent of the underlying tree, and the value of the root is 1. Then Proposition 4.8 is equivalent to the weak convergence of the normalized tour of the discrete snake to the tour of the Brownian snake. For a proof of this convergence, we refer to [29]; see also [15] for a proof in the Skorohod topology and moments convergence, and [21, 22] for extensions. □

Let us mention a consequence:

PROPOSITION 4.9. *The law of the Brownian snake with lifetime process $\mathbf{e}$ is invariant by rerooting:*

$$\text{for any } \theta \in [0, 1], \qquad \mathbb{P}_\mathcal{S} \circ (J^{(\theta)})^{-1} = \mathbb{P}_\mathcal{S}.$$

PROOF. Let $\theta \in [0,1]$ and $\theta_n = \lfloor 2n\theta \rfloor / 2n$; since $J^{(\theta_n)}$ is a bijection in $\mathbb{T}_n$, $(\mathbf{r}_n, \mathbf{v}_n)$ and $(\mathbf{r}_n, \mathbf{v}_n)^{(\theta_n)}$ have the same law in $\mathbb{T}_n$, and thus, by Proposition 4.8, $(\mathbf{r}_n, \mathbf{v}_n)^{(\theta_n)} \xrightarrow[n]{weakly} (\mathbf{r}, \mathbf{v})$. Moreover, since $J^{(\theta)}$ is continuous, $(\mathbf{r}_n, \mathbf{v}_n)^{(\theta)} \xrightarrow[n]{weakly} (\mathbf{r}, \mathbf{v})^{(\theta)}$. We conclude since $d_\mathbb{T}((\mathbf{r}_n, \mathbf{v}_n)^{(\theta_n)}, (\mathbf{r}_n, \mathbf{v}_n)^{(\theta)}) \longrightarrow 0$. □

Aldous [1], page 40, proves that if $U$ is uniform random on $[0, 1]$ and independent of $\mathbf{v}$, then $\mathbf{v}^{(U)} \stackrel{(d)}{=} \mathbf{v}$. Proposition 4.9 allows us to consider the case $U$ nonuniformly distributed and concerns also the label process.

4.3.1. *Embedding of $\mathcal{Q}_n^\bullet$ in a quotient space of $\mathbb{T}$ and convergence.* We have encoded the normalized rooted quadrangulations by elements of $\mathbb{T}$. We now encode pointed quadrangulations by classes of elements of $\mathbb{T}$.

Since $\mathcal{W}_n \simeq \mathbb{T}_n$, the action of $G_n$ on $\mathcal{W}_n$ defines an action on $\mathbb{T}_n$. More precisely, let $\mathcal{O} = ([0,1), \oplus)$, where $\oplus$ is the addition modulo 1 and let $(\mathcal{O}_n, \oplus)$ be the cyclic subgroup generated by $1/(2n)$ [i.e., the set $\{0, (2n)^{-1}, \ldots, (2n-1)(2n)^{-1}\}$]. It is straightforward to see that the action of $G_n$ on $\mathcal{W}_n$ induces an action of $\mathcal{O}_n$ on $\mathbb{T}_n$ which is just the restriction of the action (to $\mathcal{O}_n$ on $\mathbb{T}_n$) of $\mathcal{O}$ on $\mathbb{T}$ defined in (3.1) and (3.2): indeed, if $(f, \zeta) \in \mathbb{T}_n$ and $\theta \in \mathcal{O}_n$ then $J^{(\theta)}(f, \zeta) \in \mathbb{T}_n$. Since $\mathcal{Q}_n^\bullet \simeq \mathcal{W}_n / G_n$, we have

(4.9) $$\mathcal{Q}_n^\bullet \simeq \mathbb{T}_n / \mathcal{O}_n.$$

This bijection defines an encoding of elements of $\mathcal{Q}_n^\bullet$ by elements of $\mathbb{T}_n / \mathcal{O}_n$. We now embed $\mathbb{T}_n / \mathcal{O}_n$ in $\mathbb{T}/\mathcal{O}$ as follows. Since $\mathcal{O}_n$ is a subgroup of $\mathcal{O}$, the canonical surjection

$$\pi : \mathbb{T} \longrightarrow \mathbb{T}/\mathcal{O}, \qquad x \longmapsto \pi(x) = \bar{x}$$

leads to a natural application

$$\mathcal{I} : \mathbb{T}_n / \mathcal{O}_n \longrightarrow \mathbb{T}/\mathcal{O},$$

$$x^\bullet \longmapsto \bar{x}.$$

Now, if $x$ and $y$ are two elements of $\mathbb{T}_n$ such that $y = x^{(\theta)}$, then necessarily $\theta \in \mathcal{O}_n$. Thus, $\bar{x} = \bar{y} \Rightarrow x^\bullet = y^\bullet$, that is, $\mathcal{I}$ is an injection and $\mathbb{T}_n / \mathcal{O}_n$ is identified with the subset $\mathbb{T}_n / \mathcal{O}$. Thus, in view of (4.9),

$$\mathcal{Q}_n^\bullet \text{ is embedded in } \mathbb{T}/\mathcal{O}.$$

We endow $\mathbb{T}/\mathcal{O}$ with the quotient topology. The law $\mathbb{P}_\mathcal{S}^n$ on $\mathcal{Q}_n^\bullet$ is transported on $\mathbb{T}/\mathcal{O}$; we still call it by the same name:

(4.10) $$\mathbb{P}_{\bar{\mathcal{S}}}^n = \mathbb{P}_\mathcal{S}^n \circ \pi^{-1}.$$

This is the distribution of $\pi(\mathbf{r}_n, \mathbf{v}_n) = \overline{(\mathbf{r}_n, \mathbf{v}_n)}$. Its support is $\mathbb{T}_n/\mathcal{O}$.

Denote by $\mathbb{P}_{\bar{\mathcal{S}}}$ the distribution on $\mathbb{T}/\mathcal{O}$ defined by $\mathbb{P}_\mathcal{S} \circ \pi^{-1}$, that is, the law of $\pi(\mathbf{r}, \mathbf{v}) = \overline{(\mathbf{r}, \mathbf{v})}$. Since $\pi$ is continuous, Proposition 4.8 gives the following:

THEOREM 4.10. *The following weak convergence holds:*
$$\pi(\mathbf{r}_n, \mathbf{v}_n) \xrightarrow[n]{weakly} \pi(\mathbf{r}, \mathbf{v}) \qquad \text{that is, } \mathbb{P}_{\bar{\mathcal{S}}}^n \xrightarrow[n]{weakly} \mathbb{P}_{\bar{\mathcal{S}}}.$$

4.4. *Topology and metric on* $\mathbb{T}/\mathcal{O}$. The space $\mathbb{T}/\mathcal{O}$ is endowed with the quotient topology and the canonical surjection
$$\pi : \mathbb{T} \longrightarrow \mathbb{T}/\mathcal{O}$$
is continuous for this topology. Set
$$\delta : (\mathbb{T}/\mathcal{O})^2 \longrightarrow \mathbb{R}^+,$$
$$(\bar{x}, \bar{y}) \longmapsto \inf_{\substack{\theta_1 \in \mathcal{O} \\ \theta_2 \in \mathcal{O}}} d_{\mathbb{T}}(x^{(\theta_1)}, x^{(\theta_2)})$$

and

(4.11)
$$D : (\mathbb{T}/\mathcal{O})^2 \longrightarrow \mathbb{R}^+,$$
$$(\bar{x}, \bar{y}) \longmapsto \inf_{p \in \mathbb{N}} \inf_{(\bar{z}_i)_{i \in [\![0,p]\!]}} \sum_{i=0}^{p-1} \delta(\bar{z}_i, \bar{z}_{i+1}),$$

where the second infimum is taken on all sequence $\bar{z}_0, \ldots, \bar{z}_p \in \mathbb{T}/\mathcal{O}$ such that $\bar{z}_0 = \bar{x}$ and $\bar{z}_p = \bar{y}$.

PROPOSITION 4.11. *$D$ is a metric on $\mathbb{T}/\mathcal{O}$ inducing the quotient topology.*

For the proof of Proposition 4.11, we shall need the following straightforward inequality, valid for any $\theta, \theta' \in \mathcal{O}$:

(4.12) $$d_{\mathbb{T}}(x^{(\theta)}, y^{(\theta)}) \leq 4 d_{\mathbb{T}}(x^{(\theta')}, y^{(\theta')}).$$

The proof of the following lemma is postponed to the Appendix.

LEMMA 4.12. *$\delta$ enjoys the four following properties:*

(i) $\delta(\bar{x}, \bar{y}) = 0 \iff \bar{x} = \bar{y}$.
(ii) *For any $\bar{x}, \bar{y}, \bar{z} \in \mathbb{T}/\mathcal{O}$,*
$$\delta(\bar{x}, \bar{z}) \leq \delta(\bar{x}, \bar{y}) + 4\delta(\bar{y}, \bar{z}).$$

(iii) *For any $\rho > 0$ and $\bar{x} \in \mathbb{T}/\mathcal{O}$, set $B_\delta(\bar{x}, \rho) = \{\bar{y}, \delta(\bar{x}, \bar{y}) < \rho\}$. The family $\{B_\delta(\bar{x}, \rho), \rho > 0, \bar{x} \in \mathbb{T}/\mathcal{O}\}$ is a base for the quotient topology. We say that $\delta$ induces the quotient topology on $\mathbb{T}/\mathcal{O}$.*

(iv) *For any $\bar{x}, \bar{y} \in \mathbb{T}/\mathcal{O}$,*

$$\delta(\bar{x}, \bar{y})/4 \leq D(\bar{x}, \bar{y}) \leq \delta(\bar{x}, \bar{y}).$$

PROOF OF PROPOSITION 4.11. From its definition, $D$ is obviously a pseudometric, and by Lemma 4.12(i) and (iv), it is a metric. By (4.11), and Lemma 4.12(iv), $\delta/4 \leq D \leq \delta$ and thus, $D$ induces the same topology as $\delta$, that is, the quotient topology on $\mathbb{T}/\mathcal{O}$ by Lemma 4.12(iii). □

4.5. *Embedding in the space of "positive snakes."* Denote by $\mathbb{T}^+$ (resp. $\mathbb{T}^{+\star}$) the nonempty subset of $\mathbb{T}$ of elements $(f, \zeta)$ that satisfy $f(x) \geq 0$ [resp. $f(x) > 0$] for all $x \in (0, 1)$. We denote by $\overline{(f, \zeta)}^+$ the set of nonnegative representatives of $\overline{(f, \zeta)}$:

$$(4.13) \quad \overline{(f,\zeta)}^+ = \overline{(f,\zeta)} \cap \mathbb{T}^+ = \{J^{(\theta)}(f,\zeta),\ \theta \in A(f,\zeta)\},$$

where

$$A(f, \zeta) = \{\theta \in [0, 1),\ f(\theta) = \min f\}$$

is the subset of $[0, 1)$ where $f$ reaches its minimum. We denote by $\mathbb{T}^+/\mathcal{O}$ the set $\{\overline{(f,\zeta)}^+, \overline{(f,\zeta)} \in \mathbb{T}/\mathcal{O}\}$. The set $\mathbb{T}^+/\mathcal{O}$ is the quotient set of $\mathbb{T}^+$ by the equivalence relation $\underset{+}{\sim}$ between elements of $\mathbb{T}^+$:

$$(f_1, \zeta_1) \underset{+}{\sim} (f_2, \zeta_2) \Leftrightarrow \exists \theta \in \mathcal{O} \text{ s.t. } (f_1, \zeta_1) = J^{(\theta)}(f_2, \zeta_2).$$

The application $(f, \zeta) \to \overline{(f, \zeta)}^+$ is surjective from $\mathbb{T}$ on $\mathbb{T}^+/\mathcal{O}$ and

$$\overline{(f_1, \zeta_1)}^+ = \overline{(f_2, \zeta_2)}^+ \Leftrightarrow \exists \theta \in \mathcal{O} \text{ s.t. } (f_1, \zeta_1) = J^{(\theta)}(f_2, \zeta_2).$$

We thus have defined a bijection:

$$\mathrm{Proj} : \mathbb{T}/\mathcal{O} \longrightarrow \mathbb{T}^+/\mathcal{O},$$
$$\overline{(f, \zeta)} \longmapsto \overline{(f, \zeta)}^+.$$

This bijection transports the metric of $\mathbb{T}/\mathcal{O}$ on $\mathbb{T}^+/\mathcal{O}$.

4.5.1. *Topology of $\mathbb{T}^+/\mathcal{O}$.* The image topology by Proj on $\mathbb{T}^+/\mathcal{O}$ is induced by $\delta$ which takes into account all the representatives. We show now that it suffices to consider only the nonnegative representatives. We define the function

$$\delta^+ : (\mathbb{T}^+/\mathcal{O})^2 \longrightarrow \mathbb{R}^+,$$
$$(\bar{x}, \bar{y}) \longmapsto \inf_{\substack{\theta_1 \in A(x) \\ \theta_2 \in A(y)}} d_{\mathbb{T}}(x^{(\theta_1)}, y^{(\theta_2)}).$$

LEMMA 4.13. *The topologies induced by $\delta$ and by $\delta^+$ on $\mathbb{T}^+/\mathcal{O}$ are identical.*

PROOF. In other words, we have to prove that

$$\mathrm{Id} : (\mathbb{T}^+/\mathcal{O}, \delta) \longrightarrow (\mathbb{T}^+/\mathcal{O}, \delta^+),$$
$$x \longmapsto x$$

is a homeomorphism. Since $\delta(\bar{x}, \bar{y}) \leq \delta^+(\bar{x}, \bar{y})$, it is sufficient to show that Id is continuous. Let $y = (f, \zeta) \in \mathbb{T}^+$. Let us prove that, for any $\varepsilon > 0$, there exists $\rho > 0$, such that if $x \in \mathbb{T}^+$ and $\delta(\bar{x}, \bar{y}) < \rho$, then $\delta^+(\bar{x}, \bar{y}) < \varepsilon$.

For any $x \in \mathbb{T}^+$, there exists $\theta(x) \in [0, 1]$ such that

$$(4.14) \qquad d_{\mathbb{T}}(x, y^{(\theta(x))}) \leq 4\delta(\bar{x}, \bar{y}).$$

Set $A'(y) = \{\theta \in [0, 1], f(\theta) = \min f\}$ and $\alpha(\theta) = \inf\{|\theta - s|, s \in A'(y)\}$. Due to the uniform continuity of $(f, \zeta)$, there exists $\eta > 0$ such that

$$\alpha(\theta(x)) < \eta \Rightarrow d_{\mathbb{T}}(y^{(\theta(x))}, y^s) < \varepsilon/2$$

for some $s \in A'(y)$. Since $x, y \in \mathbb{T}^+$, there exists $\eta' > 0$ such that $d_{\mathbb{T}}(x, y^{(\theta(x))}) < \eta'$ implies $\alpha(\theta(x)) \leq \eta$. Choose a $\rho < \eta'/4 \wedge \varepsilon/8$,

$$\delta^+(\bar{x}, \bar{y}) \leq d_{\mathbb{T}}(x, y^s) \leq d_{\mathbb{T}}(x, y^{(\theta(x))}) + d_{\mathbb{T}}(y^{(\theta(x))}, y^s) \leq \varepsilon/2 + \varepsilon/2. \qquad \square$$

As a corollary, we have the following:

COROLLARY 4.14. *The function* Proj *is a homeomorphism from* $(\mathbb{T}/\mathcal{O}, D)$ *on* $(\mathbb{T}^+/\mathcal{O}, \delta^+)$.

4.5.2. *Positive Brownian snake and notation.* Here is a result concerning the minimum of the Brownian snake due to T. Duquesne:

LEMMA 4.15 (Duquesne). $\#A(\mathbf{r}, \mathbf{v}) = 1$ *a.s.*

Duquesne provides us a first proof of this lemma in a personal communication. For a published proof, see [27], Proposition 2.5.

This lemma implies that $\mathbb{P}_{\mathcal{S}}(\mathbb{T}^{+\star}) = 1$ and:

COROLLARY 4.16.  *We have* $\#\overline{(\mathbf{r},\mathbf{v})}^+ = 1$ *a.s.* $[i.e.\ \mathbb{P}_{\bar{\mathcal{S}}}(\mathbb{T}^{+\star}/\mathcal{O}) = 1]$.

We now introduce notation. We denote by $\bar{\mathbb{P}}_n^+ = \mathbb{P}_{\bar{\mathcal{S}}}^n \circ \text{Proj}^{-1}$ the law of $\overline{(\mathbf{r}_n, \mathbf{v}_n)}^+$ and by $\bar{\mathbb{P}}^+ = \mathbb{P}_{\bar{\mathcal{S}}} \circ \text{Proj}^{-1}$ the law of $\overline{(\mathbf{r},\mathbf{v})}^+$. According to Lemma 4.15, $\#\overline{(\mathbf{r},\mathbf{v})^+} = 1$, $\bar{\mathbb{P}}^+$-a.s. Since $(x \in \mathbb{T}^{+\star}) \Leftrightarrow (\bar{x}^+ = \{x\}) \Leftrightarrow (\#\bar{x}^+ = 1)$, we identify $\mathbb{T}^{+\star}$ and $\mathbb{T}^{+\star}/\mathcal{O}$, and we define by this way a law $\mathbb{P}^+$ on $\mathbb{T}^+$ by setting

$$\mathbb{P}^+(A) = \bar{\mathbb{P}}^+(A \cap \mathbb{T}^{+\star}) \qquad \text{for any Borel set } A \in \mathbb{T}^+.$$

Once again, $\mathbb{P}^+(\mathbb{T}^{+\star}) = 1$. In the sequel we will write $(\mathbf{r},\mathbf{v})^+$ for a $\mathbb{T}^+$-valued and $\mathbb{P}^+$-distributed random variable.

The random variables $(\mathbf{r},\mathbf{v})^+$ and $\overline{(\mathbf{r},\mathbf{v})}^+$ are not really different since one can consider that their distributions $\mathbb{P}^+$ and $\bar{\mathbb{P}}^+$ are the same distribution concentrated on $\mathbb{T}^{+\star}$. However, we must keep a distinction because $(\mathbf{r}_n, \mathbf{v}_n)^+$ and $\overline{(\mathbf{r}_n, \mathbf{v}_n)}^+$ take respectively their values in $\mathbb{T}^+$ and $\mathbb{T}^+/\mathcal{O}$.

4.5.3. *Convergence in* $\mathbb{T}^+/\mathcal{O}$. Since Proj is continuous, the weak convergence of $\overline{(\mathbf{r}_n,\mathbf{v}_n)}$ to $\overline{(\mathbf{r},\mathbf{v})}$ implies that $\overline{(\mathbf{r}_n, \mathbf{v}_n)}^+$ converges weakly to $\overline{(\mathbf{r},\mathbf{v})}^+$. Hence, we have the following:

PROPOSITION 4.17.  *We have* $\bar{\mathbb{P}}_n^+ \xrightarrow{weakly}_n \bar{\mathbb{P}}^+$.

4.5.4. *A.s. convergence to an element of* $\mathbb{T}^{+\star}$. Since the metric space $\mathbb{T}$ is separable, the metric space $\mathbb{T}^+/\mathcal{O}$ is also separable. Thus, we can apply the Skorohod representation theorem [23], Theorem 4.30, page 79. Thanks to Corollary 4.16 and to Proposition 4.17, there exists a probability space $(\Xi, \mathbb{P}_\Xi)$ on which are defined the random variables $(\overline{(\mathsf{r},\mathsf{v})}^+), (\overline{(\mathsf{r}_n,\mathsf{v}_n)}^+)_{n \geq 1}$ such that:

- $\overline{(\mathsf{r},\mathsf{v})}^+$ is $\mathbb{T}^+/\mathcal{O}$-valued and $\bar{\mathbb{P}}^+$-distributed.
- $\overline{(\mathsf{r}_n,\mathsf{v}_n)}^+$ is $\mathbb{T}^+/\mathcal{O}$-valued and $\bar{\mathbb{P}}_n^+$-distributed, for any $n \geq 1$.
- $\#\overline{(\mathsf{r},\mathsf{v})}^+(\omega) = 1$ for all $\omega \in \Xi$.
- $\overline{(\mathsf{r}_n,\mathsf{v}_n)}^+$ converges $\mathbb{P}_\Xi$-a.s. to $\overline{(\mathsf{r},\mathsf{v})}^+$.

Hence,

$$(4.15) \qquad \delta^+(\overline{(\mathsf{r}_n,\mathsf{v}_n)}^+, \overline{(\mathsf{r},\mathsf{v})}^+) \longrightarrow 0 \qquad \text{a.s.}$$

This means that the distance between the closest elements of the two classes $\overline{(r_n, v_n)}^+$ and $\overline{(r, v)}^+$ goes to 0. One has $\overline{(r, v)}^+ = \{(r, v)^+\}$. Let

$$(4.16) \quad d^+(\overline{(r_n, v_n)}^+, (r, v)^+) = \max\{d_{\mathbb{T}}((f, \zeta), (r, v)^+); (f, \zeta) \in \overline{(r_n, v_n)}^+\}$$

be the maximal distance between the elements of $\overline{(r_n, v_n)}^+$ and $(r, v)^+$. Since $\#\overline{(r, v)}^+(\omega) = 1$ for all $\omega \in \Xi$, equation (4.15) implies that the diameter of the classes $\overline{(r_n, v_n)}^+$ goes to 0 a.s. And thus we have the following:

PROPOSITION 4.18.  *On* $(\Xi, \mathbb{P}_\Xi)$, $d^+(\overline{(r_n, v_n)}^+, (r, v)^+) \xrightarrow[n]{a.s.} 0$.

4.6. *Proof of Theorem 3.37.* This proof is a recap of the construction done in the previous section. Using Proposition 4.4, pointed quadrangulations and abstract pointed quadrangulations are identified. Under $\mathbb{P}_{\bar{\mathcal{S}}}^n$, the class $\overline{(\mathbf{r}_n, \mathbf{v}_n)}^+$, that encodes well-labeled trees equal up to a rerooting, is $\mathbb{P}_n^+$-distributed. The class $\overline{(\mathbf{r}_n, \mathbf{v}_n)}^+$ encodes, via the passage by corresponding pair of doddering and gluer trees, a class of rooted abstract maps corresponding to a pointed abstract map $\bar{\mathbf{M}}_n^\bullet$ (which is a random variable $\mathfrak{M}_n^\bullet$ valued and $\mathbb{P}_{\bar{\mathcal{S}}}^n$-distributed). The convergence of the diameter of $\overline{(\mathbf{r}_n, \mathbf{v}_n)}^+$ to 0, and the convergence of $\overline{(\mathbf{r}_n, \mathbf{v}_n)}^+$ under $\mathbb{P}_{\bar{\mathcal{S}}}^n$ to $\overline{(\mathbf{r}, \mathbf{v})}^+$—which is identified with $(\mathbf{r}^+, \mathbf{v}^+)$—and the convergence of associated discrete corners measures to (deterministic) Lebesgue measure, allow to conclude that $\bar{\mathbf{M}}_n^\bullet \xrightarrow[n]{weakly} \mathbf{M}_\infty^\bullet$, and then thanks to Proposition 4.7, $\mathbf{M}_n^\bullet \xrightarrow[n]{weakly} \mathbf{M}_\infty^\bullet$.

4.7. *Proof of Theorem 3.3.* Set $T \in \mathcal{W}_n$ and let $(f_n, \zeta_n)$ be its encoding in $\mathbb{T}_n$. The different elements of $\overline{(f_n, \zeta_n)}^+$ are the encodings of the well-labeled trees obtained from $T$ by a rerooting. In the terminology of Section 4.2.1, the elements of $\overline{(f_n, \zeta_n)}^+$ are the encodings of the well-labeled trees $\{L^{\circ k}(T), k \in \mathbb{N}\}$ and so of the rooted quadrangulations associated with $\widetilde{Q}(T)$ [i.e., $K^{-1}(\widetilde{Q}(T))$]. Hence, if $\widetilde{Q}(T)$ has no symmetry, $\#\overline{(f_n, \zeta_n)}^+$ is the degree of the root of $\widetilde{Q}(T)$ and it is also $\#A(f_n, \zeta_n)$, the number of minima of $f_n$ in $[0, 1)$.

The encoding $(\mathbf{r}_n^+, \mathbf{v}_n^+)$ of normalized well-labeled trees (under $\mathbb{P}_D^n$) are random variables in $\mathbb{T}^+$. For each element $\bar{x}^+ \in \mathbb{T}_n^+/\mathcal{O}$, we define $N(\bar{x}^+) = \#\bar{x}^+$ and we fix a numbering $x^1, \ldots, x^{N(x)}$ of the representatives of $\bar{x}^+$ that are element of $\mathbb{T}_n^+$. Hence, $\mathbb{T}_n^+$ is the disjoint union of the sets $\{x^1, \ldots, x^{N(x)}\}$ for $\bar{x}^+ \in \mathbb{T}_n^+/\mathcal{O}$ (recall that $\mathbb{T}_n^+/\mathcal{O} \simeq \mathbb{T}_n^+/\mathcal{O}_n$). The following lemma gives the relation we announced in the Introduction, between the two distributions $(\mathcal{Q}_n^\bullet, \mathbb{P}_{\bar{\mathcal{S}}}^n)$ and $(\vec{\mathcal{Q}}_n, \mathbb{P}_D^n)$.

LEMMA 4.19. *Let $\bar{X}_n^+$ be $\bar{\mathbb{P}}_n^+$-distributed and $\bar{U}_n$ be uniformly distributed on $[\![1, N(\bar{X}_n^+)]\!]$ conditionally on $\bar{X}_n^+$. The $\mathbb{T}^+$-valued random variable $X_n^{\bar{U}_n}$ is $\mathbb{P}_D^n$-distributed. Moreover, for each $q \in \vec{\mathcal{Q}}_n$ with root-degree $\deg(q)$,*

$$\mathbb{P}_D^n(q) = \frac{2n}{C_n 3^n} \frac{1}{\deg(q)}.$$

PROOF. Let $x_n$ be any element of $\mathbb{T}_n$, and $\bar{x}_n^+ = \text{Proj}(\bar{x}_n)$. One has

$$\bar{\mathbb{P}}_n^+(\bar{x}_n^+) = \mathbb{P}_{\bar{\mathcal{S}}}^n(\bar{x}_n) = \frac{2n}{\#\text{Stab}(x_n)} \frac{1}{C_n 3^n},$$

where $\text{Stab}(x_n) = \{\theta \in \mathcal{O}_n, x_n^{(\theta)} = x_n\}$. Given $\bar{X}_n^+$, the random variable $X_n^{\bar{U}_n}$ is uniform on the set $\{X_n^1, \ldots, X_n^{N(\bar{X}_n^+)}\}$. The cardinality $N(\bar{X}_n^+)$ of the class $\bar{X}_n^+$ is equal to $\#A(f_n, \zeta_n)/\#\text{Stab}(f_n, \zeta_n)$ for any $(f_n, \zeta_n) \in \bar{X}_n^+$. For any $(f_n, \zeta_n) \in \mathbb{T}_n^+$, $\mathbb{P}(X_n^{\bar{U}_n} = (f_n, \zeta_n))$ is equal to

$$\mathbb{P}(X_n^{\bar{U}_n} = (f_n, \zeta_n) | \bar{X}_n^+ = \overline{(f_n, \zeta)}^+) \mathbb{P}(\bar{X}_n^+ = \overline{(f_n, \zeta)}^+)$$

$$= \frac{\#\text{Stab}(f_n, \zeta_n)}{\#A(f_n, \zeta_n)} \frac{2n}{\#\text{Stab}(f_n, \zeta_n)} \frac{1}{C_n 3^n} = \frac{2n}{C_n 3^n} \frac{1}{\#A(f_n, \zeta_n)}.$$

Now, $\mathbb{P}_D^n(f_n, \zeta_n)$ is also proportional to $1/\#A(f_n, \zeta_n)$ for any $(f_n, \zeta_n) \in \mathbb{T}_n^+$ [the degree of the root-vertex of the quadrangulation encoded by $(f_n, \zeta_n)$ is $\#A(f_n, \zeta_n)$], thus $X_n^{\bar{U}_n}$ is $\mathbb{P}_D^n$-distributed. □

PROOF OF THEOREM 3.3. The idea of the proof is to extract from the sequence of classes $\overline{(r_n, v_n)}^+$ that converges a.s. in $\Xi$ to $\overline{(r, v)}^+ = \{(r, v)^+\}$, a sequence of $\mathbb{T}_n^+$ valued and $\mathbb{P}_D^n$-distributed random variables, $(r_n^+, v_n^+)$, that converges a.s. to $(r, v)^+$. For this, we construct a new space $\Xi'$.

Since we have identified $\mathbb{T}^{+\star}$ with $\mathbb{T}^{+\star}/\mathcal{O}$, we consider $\overline{(r, v)}^+$ as equal to its unique representative $(r, v)^+$ and thus as a $\mathbb{T}^+$ valued and $\mathbb{P}^+$-distributed random variable. We consider the probability space

$$\Xi' = \Xi \times [0, 1]^{\mathbb{N}^\star} \qquad \text{endowed with } \mathbb{P}_{\Xi'} = \mathbb{P}_\Xi \otimes (\text{Leb}^{[0,1]})^{\otimes \mathbb{N}^\star}.$$

Let $(U_i)_{i \in \mathbb{N}}$ be the coordinate function in $[0, 1]^{\mathbb{N}^\star}$; the random variable $(((r, v)^+, (\overline{(r_n, v_n)}^+)_{n \geq 1}), (U_i)_{i \geq 1})$ is defined on $\Xi'$ and takes its values in $(\mathbb{T}^+ \times (\mathbb{T}^+/\mathcal{O})^{\mathbb{N}^\star}) \times [0, 1]^{\mathbb{N}^\star}$. The random variables $U_i$ are i.i.d., uniform on $[0, 1]$, independent of the sequence $(\overline{(r_n, v_n)}^+)_{n \geq 1}$; the random variable

$$(r_n^+, v_n^+) := (\overline{(r_n, v_n)}^+)^{\lceil U_n N(\overline{(r_n, v_n)}^+) \rceil}$$

is $\mathbb{T}^+$-valued and is $\mathbb{P}_D^n$-distributed according to Lemma 4.19. To conclude the proof, it suffices to show that $(r_n^+, v_n^+)$ converges a.s. to $(r, v)^+$. One has

$$d_{\mathbb{T}}((r_n^+, v_n^+), (r, v)^+) \leq \sup\{d_{\mathbb{T}}((f, \zeta), (r, v)^+); (f, \zeta) \in \overline{(r_n, v_n)}^+\};$$

the term in the right-hand side goes to 0 a.s. thanks to Proposition 4.18. □

REMARK 4.20. We conjecture that Theorem 3.3 is also true if $\overrightarrow{\mathcal{Q}}_n$ is endowed with the uniform law. If this conjecture is right, this model of rooted quadrangulations converges weakly to the Brownian map with the same normalization as in the model $\mathbb{P}_D^n$.

**5. Asymptotic of functionals of quadrangulations.** We may deduce from Theorem 3.3 (resp. Theorem 3.37) the weak convergence of continuous functions of rooted (resp. pointed) quadrangulations under $\mathbb{P}_D^n$ (resp. $\mathbb{P}_U^n$). Among them, is the radius:

*Convergence of the radius.* The radius of a map is the largest distance between the root and a node. On the encoding, it is the largest label of $\mathbf{R}_n^+$ for the rooted map, and the range of the label process for pointed maps:

- In the model $(\overrightarrow{\mathcal{Q}}_n, \mathbb{P}_D^n)$, the normalized radius (i.e., the radius divided by $n^{1/4}$) is $\mathrm{Rad}_n = \max \mathbf{r}_n^+$. We have $\mathrm{Rad}_n \xrightarrow[n]{(law)} \max \mathbf{r}^+$.
- In the model $(\mathcal{Q}_n^\bullet, \mathbb{P}_{\tilde{\mathcal{S}}}^n)$, the normalized radius is $\mathrm{Rad}_n = \max \mathbf{r}_n - \min \mathbf{r}_n$; we have $\mathrm{Rad}_n \xrightarrow[n]{(law)} \max \mathbf{r} - \min \mathbf{r}$.

According to the construction of Section 4.7, the two limit laws are equal. They are also equal to the limit radius when $\overrightarrow{\mathcal{Q}}_n$ is endowed with the uniform distribution with the same normalization (see [15]).

*Convergence of the profile.* We follow the steps of Chassaing and Schaeffer [15], Section 6.4. Let

$$\mathbf{L}^{(n)}(j) = \frac{\#\{k, \mathbf{r}_n^+(k) \leq j/n^{1/4}\}}{2n}$$

be the proportion of edges incident to nodes at distance smaller than $j - 1$ of the root in a rooted quadrangulation under $\mathbb{P}_D^n$. Assume that $\mathbf{L}^{(n)}$ is interpolated between integer points. For $\lambda \in \mathbb{R}$, let $\mathbf{l}^{(n)}(\lambda) = \mathbf{L}^{(n)}(\lambda n^{1/4})$.

Using the arguments of [15], we can consider our random variables defined on a probability space $\widetilde{\Xi}$, on which ISE shifted is absolutely continuous and on which $\mathbf{r}_n^+$ converges to $\mathbf{r}^+$, such that, a.s., for any fixed $\lambda$,

$$\mathbf{l}^{(n)}(\lambda) \longrightarrow \mathbf{l}(\lambda) = \int_0^1 \mathbb{1}_{[0,\lambda]}(\mathbf{r}^+(s))\, ds.$$

Since the involved functions are increasing, continuous and bounded, this convergence is uniform in $\lambda$. Finally,

$$\mathbf{l}^{(n)} \xrightarrow[n]{weakly} \mathbf{l}$$

on $C((-\infty,+\infty))$ endowed with the topology of uniform convergence on the compact sets. By Proposition 4.18, it is straightforward that the same result holds for pointed quadrangulations under $\mathbb{P}^n_U$.

*Open question.* Does $d_{\mathbf{M}_n}(S(F_{\mathbf{D}_n}(s)), S(F_{\mathbf{D}_n}(t)))$ weakly converge to $d_{\mathbf{M}_\infty}(S(F_{\mathbf{D}_\infty}(s)), S(F_{\mathbf{D}_\infty}(t)))$ for any fixed $s$ and $t$ in $[0,2]$? If it is true, is it also true for the process $(d_{\mathbf{M}_n}(S(F_{\mathbf{D}_n}(s)), S(F_{\mathbf{D}_n}(t))))_{(s,t)\in[0,2]^2}$?

The considered topology on abstract maps does not allow to prove convergence of functionals of the node-degrees to the ones of the Brownian map. We conjecture the following fact about the nodes-degrees in the Brownian map:

CONJECTURE 5.1. *In the Brownian map, a.s.,*

$$\max_{\hat{x}\in M_\infty} \deg(\hat{x}) = 3.$$

## 6. Conclusion.

6.1. *A bijection between maps and quadrangulations.* We denote by $\overrightarrow{M_n}$ (resp. by $M_n^\bullet$) the set of rooted (resp. pointed) maps with $n$ edges.

PROPOSITION 6.1. *The following bijections hold:*

$$M_n^\bullet \simeq \mathcal{Q}_n^\bullet, \qquad \overrightarrow{M_n} \simeq \overrightarrow{\mathcal{Q}}_n.$$

These bijections are classical. For reader convenience, we give a short description of them below (see also illustrations in Figure 19):

• Take a map $m \in M_n^\bullet$ pointed in $u$. Color this map in blue. Add a red node in each face. In each face, add a red edge between the red node and each of the blue nodes adjacent to this face. Denote by $q$ the map pointed in $u$ that has for the set of nodes, the red nodes and the blue nodes, and that has for edges, the red edges. The map $q$ is a quadrangulation that has $n$ faces. Indeed, each face of $q$ contains exactly one edge of $M_n^\bullet$, and each face of $q$ has degree 4. This construction is invertible; let us start from $q$ pointed in $u$ and rebuild $m$. The quadrangulation $q$ is bipartite. Color in blue the nodes at even distance from $u$, in red the other ones. In each face add a blue edge between the blue nodes. The blue map with distinguished point $u$ is $m$.

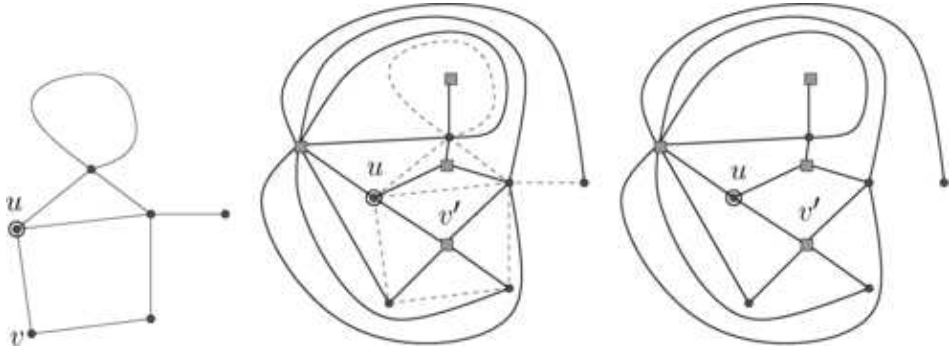

FIG. 19. *From a map to a quadrangulation.*

• The construction of the bijection from $\overrightarrow{M_n}$ onto $\overrightarrow{\mathcal{Q}}_n$ is the same except that we have to consider the root (instead of the distinguished node). Let $\overrightarrow{uv}$ be the root of $m$. The root $\overrightarrow{u'v'}$ of $q$ is chosen as follows: it is the first red edge on the left of $\overrightarrow{uv}$ with origin $u$. Starting from $q$ rooted on $\overrightarrow{u'v'}$, the root $\overrightarrow{uv}$ of $m$ is chosen as the blue edge added in the face adjacent to $\overrightarrow{u'v'}$ at the right of $\overrightarrow{u'v'}$ (with origin $u$).

6.2. *Limit of other models of maps.* The results we obtained suggest that (some) other models of maps can be treated in the same way as quadrangulations and that the Brownian map should be a natural limit for other normalized random maps. From our work, two main approaches can be proposed in order to prove the weak convergence other models, of random maps:

• The first one is to find a representation of the maps of the considered model by a product of two trees (or by a product of a tree and a forest). Indeed, from each map, one can extract a tree that contains all the edges once. The gluings needed to build the map back can be encoded by a noncrossing partition (or a parenthesis system) that should be encoded by a forest (or a tree).

Schaeffer ([33], Chapter I) gives several representations of different models of maps that are obtained from marked trees with the help of more or less complex procedures of gluing.

• The second approach uses Proposition 6.1 to encode any map with the help of a quadrangulation. A random model of rooted maps is a probability $\mathbb{P}$ on $\overrightarrow{M} = \bigcup_k \overrightarrow{M_k}$. By the bijection presented in Proposition 6.1, $\mathbb{P}$ induces a law $\mathbb{P}'$ on $\overrightarrow{\mathcal{Q}}$ (where $\overrightarrow{\mathcal{Q}} = \bigcup_k \overrightarrow{\mathcal{Q}_k}$) (e.g., if $\mathbb{P}$ is the uniform distribution on $\overrightarrow{M_n}$, $\mathbb{P}'$ is the uniform distribution on $\overrightarrow{\mathcal{Q}}_n$).

A sequence of models of maps is a sequence of probabilities $\mathbb{P}_n$ on $\overrightarrow{M}$; it is transported (by the bijection) as a sequence of probabilities $\mathbb{P}'_n$ on $\overrightarrow{\mathcal{Q}}$. The study of the asymptotics of $\mathbb{P}_n$ reduces to the study of the asymptotic of $\mathbb{P}'_n$.

One may hope that, for simple models of maps (maps defined by restriction on the degree of the faces or degree of the vertices as triangulations), the $\mathbb{P}'_n$-distributed process $(r_n, v_n)$ would converge to the $\mathbb{P}^+$-distributed process $(r, v)^+$ (up to some scaling).

Notice that $\mathbb{P}'_n$ can give weight to quadrangulations with different sizes. Even in this case, it may exist a normalization $(r_n, v_n)$ of $(R_n, V_n)$ which converges.

To conclude, one can see that the convergence of $(r_n, v_n)$ to $(r, v)$ (given in Proposition 4.8) is quite robust: one can change the law of the underlying tree (take uniform binary trees, uniform ternary trees, other models of simply generated trees) and the law of the increments (here the increments $X$ are uniform in $\{-1, +1, 0\}$, one can take any $X$ with symmetric law in $\{-1, +1, 0\}$): in these cases the corresponding encoding $(r_n, v_n)$ converges to the same limit $(r, v)$ (up to constant scales). Each of these models of law [of $(r_n, v_n)$] corresponds to a model of finite quadrangulations.

## APPENDIX

### A.1. Relation between CTrees and RTrees.

PROPOSITION A.1. *Let $(g, \mu)$ be a tree-encoding, and $\mathcal{T} = \mathrm{CTree}(g, \mu)$. Let $F_\mathcal{T}$ be the CDFT of $\mathcal{T}$, $G_\mathcal{T}$ be its RDFT, $l_\mathcal{T}$ be its RDFW and $\mu'$ be its reverse measure. Set $\widetilde{\mathcal{T}} = \mathrm{RTree}(l_\mathcal{T}, \mu')$. The application $\Phi$ from $E_{\widetilde{\mathcal{T}}}$ on $E_\mathcal{T}$ defined by*

$$\Phi(\dot{x}) = F_\mathcal{T}(\Psi_T(x)) = G_\mathcal{T}(x) \tag{A.1}$$

*is a bijective isometry. Moreover, $E_\mu = \Phi(E_{\mu'})$.*

PROOF. The RDFW of $\widetilde{\mathcal{T}}$ is $l_\mathcal{T}$ and its reverse measure is $\mu'$. We have

$$\begin{aligned}d_{\widetilde{\mathcal{T}}}(\dot{x}, \dot{y}) &= l_\mathcal{T}(x) + l_\mathcal{T}(y) - 2\check{l}_\mathcal{T}(x, y) \\ &= g(\Psi_\mathcal{T}(x)) + g(\Psi_\mathcal{T}(y)) - 2\check{g}(\Psi_\mathcal{T}(x), \Psi_\mathcal{T}(y)).\end{aligned} \tag{A.2}$$

The relation (A.2) implies that $(x \underset{\widetilde{\mathcal{T}}}{\sim} y) \Leftrightarrow (\Psi_\mathcal{T}(x) \underset{\mathcal{T}}{\sim} \Psi_\mathcal{T}(y))$.

Now, since $\mu'(A) = \mu(\Psi_\mathcal{T}(A))$, we have $\mathrm{supp}(\mu) = \Psi_\mathcal{T}(\mathrm{supp}(\mu'))$ and thus, $E_\mu = F_\mathcal{T}(\Psi_\mathcal{T}(\mathrm{supp}\,\mu'))$, which, in view of (A.1), implies that $E_\mu = \Phi(E_{\mu'})$. □

According to the description of the CTree and RTree,

$$\left(\dot{x} \underset{CO}{\preccurlyeq} \dot{y} \text{ in } \widetilde{\mathcal{T}}\right) \Leftrightarrow \left(\Phi(\dot{x}) \underset{CO}{\preccurlyeq} \Phi(\dot{y}) \text{ in } \mathcal{T}\right),$$

$$\left(\dot{x} \underset{RO}{\preccurlyeq} \dot{y} \text{ in } \widetilde{\mathcal{T}}\right) \Leftrightarrow \left(\Phi(\dot{x}) \underset{RO}{\preccurlyeq} \Phi(\dot{y}) \text{ in } \mathcal{T}\right).$$

One may say that $T$ and $\widetilde{T}$ represent the same tree.

REMARK A.2. Let $T = \text{CTree}(g, \mu)$ and $T_S = \text{RTree}(g, \mu)$. Obviously the function $\text{Id}: x \mapsto x$ in $[0, a]$ induces an isometry $\dot{x} \mapsto \phi(\dot{x})$ between $E_T$ and $E_{T_S}$. But now,

$$\left(\dot{x} \underset{RO}{\preccurlyeq} \dot{y} \text{ in } T\right) \Leftrightarrow \left(\phi(\dot{x}) \underset{CO}{\preccurlyeq} \phi(\dot{y}) \text{ in } T_S\right),$$

$$\left(\dot{x} \underset{CO}{\preccurlyeq} \dot{y} \text{ in } T\right) \Leftrightarrow \left(\phi(\dot{x}) \underset{RO}{\preccurlyeq} \phi(\dot{y}) \text{ in } T_S\right).$$

We say that $T$ and $T_S$ are symmetric. In the discrete case, a planar representation of the tree $T_S$ is obtained from a planar representation of $T$ by any axial symmetry of the plane.

**A.2. A bijection between $\mathcal{Q}_n^\bullet$ and a set of unrooted marked trees.** Consider a tree $T$ in $\mathcal{W}_n$; let $F$ be its CDFT and $R$ its label process. Each edge of $T$ is traversed twice by $F$. Then, we can define the two sides of an edge: if $k$ and $l$ are such that $k < l$, $F(k) = F(l+1)$ and $F(k+1) = F(l)$, one says that $\overrightarrow{F(k)F(k+1)}$ is the first side and $\overrightarrow{F(l)F(l+1)}$ is the second side. We mark now the sides of the edges: the mark of the side $\overrightarrow{F(j)F(j+1)}$ is $R(j+1) - R(j) \in \{+1, -1, 0\}$. We have thus defined an application $X$ from $\mathcal{W}_n$ on $\mathcal{U}_n$, the set of rooted plane trees with $n$ edges where the edges-sides are marked $\mathcal{U}_n \simeq \{+1, -1, 0\}^n \times \Omega_n$ (since the mark of a side determines the mark of the other side). The application $X$ is a bijection since one can recover all the labels with the marks and set the root vertex label equal to 1.

Now, consider two trees $T_1$ and $T_2$ in $\mathcal{W}_n$ being in the same class modulo $G_n$. Since the marks do not depend on the position of the root, and are the same if one adds a constant to all the labels of $T$, the two marked trees $X(T_1)$ and $X(T_2)$ differ only by their root positions. Conversely, if $X(T_1)$ and $X(T_2)$ differ only by their root position, it is straightforward to see that $T_1$ and $T_2$ are in the same class modulo $G_n$.

Thus, $G_n$ acts on $\mathcal{U}_n$ by rerooting the branching structure without moving the marks. Hence, $X$ induces a bijection from $\mathcal{W}_n/G_n$ on $\mathcal{U}_n/G_n$. The set $\mathcal{U}_n/G_n$ is naturally identified with the set of unrooted plane trees with $n$ edges, where each edge side is marked by $+1, -1$ or $0$ (the mark of a side being the opposite of the mark of the other side).

**A.3. Proof of Lemma 3.19.** The second statement of this lemma is a consequence of the first one. To prove the first one, we use the argument of [28] (Section 3.1). The reverse DFW and the height process both visit the nodes of a given tree $\tau$ in the reverse depth first order. If one sorts the nodes of $\tau$ according to the reverse depth first order $(v_0, v_1, \ldots, v_{n-1})$, then the reverse height process $\widetilde{H}_n$ can be expressed by

$$\widetilde{H}_n(l) = d(root, v_l)$$

and the reverse DFW is given by

$$\widetilde{V}_n(k) = \begin{cases} \widetilde{H}_n(l), & \text{if } k = m(l) \text{ for a given } l, \\ \widetilde{H}_n(l) - (k - m(l)), & \text{if } k \in [\![m(l)+1, m(l+1)-1]\!] \\ & \text{for a given } l, \end{cases}$$

where

$$m(l) = \inf\{k | \widetilde{F}(k) = v_l\};$$

moreover, $m(l)$ satisfies

(A.3) $$m(l) + \widetilde{H}_n(l) = 2l.$$

Let $j$ be the integer such that $m(j) + 1 \leq 2l \leq m(j+1) - 1$. Thanks to (A.3), $j \in [l-1, l + \max_k \widetilde{H}_n(k)]$. As a consequence,

$$\sup_l \frac{|\widetilde{H}_n(l) - \widetilde{V}_n(2l)|}{c_n} \leq \sup_l \sup_{j \in [l-1, l+\max_k \widetilde{H}_n(k)]} \frac{|\widetilde{H}_n(l) - \widetilde{H}_n(j)| + 1}{c_n}$$

and, thus,

$$\sup_{t \in [0,1]} \frac{|\widetilde{V}_n(2nt) - \widetilde{H}_n(nt)|}{c_n} \leq \sup_{t \in [0,1]} \sup_{|s-t| \leq (\max \widetilde{H}_n + 1)/n} \frac{|\widetilde{H}_n(nt) - \widetilde{H}_n(ns)| + 1}{c_n}.$$

Let us denote by $\delta : C([0,1]) \times [0,1] \longrightarrow +\infty$ the continuity modulus

$$\delta(f, \varepsilon) = \sup_{|t-s| \leq \varepsilon} |f(t) - f(s)|;$$

we then have

(A.4) $$\sup_{t \in [0,1]} \frac{|\widetilde{V}_n(2nt) - \widetilde{H}_n(nt)|}{c_n} \leq \delta\left(\frac{\widetilde{H}(n\cdot)}{c_n}, \frac{\max \widetilde{H}_n(n\cdot) + 1}{c_n}\right).$$

In view of the assumptions of the lemma and since $\delta(\cdot, \cdot)$ is continuous, it follows that the right-hand side of (A.4) converges weakly to $\delta(h, 0) = 0$, which proves the first part of the lemma.

**A.4. Proof of Theorem 3.34.** We use Theorem 3.33. Let $(\mathbf{c}_n^+, \mathbf{v}_n^+)$ [resp. $(\widetilde{\mathbf{c}}_n^+, \mathbf{v}_n^+)$] be the normalized DFW of $(\mathbf{D}_n, \mathbf{G}_n)$ [resp. $(\widetilde{\mathbf{D}}_n, \mathbf{G}_n)$]. Since $\widetilde{\mathbf{M}}_n$ and $\mathbf{M}_n$ differ only by the RDFW of their doddering trees,

$$d_{\mathfrak{M}(1)}(\widetilde{\mathbf{M}}_n, \mathbf{M}_n) = \|\widetilde{\mathbf{c}}_n^+ - \mathbf{c}_n^+\|_\infty.$$

We prove that there exists $\varepsilon > 0$ such that

(A.5) $$\mathbb{P}\left(\max_{k \in [\![0, 4n]\!]} \left|\mathbf{c}_n^+\left(\frac{k}{2n}\right) - \widetilde{\mathbf{c}}_n^+\left(\frac{k}{2n}\right)\right| \geq n^{-\varepsilon}\right) \longrightarrow 0.$$

Hence, $\|\widetilde{\mathbf{c}_n^+} - \mathbf{c}_n^+\|_\infty \xrightarrow[n]{proba} 0$ and by Theorem 3.33, this allows us to conclude. Consider a node $u$ in the underlying normalized doddering tree $\mathbf{D}_n$ and let $u_0 = root, u_1, \ldots, u_{d(u,root)}$ be the branch in $\mathbf{D}_n$ between $u$ and the root. Denote by $L(i)$ the edge length between $u_i$ and $u_{i-1}$. If $k$ is a representative of $u$ (in $\mathbf{c}_n^+$),

$$\widetilde{\mathbf{c}_n^+}(k/(2n)) = \frac{L(1) + \cdots + L(n^{1/4}\mathbf{c}_n^+(k/(2n)))}{n^{1/4}}. \tag{A.6}$$

Hence, for any $k \in [\![0, 2n]\!]$, we have the following representation:

$$\widetilde{\mathbf{c}_n^+}(k/(2n)) - \mathbf{c}_n^+(k/(2n)) \stackrel{(d)}{=} \frac{1}{n^{1/4}} \sum_{i=1}^{n^{1/4}\mathbf{c}_n^+(k/(2n))} (L(i) - 1), \tag{A.7}$$

where the $L(i)$ are i.i.d. $\mu$-distributed. Set $\varepsilon \in (0, 1/4)$. Thanks to Fuk and Nagaev's inequality (see [30], Addenda 2.6.5, page 78), for any $p \geq 2$, such that $\mathbb{E}(|L|^p) < +\infty$, there exist two constants $c(p)$ and $c'(p)$ depending only on $p$ such that

$$\begin{aligned}
&\max_{l \in [\![1, n^{1/4+\varepsilon}]\!]} \mathbb{P}\left(\left|\sum_{j=1}^{l} L(j) - 1\right| \geq n^{1/4-\varepsilon}\right) \\
&\leq \max_{l \in [\![1, n^{1/4+\varepsilon}]\!]} \left\{\frac{c(p)l}{n^{p/4-p\varepsilon}} + \exp\left(-\frac{c'(p)n^{1/2-2\varepsilon}}{l \operatorname{var}(L)}\right)\right\}.
\end{aligned} \tag{A.8}$$

This is $o(1/n)$ for $p = 5 + \delta$ and a well chosen $\varepsilon > 0$. We fix from now on this $\varepsilon > 0$. The left-hand side of (A.5) is bounded by

$$\mathbb{P}(\|\mathbf{c}_n^+\|_\infty > n^\varepsilon) + \mathbb{P}\left(\max_{k \in [\![0, 4n]\!]} \left|\mathbf{c}_n^+\left(\frac{k}{2n}\right) - \widetilde{\mathbf{c}_n^+}\left(\frac{k}{2n}\right)\right| \geq n^{-\varepsilon}, \|\mathbf{c}_n^+\|_\infty \leq n^\varepsilon\right)$$

$$\leq \mathbb{P}(\|\mathbf{c}_n^+\|_\infty > n^\varepsilon)$$
$$+ \sum_{k \in [\![0, 4n]\!]} \mathbb{P}\left(\left|\mathbf{c}_n^+\left(\frac{k}{2n}\right) - \widetilde{\mathbf{c}_n^+}\left(\frac{k}{2n}\right)\right| \geq n^{-\varepsilon}, \|\mathbf{c}_n^+\|_\infty \leq n^\varepsilon\right),$$

which goes to 0 thanks to (A.7), (A.6) and (A.8).

**A.5. Proof of Proposition 4.7.** Since the orbits in $\mathcal{W}_n$ under $G_n$ have not a constant size, $\mathbb{P}_{\widetilde{\mathcal{S}}}^n$ is not the uniform law in $\mathcal{Q}_n^\bullet$. For $T \in \mathcal{W}_n$, denote by

$$\operatorname{Stab}(T) = \{\theta \in [\![0, 2n - 1]\!] \text{ s.t. } T^{(\theta)} = T\}$$

the stabilizer of $T$; the size of the orbit of $T$ is $2n/\#\operatorname{Stab}(T)$. Let $q$ be an element of $\mathcal{Q}_n^\bullet$ and let $T$ be an element of $\mathcal{W}_n$ such that $\widetilde{Q}(T) = q$, then

$$\mathbb{P}_{\widetilde{\mathcal{S}}}^n(q) = \frac{2n}{\#\operatorname{Stab}(T)} \frac{1}{C_n 3^n}.$$

Denote by $\mathcal{W}_n^\star$ (resp. $\Omega_n^\star$) the set of elements of $\mathcal{W}_n$ (resp. $\Omega_n$) whose orbit size by the action of $G_n$ is $2n$. Set $\mathcal{Q}_n^{\bullet\star} = \widetilde{Q}(\mathcal{W}_n^\star)$. For any $q \in \mathcal{Q}_n^{\bullet\star}$,
$$\mathbb{P}_{\tilde{\mathcal{S}}}^n(q) = 2n/(C_n 3^n).$$
We first show that $\mathbb{P}_{\tilde{\mathcal{S}}}^n(\mathcal{Q}_n^{\bullet\star}) \to 1$.

Since any labeled tree whose underlying tree belongs to $\Omega_n^\star$ is in $\mathcal{W}_n^\star$,

(A.9) $$3^n \#(\Omega_n^\star/G_n) = 3^n \frac{\#\Omega_n^\star}{2n} \leq \frac{\#\mathcal{W}_n^\star}{2n} \leq \#\mathcal{Q}_n^\bullet \leq 3^n \#(\Omega_n/G_n).$$

The last inequality comes from the fact that there are less than $3^n$ ways to mark the edges of an unrooted plan tree by $\{+1, -1, 0\}$. Walkup [36] shows that $\#(\Omega_n/G_n)$, the number of unrooted trees with $n$ edges, satisfies, for $n \geq 1$,

(A.10) $$\#(\Omega_n/G_n) = \frac{C_n}{2n} + \frac{1}{4n}\binom{n+1}{(n+1)/2} + \frac{1}{n}\phi(n)$$
$$+ \frac{1}{2n} \sum_{\substack{s|n \\ 1 < s < n}} \phi\left(\frac{n}{s}\right)\binom{2s}{s},$$

where the second term in the right-hand side is understood to be zero if $n$ is even, and where $\phi$ is the Euler totient function. The elements of $\Omega_n$ can be sorted according to their orbit sizes under the action of $G_n$. Denote by
$$\Omega_n^{[k]} = \{t \in \Omega_n \text{ such that } \#\{t^{(\theta)}, \theta \in G_n\} = k\}.$$
Notice that $\Omega_n^\star = \Omega_n^{[2n]}$ and that

(A.11) $$\frac{1}{C_n} \sum_{k|2n} \#\Omega_n^{[k]} = 1,$$

(A.12) $$\frac{1}{\#(\Omega_n/G_n)} \sum_{k|2n} \frac{\#\Omega_n^{[k]}}{k} = 1,$$

since the orbit sizes $k$ divide $2n$. Now, since in view of (A.10),

(A.13) $$\#(\Omega_n/G_n) \sim C_n/(2n),$$

formula (A.12) leads to
$$\frac{\#\Omega_n^{[2n]}}{C_n} + \sum_{\substack{k|2n \\ k \leq n}} \frac{2n \#\Omega_n^{[k]}}{kC_n} \xrightarrow[n]{} 1.$$

Subtracting (A.11), we obtain
$$\frac{1}{C_n} \sum_{\substack{k|2n \\ k \leq n}} \left(\frac{2n}{k} - 1\right) \#\Omega_n^{[k]} \longrightarrow 0.$$

Since, for the involved $k$, $(2n/k - 1) \geq n/k$, one also has

$$\frac{2}{C_n} \sum_{\substack{k|2n \\ k \leq n}} \frac{n}{k} \# \Omega_n^{[k]} \longrightarrow 0.$$

Finally, this gives

$$(A.14) \qquad \frac{2n}{C_n}\left(\#(\Omega_n/G_n) - \frac{\#\Omega_n^{[2n]}}{2n}\right) = \frac{2n}{C_n}\left(\sum_{\substack{k|2n \\ k \leq n}} \frac{\#\Omega_n^{[k]}}{k}\right) \longrightarrow 0.$$

Formulas (A.14), (A.9) and (4.7) yield

$$(A.15) \qquad \mathbb{P}_{\bar{\mathcal{S}}}^n(\mathcal{Q}_n^\bullet \setminus \mathcal{Q}_n^{\bullet\star}) \leq \frac{2n\,3^n(\#(\Omega_n/G_n) - \#(\Omega_n^\star/G_n))}{3^n C_n} = o(1)$$

and, thus, $\mathbb{P}_{\bar{\mathcal{S}}}^n(\mathcal{Q}_n^{\bullet\star}) \to 1$. A consequence of (A.9) and (A.14) is that

$$(A.16) \qquad \#\mathcal{Q}_n^\bullet \sim 3^n C_n/(2n).$$

(The exact enumeration of pointed quadrangulations is useless for our work.)

On the other hand, for any $q$ in $\mathcal{Q}_n^\bullet$, $\mathbb{P}_U^n(q) = (\#\mathcal{Q}_n^\bullet)^{-1}$. We have

$$\mathbb{P}_U^n(\mathcal{Q}_n^{\bullet\star}) = \frac{\#\mathcal{Q}_n^{\bullet\star}}{\#\mathcal{Q}_n^\bullet} \geq \frac{3^n \#\Omega_n^\star/G_n}{\#\mathcal{Q}_n^\bullet}.$$

Using (A.16), (A.14) and (A.13), we obtain that

$$(A.17) \qquad \mathbb{P}_U^n(\mathcal{Q}_n^{\bullet\star}) \longrightarrow 1.$$

For any $A$ subset of $\mathcal{Q}_n^\bullet$, $|\mathbb{P}_U^n(A) - \mathbb{P}_{\bar{\mathcal{S}}}^n(A)|$ is bounded by

$$(A.18) \qquad \begin{aligned} &|\mathbb{P}_U^n(A \setminus \mathcal{Q}_n^{\bullet\star})| + |\mathbb{P}_{\bar{\mathcal{S}}}^n(A \setminus \mathcal{Q}_n^{\bullet\star})| \\ &\quad + |\mathbb{P}_U^n(A \cap \mathcal{Q}_n^{\bullet\star}) - \mathbb{P}_{\bar{\mathcal{S}}}^n(A \cap \mathcal{Q}_n^{\bullet\star})| \\ &\leq |\mathbb{P}_U^n(\mathcal{Q}_n^\bullet \setminus \mathcal{Q}_n^{\bullet\star})| + |\mathbb{P}_{\bar{\mathcal{S}}}^n(\mathcal{Q}_n^\bullet \setminus \mathcal{Q}_n^{\bullet\star})| \\ &\quad + |\mathbb{P}_U^n(A \cap \mathcal{Q}_n^{\bullet\star}) - \mathbb{P}_{\bar{\mathcal{S}}}^n(A \cap \mathcal{Q}_n^{\bullet\star})|. \end{aligned}$$

Thanks to (A.17) and (A.15), the two first terms in the right-hand side of (A.19) go to 0. It remains to show that the last one is uniformly negligible:

$$|\mathbb{P}_U^n(A \cap \mathcal{Q}_n^{\bullet\star}) - \mathbb{P}_{\bar{\mathcal{S}}}^n(A \cap \mathcal{Q}_n^{\bullet\star})| \leq \#(A \cap \mathcal{Q}_n^{\bullet\star})\left|\frac{1}{\#\mathcal{Q}_n^\bullet} - \frac{2n}{C_n 3^n}\right|$$

$$\leq \#\mathcal{Q}_n^\bullet\left|\frac{1}{\#\mathcal{Q}_n^\bullet} - \frac{2n}{C_n 3^n}\right| = \left|1 - \frac{2n\#\mathcal{Q}_n^\bullet}{C_n 3^n}\right|.$$

Thanks to (A.16), this last term goes to 0.

**A.6. Proof of Proposition 3.2.** We first give two lemmas.

LEMMA A.3. *Set $(f, \zeta) \in \mathbb{T}$; there exists a sequence of functions $(f_n, \zeta_n) \in \mathbb{T}_n$ such that $\lim \|\zeta_n - \zeta\|_\infty = 0$.*

PROOF. By a density argument, we may take $\zeta$ is $C^1$. The function $\zeta_n$ must be nonnegative, piecewise linear on the intervals $[i/(2n), (i+1)/(2n)]$ and must satisfy $\zeta_n((i+1)/(2n)) - \zeta_n(i/(2n)) = \pm 1/\sqrt{n}$ (for $i \in [\![0, 2n-1]\!]$). A construction of this approximating sequence is achieved as follows: we write $x_j$ for $j/(2n)$ and we just define $(\zeta_n(x_j))$ since the other values are defined by linear interpolation; we set $\zeta_n(0) = 0$, and, for $0 \leq i \leq 2n - 1$,

$$\zeta_n(x_{i+1}) = \begin{cases} \zeta_n(x_i) + 1/\sqrt{n}, & \text{if } \zeta_n(x_i) \leq \zeta(x_i), \\ \zeta_n(x_i) - 1/\sqrt{n}, & \text{if } \zeta_n(x_i) > \zeta(x_i). \end{cases}$$

Set $b = \sup\{|\zeta'(x)|, x \in [0,1]\}$; one can iteratively establish the following formula, valid for all $i \in [\![0, 2n]\!]$ and for $n$ large enough:

$$|\zeta(x_i) - \zeta_n(x_i)| \leq 1/\sqrt{n} + b/(2n). \qquad \square$$

LEMMA A.4. *For $x$ and $y$ in $[0,1]$, we set $d_\zeta(x, y) = \zeta(x) + \zeta(y) - 2\check{\zeta}(x, y)$; then for any $\theta \in \mathcal{O}$,*

$$d_{\zeta^{(\theta)}}(x, y) = d_\zeta(x \oplus \theta, y \oplus \theta).$$

PROOF. We set $\theta_n = \lfloor 2n\theta \rfloor / 2n$ and let $(\zeta_n)$ be the sequence given by Lemma A.3; we have

$$\begin{aligned}|d_{\zeta^{(\theta)}}(x, y) &- d_\zeta(x \oplus \theta, y \oplus \theta)| \\ &\leq |d_{\zeta^{(\theta)}}(x, y) - d_{\zeta^{(\theta_n)}}(x, y)| \\ &\quad + |d_{\zeta^{(\theta_n)}}(x, y) - d_{\zeta_n^{(\theta_n)}}(x, y)| \\ &\quad + |d_{\zeta_n^{(\theta_n)}}(x, y) - d_{\zeta_n}(x \oplus \theta_n, y \oplus \theta_n)| \\ &\quad + |d_{\zeta_n}(x \oplus \theta_n, y \oplus \theta_n) - d_{\zeta_n}(x \oplus \theta, y \oplus \theta)| \\ &\quad + |d_{\zeta_n}(x \oplus \theta, y \oplus \theta) - d_\zeta(x \oplus \theta, y \oplus \theta)|.\end{aligned}$$

Each term in the right-hand side goes to 0 (the first one since $\zeta$ is uniformly continuous, the second one and the fifth one because $\zeta_n \longrightarrow \zeta$ uniformly, the fourth one is bounded by $4/\sqrt{n}$ and the third one is 0 by the properties of the rerooting operator $J^{(\theta_n)}$ on $\Omega_n$). $\square$

PROOF OF PROPOSITION 3.2. (i) Let $(f, \zeta) \in \mathbb{T}$, $\theta \in \mathcal{O}$ and $(f^{[\theta]}, \zeta^{(\theta)}) = J^{(\theta)}(f, \zeta)$; we just have to show that $(f^{[\theta]}, \zeta^{(\theta)})$ has the snake property:

$$\zeta^{(\theta)}(x) = \zeta^{(\theta)}(y) = \check{\zeta}^{(\theta)}(x, y) \Rightarrow d_{\zeta^{(\theta)}}(x, y) = 0$$

$$\Rightarrow d_\zeta(x \oplus \theta, y \oplus \theta) = 0$$
$$\Rightarrow f(x \oplus \theta) = f(y \oplus \theta)$$
$$\Rightarrow f^{[\theta]}(x) = f^{[\theta]}(y).$$

(ii) We use the sequence $\zeta_n$ defined in Lemma A.3 and converging uniformly to $\zeta$. Let $\theta$ and $\theta' \in \mathcal{O}$ and set $\theta_n = \lfloor 2n\theta \rfloor/(2n)$ and $\theta'_n = \lfloor 2n\theta' \rfloor/(2n)$ be approximating sequences belonging to $\mathcal{O}_n$. By simple properties of the discrete rerooting operator,

$$(A.19) \qquad (\zeta_n^{(\theta_n)})^{(\theta'_n)} = \zeta_n^{(\theta_n \oplus \theta'_n)}.$$

Let us prove that

$$(A.20) \qquad (\zeta^{(\theta)})^{(\theta')} = \zeta^{(\theta \oplus \theta')}.$$

On the first hand,

$$|\zeta^{(\theta \oplus \theta')} - \zeta_n^{(\theta_n \oplus \theta'_n)}| \leq |\zeta^{(\theta \oplus \theta')} - \zeta_n^{(\theta \oplus \theta')}| + |\zeta_n^{(\theta \oplus \theta')} - \zeta_n^{(\theta_n \oplus \theta'_n)}|.$$

The first term in the right-hand side goes to 0 because $\zeta_n \longrightarrow \zeta$ uniformly; since $|\theta_n \oplus \theta'_n - \theta \oplus \theta'| \leq 1/n$, the second one is smaller than $8/\sqrt{n}$.
On the other hand,

$$|(\zeta^{(\theta)})^{(\theta')} - (\zeta_n^{(\theta_n)})^{(\theta'_n)}| \leq |(\zeta^{(\theta)})^{(\theta')} - (\zeta_n^{(\theta)})^{(\theta')}| + |(\zeta_n^{(\theta)})^{(\theta')} - (\zeta_n^{(\theta_n)})^{(\theta')}|$$
$$+ |(\zeta_n^{(\theta_n)})^{(\theta')} - (\zeta_n^{(\theta_n)})^{(\theta'_n)}|.$$

Each of the terms in the right-hand side goes to 0 when $n$ goes to $+\infty$. This shows (A.20) and thus the action of $\mathcal{O}$ is a group action on the second coordinate of elements of $\mathbb{T}$. Since clearly $(f^{[\theta]})^{[\theta']} = f^{[\theta \oplus \theta']}$, (ii) is proved. □

**A.7. Proof of Lemma 4.12.** (i) is clear (since $[0,1]^2$ is compact).
(ii) For any $s_1, s_2, s_3 \in \mathcal{O}$, one has

$$(A.21) \quad \delta(\bar{x}, \bar{z}) \leq d_\mathbb{T}(x^{(s_1)}, z^{(s_3)}) \leq d_\mathbb{T}(x^{(s_1)}, y^{(s_2)}) + d_\mathbb{T}(y^{(s_2)}, z^{(s_3)}).$$

By continuity and compactness, there exist $\theta_1, \theta_2, \theta_3$ and $\theta_4 \in \mathcal{O}$ such that

$$\delta(\bar{x}, \bar{y}) = d_\mathbb{T}(x^{(\theta_1)}, y^{(\theta_2)}), \qquad \delta(\bar{y}, \bar{z}) = d_\mathbb{T}(y^{(\theta_3)}, z^{(\theta_4)}).$$

Applying (A.21) with $s_1 = \theta_1, s_2 = \theta_2$ and $s_3 = \theta_4 \oplus \theta_2 \oplus (-\theta_3)$ and (4.12),

$$\delta(\bar{x}, \bar{z}) \leq \delta(\bar{x}, \bar{y}) + d_\mathbb{T}(y^{(\theta_2)}, z^{(\theta_4 \oplus \theta_2 \oplus (-\theta_3))})$$
$$\leq \delta(\bar{x}, \bar{y}) + d_\mathbb{T}(y^{(\theta_3 \oplus \theta_2 \oplus (-\theta_3))}, z^{(\theta_4 \oplus \theta_2 \oplus (-\theta_3))}) \leq \delta(\bar{x}, \bar{y}) + 4\delta(\bar{y}, \bar{z}).$$

(iii) First, let us establish that $B_\delta(\bar{x}, \rho)$ is open in the quotient topology. Set $U = \pi^{-1}(B_\delta(\bar{x}, \rho))$; we have to show that $U$ is open in $\mathbb{T}$. Consider $y \in U$; one has $\bar{y} \in B_\delta(\bar{x}, \rho)$. There exist $\theta_1, \theta_2$ such that

$$\delta(\bar{x}, \bar{y}) = d_\mathbb{T}(x^{(\theta_1)}, y^{(\theta_2)}) = \lambda < \rho.$$

Set $\varepsilon = \rho - \lambda$. Let $z \in B(y, \varepsilon/4)$ in $\mathbb{T}$. We have

$$\delta(\bar{x}, \bar{z}) \leq \delta(\bar{x}, \bar{y}) + 4\delta(\bar{y}, \bar{z}) \leq \lambda + 4d(y, z) < \lambda + \varepsilon = \rho,$$

then $z \in U$, and $B(y, \varepsilon/4) \subset U$ and $U$ is open.

Consider an open set $V$ in $\mathbb{T}/\mathcal{O}$ and $\bar{x} \in V$. The set $\pi^{-1}(V)$ is open in $\mathbb{T}$ and, thus, there exists $\rho > 0$ such that $B(x, \rho) \subset \pi^{-1}(V)$. We now prove that $B_\delta(\bar{x}, \rho/4) \subset V$. To this end, set $\bar{z} \in B_\delta(\bar{x}, \rho/4)$. There exist $\theta_1, \theta_2$ such that $d_\mathbb{T}(z^{(\theta_1)}, x^{(\theta_2)}) < \rho/4$. By (4.12), one has $d_\mathbb{T}(z^{(\theta_1 - \theta_2)}, x) < \rho$ and so $z^{(\theta_1 - \theta_2)} \in B(x, \rho)$. This implies that $\bar{z} = \bar{z}^{(\theta_1 - \theta_2)} \in V$. This says that $B_\delta(\bar{x}, \rho/4) \subset V$.

(iv) For any $p$ and any sequence $\bar{z}_0, \ldots, \bar{z}_p \in \mathbb{T}/\mathcal{O}$ such that $\bar{z}_0 = \bar{x}$, and $\bar{z}_p = \bar{y}$,

$$\sum_{i=0}^{p-1} \delta(\bar{z}_i, \bar{z}_{i+1}) = \sum_{i=0}^{p-1} d_\mathbb{T}(z_i^{(\theta_i)}, z_{i+1}^{(\theta'_i)}),$$

where the $\theta_i$ and $\theta'_i$ are elements of $\mathcal{O}$ that reach $\delta(\bar{z}_i, \bar{z}_{i+1})$. By (4.12), for any $\theta, s, s'$,

$$d_\mathbb{T}(z_i^{(s)}, z_{i+1}^{(s')}) \geq \tfrac{1}{4} d_\mathbb{T}(z_i^{(\theta)}, z_{i+1}^{(s' \oplus \theta \oplus (-s))});$$

applying this inequality in a convenient way and successively from $i = 1$ to $i = p - 1$, we find a sequence $\widetilde{\theta}_1, \ldots, \widetilde{\theta}_p$ such that $\widetilde{\theta}_1 = \theta'_0$, and such that

$$\sum_{i=0}^{p-1} \delta(\bar{z}_i, \bar{z}_{i+1}) \geq d_\mathbb{T}(z_0^{(\theta_0)}, z_1^{(\theta'_0)}) + \tfrac{1}{4} \sum_{i=1}^{p-1} d_\mathbb{T}(z_i^{(\widetilde{\theta}_i)}, z_{i+1}^{(\widetilde{\theta}_{i+1})}).$$

Thus,

$$\sum_{i=0}^{p-1} \delta(\bar{z}_i, \bar{z}_{i+1}) \geq \tfrac{1}{4} d_\mathbb{T}(z_0^{(\theta_0)}, z_p^{(\widetilde{\theta}_p)}) \geq \tfrac{1}{4} \delta(\bar{x}, \bar{y}).$$

**Acknowledgments.** We thank Mireille Bousquet-Mélou for many helpful discussions and references. We are grateful to Thomas Duquesne for his proof of Lemma 4.15. We are grateful to the referee for his careful reading, his suggestions and comments that helped us to really improve our manuscript.

UNIVERSITÉ DE VERSAILLES
45 AVENUE DES ÉTATS-UNIS
78035 VERSAILLES CEDEX
FRANCE
E-MAIL: marckert@math.uvsq.fr
mokkadem@math.uvsq.fr